\documentclass[12pt]{amsart}
\usepackage[english]{babel}
\usepackage{CJK,multicol,amsmath,amssymb,mathrsfs,epsfig,cite,graphicx}
\usepackage{color}
\usepackage{amsthm}
\usepackage{extarrows}
\usepackage{verbatim}
\usepackage{lineno}
\usepackage{appendix}
\usepackage{hyperref}
\modulolinenumbers[5]

\title[]{Harmonic  and Monogenic Representations of White Noise and its weighted harmonic Bergman estimates }

\author[Fang Yang]{Fang Yang}
\address{Fang Yang, Department of Physics\\
	Zhejiang Normal University\\
	China.}
\email{fyang1993math@163.com}

\author[Jiecheng Chen]{Jiecheng Chen}
\address{Jiecheng Chen, Department of Mathematics\\
	Zhejiang Normal University\\
	China.}
\email{jcchen@zjnu.edu.cn}

\author[Tao Qian]{Tao Qian}
\address{Tao QIAN, Macau Center for Mathematical Sciences\\
	Macau University of Science and Technology\\
	China.}
\email{tqian@must.edu.mo}
\thanks{*Corresponding author. Email: tqian@must.edu.mo}

\begin{document}
\begin{sloppypar}
	
	\newtheorem{Def}{Definition}
	\newtheorem{Thm}{Theorem}
	\newtheorem{lem}{Lemma}
	\newtheorem{prop}{Proposition}
	\newtheorem{cor}{Corollary}
	\newtheorem{rmk}{Remark}
	\newtheorem{Apdx}{Appendix}
	\newtheorem{exmp}{Example}
\newcommand{\R}{\mathbb{R}}
\newcommand{\C}{\mathbb{C}}
\newcommand{\undl }{\underline }

\maketitle
\begin{CJK*}{GBK}{song}

\section*{Abstract}
Continuing Feichtinger's generalized stochastic process theory, we analyze harmonic and monogenic representations of white noise in Euclidean spaces. We define its convolution with square-integrable functions. Using the weighted Sobolev regularity of the white noise, we establish a high-dimensional  Sokhotski-Plemelj formula for the white noise density. Furthermore, we embed these representations into new weighted Bergman spaces, proving that the weighted Bergman norms of the harmonic representations are bounded by the white noise's weighted Sobolev norms.

\bigskip

\noindent MSC 2020 Classification:  41A99; 60H15

\bigskip

\noindent {\em Keywords}:
Monogenic representation, Segal algebra,  white noise, weighted Sobolev space, {Sokhotski-Plemelj} formula, weighted Bergman space

\section{Introduction}
In this paper, we focus on harmonic/monogenic ``lifting up'' of the white noise on $\mathbb{R}^d$. Here,  monogenic functions are the hypercomplex analytic analogues of holomorphic functions of one complex variable \cite{Delanghe1992-Clifford, McIntosh1996}.  This connects to the Dirichlet problem for the Laplace equation in the upper-half space:
\begin{equation}\label{white-noise-HalfPlane}
	\left\{
	\begin{aligned}
		\left(\partial_y^2+\triangle_x\right)  u(y,x,\omega)&=0, ~(y,x) \in \mathbb{R}^{1+d}_+,\\
		u(0,x, \omega)&=S(x, \omega),~x \in \mathbb{R}^d,~\mathrm{a.s.}~\omega \in \Omega,
	\end{aligned}
	\right.
\end{equation}
with the white noise as the boundary value $S(x, \omega)$. 
In physics and engineering, the Gaussian white noise is often defined via a random orthogonal series.
Let $\mathcal{H}$ denote a separable Hilbert space of functions equipped with a complete orthonormal basis $\{\vec{e}_j\}_{j=1}^\infty$.  The
Gaussian white noise is then defined as the formal series \cite{Khader2010-Noise,Veraar2010-Regularity-GaussianNoise,Lord_2014}
\begin{equation}\label{formal-def-noise}
S = \sum_{j=1}^{\infty} \gamma_j \vec{e}_j,
\end{equation}
where the random variables $\{\gamma_j\}_{j=1}^\infty$ are independent and identically distributed according to the standard normal distribution $\text{Norm}(0,1)$.  The series \eqref{formal-def-noise} does not converge in the classical sense \cite{Lord_2014}.
However, it converges  in the distribution sense under suitable conditions, for example when $\mathcal{H}= L^2[0,1]$ \cite{Veraar2010-Regularity-GaussianNoise}.  

Returning to the rigorous definition, a  white noise process on a locally compact Abelian (LCA) group is defined as a generalized stochastic processes (GSP) developed by H. G. Feichtinger \cite{Feichtg2014-whiteNoise}.   His theory  starts from the Segal algebra $ \mathcal{S}_0(G) $ \cite{Feichtinger1981}. The latter, by definition, is a subspace of the Fourier algebra $A(G)$, being defined as the image of $L^1(\hat{G})$ under the (inverse) Fourier transform. The norm of $f\in A(G)$  is defined as  $\|f\|_{A(G)}=\|h\|_{L^1}$ for $\hat{f}=h$.
\begin{Def}\cite{Feichtg2014-whiteNoise,Jakobsen2018-FeichtingerAlgebra}
	Let a nonzero window function $\mathfrak{n}
	\in A(G) \cap \mathcal{C}_c(G)$, i.e., $\mathfrak{n}$ is continuous, with compact support and integrable Fourier transform $\widehat{\mathfrak{n}}$. Then the Segal  algebra $\mathcal{S}_0(G)$ is defined by
	\begin{equation*}
	\mathcal{S}_0(G):=\left\{f \in A(G) :~\|f\|_{S_0}:=\int_G\left\|\Rsh_{-y}\mathfrak{n} \cdot f\right\|_A d y<\infty\right\},
	\end{equation*}
	where $\Rsh_y$ is the translation operator with $\Rsh_y\mathfrak{n}(x)=\mathfrak{n}(x+y),~y\in G$.
\end{Def}
Different window functions  define the same Segal  algebra with equivalent norms. It is the the smallest Banach space among all Banach spaces which are isometrically invariant under translation and modulation  (multiplication by characters) and is contained in integrable function space $ L^1(G)$ with compactly supported Fourier transform. Therefore, the Segal algebra, serving as a test function space, can be regarded as an alternative to the Schwartz space. Moreover, the former constitutes a Banach space, whereas the latter is  a  topological vector space. The dual space $ \mathcal{S}^{\prime}_0(G) $ is treated as the distribution space on $G$. Then the  GSP on the LCA group $G$ is defined as the  Hilbert space valued distribution on the test function space $ \mathcal{S}_0(G) $.
Within the framework of GSP, Feichtinger ingeniously defined   (wide-sense)  time stationary  and  (wide-sense) frequency stationary  through the translation operator and modulation  operator on LCA group, respectively. The white noise is a GSP that is  both time stationary and frequency stationary.  Nevertheless, the white noise admits an equivalent characterization: it is a scalar multiple of an isometry from $L^2(G)$ into a Hilbert space  $H$ \cite{Feichtg2014-whiteNoise}. This definition of the  white noise via its isometric property also appears in \cite{Chauleur2025,Killip2019-whiteNoise,Mackowiak2025-Noise}. The white noise is a  Hilbert space valued distribution on  the Segal algebra  $ \mathcal{S}_0(G) $.  When the target Hilbert space $H$ of the white noise $S:~L^2(G)\rightarrow H$ is taken to be $ L^2(\Omega, \mathcal{F}, \mathbb{P})$, we can define its sample paths as in \eqref{Def-path-GSP}. Specifically, for each  $\omega \in \Omega$, the sample path is denoted by  $S_\omega(\cdot)=S(\cdot,\omega)$.  When the white noise is defined in the unit circle, its sample paths almost surely belong to the Sobolev space $  H^{-s}(\partial \mathbb{D}) $  for all $s>\frac{1}{2}$ \cite{Veraar2010-Regularity-GaussianNoise}.

The harmonic/monogenic "lift up" of the white noise can be regarded as a method of  regularization.  This  "lift up"  originates from the analytic representation theory of the distributions  \cite{Bremermann1965}. For an analytic function, its behavior at all singular points is completely determined by its values at regular points. It is natural to lift distributions on the real line to analytic functions, and to investigate such distributions through the limits of these analytic functions. Indeed,  the following theorem holds:
\begin{Thm}\cite{Bremermann1965}
	Let $f$ be a Schwartz distribution on $\mathbb{R}$. Then there exists a function $F$, analytic in $\mathbb{C}\setminus \mathbb{R}$, such that
	\[
	\lim_{\epsilon \to 0^+} \int_{-\infty}^{+\infty}
	\left[ F(x + i\epsilon) - F(x - i\epsilon) \right] \varphi(x) dx
	= \langle f, \varphi \rangle
	\]	
	for any test function $\varphi$ in the Schwartz space. We  call $F$  an analytic representation of the  distribution $f$. Furthermore, the difference $ F(x + i\epsilon) - F(x - i\epsilon) $ is called the analytic jump of $f$ with jump height $\epsilon$.	
\end{Thm}
For instance, the Cauchy kernel $\{z\mapsto\frac{1}{i 2\pi}\frac{1}{z}\}$ is an analytic representation of the Dirac distribution $\delta$. We focus on constructing the monogenic/harmonic    representation of the white noise on $\mathbb{R}^d$. In practice, the regularization of the white noise can be achieved via smooth positive kernels (of the convolution type). The approach has been applied to solving PDEs with the  white noise boundary values \cite{Chauleur2025}. Additionally, incorporating the regularized noise into certain monogenic function spaces is beneficial: the reproducing kernel structure of these spaces facilitates the study of properties (and even sparse approximations) of solutions to some stochastic  PDEs \cite{Yang2021-SPOAFD-SPDE,Yang2024-whiteNoise}.

We now recall the existing work on the analytic/harmonic representations of white noise in the unit disk setting. Let    $p^\mathbb{D}$ denote the Poisson kernel on the unit disk $ \mathbb{D} $:
\begin{equation}\label{Poiss-kernel-Disk}
p^\mathbb{D}( e^{i\theta},z)=
\frac{1-|z|^2}{|e^{i\theta}-z|^2}
=
\sum_{k=-\infty}^{\infty}z^{|k|} e^{ik\theta}
,
~z \in \mathbb{D}.
\end{equation}
It is proved that the random field obtained by applying white noise to the Poisson kernel, namely,
\begin{equation}\label{solution-unitDisk}
(z,\omega)\mapsto
\left\langle  S_{\omega},~ p^\mathbb{D}(\cdot,z)
\right\rangle
\end{equation}
is harmonic in the unit disk and converges  radially weakly to the white noise itself in the distribution sense \cite{Khader2010-Noise}.  Determining the specific function spaces for these representations of the white noise is also crucial.  To this end, we first introduce the related Bochner-type spaces.
\begin{Def}
	\cite{Lord_2014,Qian2022-SPOAFD}
	\label{Def-BochnerSpc}	Consider a probability space $(\Omega,\mathcal{F}, \mathbb{P})$. Let  $\mathcal{H}(\mathfrak{D})$
	denote the Hilbert space of functions on a domain $\mathfrak{D}$. The Bochner space $ L^2(\Omega, \mathcal{H}(\mathfrak{D}) ) $ consists of all $\mathcal{H}(\mathfrak{D})$-valued random variables $ f: \mathfrak{D}\times\Omega\rightarrow\mathbb{R} $ such that $f(\cdot, \omega)\in \mathcal{H}(\mathfrak{D})$  for a.s. $\omega\in \Omega$, and
	\begin{equation*}
	\left\|f\right\|_{L^2(\Omega, \mathcal{H}(\mathfrak{D}) )}^2
	:=
	\int_{\Omega}
	\left\|f(\cdot, \omega)
	\right\|_{\mathcal{H}(\mathfrak{D})}^2 d\mathbb{P}(\omega)
	<\infty.
	\end{equation*}
\end{Def}
The random field \eqref{solution-unitDisk} belongs to the  harmonic weighted Bochner-Bergman space $ L^2(\Omega, b^2_\alpha(\mathbb{B}))$ for each  $\alpha>0$ \cite{Yang2024-whiteNoise}, where  $b^2_\alpha(\mathbb{B})$ denotes the weighted harmonic Bergman space \cite{Axler_Harmonic-funct-thrm},
\begin{equation*}\label{}
b^2_\alpha(\mathbb{B})
=
\{u\in\text{Harm}(\mathbb{B}):~\int_{\mathbb{B} }|u(re^{i\theta})|^2 (1-r^2)^\alpha rdr d\theta<\infty \}.
\end{equation*}
Similarly, consider the analytic completion defined by
\[
(z,\omega) \mapsto
\left\langle S_{\omega},\, p^{\mathbb{D}}(\cdot,z) \right\rangle
+ i\, \mathscr{H}\!\left( \left\langle S_{\omega},\, p^{\mathbb{D}}(\cdot,z) \right\rangle \right),
\]
where $\mathscr{H}\!\left( \left\langle S_{\omega},\, p^{\mathbb{D}}(\cdot,z) \right\rangle \right)$ denotes the harmonic conjugate of $\left\langle S_{\omega},\, p^{\mathbb{D}}(\cdot,z) \right\rangle$.  This random field belongs to the  analytic weighted Bochner-Bergman space $ L^2(\Omega, B^2_\alpha(\mathbb{D}) $ for $\alpha>0$  \cite{Yang2024-whiteNoise}. This analytic weighted Bergman  space is referred to  \cite{Zhu2005_Hol-funct-spc}
\begin{equation}\label{Def-WeightedBergman-disk}
B^2_\alpha(\mathbb{D})
=
\{u\in\text{Hol}(\mathbb{D}):~\int_{\mathbb{D} }|u(z)|^2 (1-|z|^2)^\alpha dzd\bar{z}<\infty \}.
\end{equation}
The white noise on the unit circle was first studied in  \cite{Yang2024-whiteNoise} in the unit disk in  relation to the  weighted analytic and harmonic  Bergman spaces. Furthermore, the heat equation driven by Gaussian white noise on the unit circle  inspired the study of a class of parabolic Bergman spaces \cite{Yang2024-whiteNoise}.

While the unit disk provides a compact setting, it is natural to extend the investigation of these representations and Bergman-type spaces to unbounded domains.  For simplicity, we occasionally adopt the shorthand notation $y+\underline{x}$ to denote elements $(y,x)$ of the upper half-spaces. This  notation $y+\underline{x}$ is borrowed from Clifford algebra  in which $y>0$, and $\underline{x}=\sum_{j=1}^{d}x_j\textbf{e}_j$ with the basis vector $\{\textbf{e}_j\}_{j=1}^{d}$
spanning the $ d $-dimensional  Euclidean space $\mathbb{R}^d$. Denote by
  \[
   P_y(\underline{x}-\underline{t})=
 P_{y+\underline{x}}(\underline{t})=
 P_{y+\underline{t}}(\underline{x})=
c_d
 \frac{y}{(y^2+|\underline{x}-\underline{t}|^2)^{\frac{1+d}{2}}},\quad
 c_d=\frac{\Gamma\left(\frac{d+1}{2}\right)}{\pi^{\frac{d+1}{2} }}
   \]
the Poisson kernel in  $\mathbb{R}^{1+d}_+$. For almost every $\omega \in \Omega$, the sample path $S_\omega$ is a tempered distribution, as in Proposition~\ref{Prop-Sobolev-rglrt-Noise}. Consider the formal solution of Problem \eqref{white-noise-HalfPlane}, namely the following dual pairing,
\begin{equation}\label{Poisson-WhiteNoise}
u_\omega(y+\underline{x}):=
\left\langle S_\omega,
P_{y+\underline{x}}\right\rangle.
\end{equation}
Its harmonic property in $ \mathbb{R}^{1+d}_+$ for each $\omega$ fixed, as well as whether it can converge vertically to the white noise itself in the distribution sense, is of particular interest. The significance of this  paper is to provide the monogenic and harmonic representations  for the  white noise on $\R^d$. It  is an analogy to the case of the unit circle in  \cite{Yang2024-whiteNoise}. The contributions of this paper include:
\begin{itemize}
\item
We study  sample paths of the white noise {as distributions} in weighted Sobolev spaces, establishing the embedding relationship between these weighted Sobolev spaces and the Segal-Feichtinger algebra. The convolution of the white noise with $ L^2(\mathbb{R}^d) $ function is defined, and the stochastic weighted Sobolev norm of the white noise is computed-expressible in terms of the Beta function.
\item
In the pathwise setting, we  propose   harmonic  and  monogenic functions representations of the white noise in   $\mathbb{R}^{1+d}_+$. This extends the classical monogenic and harmonic representations of the white noise in the unit disk  \cite{Yang2024-whiteNoise} to this unbounded domain. In the unbounded domains it is necessary to work with suitable weighted spaces.
In particular, we prove that the  harmonic  representation of the white noise  converges weakly to the boundary white noise.
\item
We embed the harmonic/monogenic representation into a weighted harmonic/monogenic Bergman space equipped with a weight function derived from a conformal transformation.  We derive new norm estimates that relate the weighted Bergman norms in  $\mathbb{R}^{1+d}_+$ to the weighted Sobolev norms of the boundary white noise pathwise.
\end{itemize}
The paper is organized as follows. In \S\ref{preliminary}, we review the necessary background and the notation  used in the paper. Some new properties of the white noise are derived in this section.
\S\ref{Sec-Sokhotski-Plemelj formula} discusses the harmonic and monogenic representations of the white noise  and their vertical weak convergence. \S\ref{Sec_A new Bergman spac} constructs the corresponding weighted Bergman spaces. Finally, \S\ref{Sec_Priori estimates} establishes norm estimates linking the weighted Bergman norms in the upper half-spaces and the weighted Sobolev norms of the boundary white noise pathwise.

\section*{Notation}\label{Notation}
We denote by $\mathscr{S}(\mathbb{R}^d)$ the space of smooth vertical decaying functions  and by $ \mathscr{S}^{\prime}(\mathbb{R}^d) $ its dual, the space of the tempered distributions. We employ the following standard notation:
\[ 
\mathbb{N} = \{1, 2, 3, \dots\}, \quad \mathbb{N}_0 = \{0, 1, 2, \dots\} = \{0\} \cup \mathbb{N}.
\]
We  let $d \in \mathbb{N}$ be the dimension of the {configuration space}. For any {multi-index} $\alpha = (\alpha_1, \dots, \alpha_d) \in \mathbb{N}_0^d$ we define its {length} and {factorial} as
\[ 
|\alpha| = \alpha_1 + \dots + \alpha_d, \quad \alpha! = (\alpha_1!) \cdots (\alpha_d!),
\]
respectively. In addition, for any $\alpha, \beta \in \mathbb{N}_0^d$ we let
\[ 
\alpha \leq \beta \stackrel{~}{\iff} \alpha_j \leq \beta_j, \quad \forall~ j = 1, \dots, d,
\]
and define the {binomial coefficient} as
\[ 
\binom{\alpha}{\beta} = \frac{\alpha!}{\beta!(\alpha - \beta)!} \quad \text{if } 0 \leq \beta \leq \alpha, \quad \binom{\alpha}{\beta} = 0 \quad \text{otherwise},
\]
where $\alpha - \beta = (\alpha_1 - \beta_1, \dots, \alpha_d - \beta_d)$. For any $x = (x_1, \dots, x_d) \in \mathbb{R}^d$ and $\alpha = (\alpha_1, \dots, \alpha_d) \in \mathbb{N}_0^d$ we write
\[ 
x^\alpha = x_1^{\alpha_1} \cdots x_d^{\alpha_d}, \quad \partial^\alpha = (\partial_1)^{\alpha_1} \cdots (\partial_d)^{\alpha_d}, \quad \partial_j = \partial_{x_j} = \frac{\partial}{\partial x_j}.
\]
We recall the  generalized product rule, commonly known as the  Leibniz rule: for any $\alpha \in \mathbb{N}_0^d$ and $f, g \in C^{|\alpha|}(\mathbb{R}^d)$,
\begin{equation}\label{LeibnizRule}
	\partial^\alpha (fg) =
	\sum_{\beta \in \mathbb{N}_0^d:~|\beta|\leq|\alpha|} \binom{\alpha}{\beta} (\partial^{\alpha-\beta} f) (\partial^\beta g).
\end{equation}
\section{Preliminaries}\label{preliminary}
\subsection{Generalized stochastic process}
H. G. Feichtinger defined generalized stochastic processes, wide-sense stationary processes, and the white noise over LCA groups \cite{Feichtg2014-whiteNoise}. The definition of generalized stochastic processes over LCA groups is based on Segal algebras which are studied in harmonic analysis on groups \cite{Feichtinger1981,Feichtg2014-whiteNoise}.
\subsubsection{Segal-Feichtinger algebra}
For the LCA group $G=\mathbb{R}^d$, a more convenient (equivalent) characterization of the Segal algebra norm is via the short-time Fourier transform:
\begin{Def}\cite{Feichtg2014-whiteNoise}
Let $f, g \in L^2\left(\mathbb{R}^d\right)$. The short-time Fourier transform (STFT) $\mathcal{V}_g(f): \mathbb{R}^d \times \mathbb{R}^d \longrightarrow \mathbb{C}$ of $f$ with window $g$ is defined by
\[
 \left(\mathcal{V}_g f\right)(x, \xi):=\int_{\mathbb{R}^d} f(t) \overline{g(t-x)} e^{-2 \pi i \xi \cdot t} d t.
  \]
\end{Def}
Using the STFT, the Segal-Feichtinger algebra can be defined as follows.
\begin{Def} \cite{Feichtg2014-whiteNoise}
	Let the window be the Gaussian function, i.e.,  $\mathfrak{n}(x):=e^{-\pi|x|^2}$. Segal algebra $\mathcal{S}_0
 \left(\mathbb{R}^d\right)$ is then given by
\[
 \mathcal{S}_0
 \left(\mathbb{R}^d\right):=\left\{f \in L^2\left(\mathbb{R}^d\right) :~\|f\|_{\mathcal{S}_0
 	}=\left\|\mathcal{V}_{\mathfrak{n}} f\right\|_{L^1\left(\mathbb{R}^{2 d}\right)}<\infty\right\} .
  \]
\end{Def}
The Banach space $\mathcal{S}_0
(\mathbb{R}^d)$ can be characterized as the smallest Banach space among all the Banach spaces that contains all $L^1$ functions  with compactly supported Fourier transform and is isometrically invariant under translation $\Rsh_{x}:~f(y)\mapsto f(y+x)$ and modulation $\circlearrowleft_x:~f(y)\mapsto e^{i2\pi y\cdot x}f(y)$.

There is an equivalent norm on the Segal algebra $\mathcal{S}_0
(\mathbb{R}^d)$.   Let \(\varrho \in {\mathscr {S}}(\mathbb {R}^d)\) be a smooth, symmetric bump function, that is \(0 \le \varrho \le 1\), \(\varrho (\xi ) = 1\) if \(|\xi |_\infty\le 1/2\), \(\varrho (\xi ) = 0\) if \(|\xi |_\infty \ge 1\). Let
\[ \varrho _k(\xi ) = \varrho (\xi -k), \quad k \in {\mathbb {Z}^d}. \]
Define \(Q_0 = [-1/2,1/2)^d\) and \(Q_k = k + Q_0\). Define
\[
\sigma _k(\xi ) = \varrho_k(\xi )\left( \sum _{k\in {\mathbb {Z}}}\varrho_k(\xi )\right) ^{-1}, \quad k \in {\mathbb {Z}^d}.  \]
Then, for some \(c > 0\), the smooth cut-off functions \(\sigma _k\) satisfy
\begin{equation}\label{isometric decomposition bump}
\begin{aligned} {\left\{ \begin{array}{ll} |\sigma _k(\xi )| \ge c, \quad \forall~ \xi \in Q_k,\\ {\operatorname{Spt}}\{\sigma _k\}
	\subset \{\xi:~|\xi -k|_\infty \le 1\},\\
\sum _{k\in {\mathbb {Z}}} \sigma _k(\xi ) = 1, \quad \forall~ \xi \in {\mathbb {R}^d},\\
 |\partial^{\alpha }\sigma _k(\xi )| \le C_m, \quad \forall~ \xi \in {\mathbb {R}^d}, |\alpha | \le m. \end{array}\right. } \end{aligned}
\end{equation}
Given a sequence of functions \(\sigma _k\) satisfying \eqref{isometric decomposition bump}, the sequence of operators
	\[
	\begin{aligned} \Box _k = {\mathscr {F}}^{-1}\sigma _k {\mathscr {F}}, \quad k \in {\mathbb {Z}^d},
	\end{aligned}
	\]
	is called a family of \textbf{frequency-uniform decomposition} operators.  This construction is widely used in harmonic analysis \cite{Beyni2020}. An important role  will be played by the polynomial weights
\[
\langle x\rangle^\mu: =\left(1+|x|^2\right)^{\mu / 2}
\]
with $\mu \in \mathbb{R}$. For $p \in[1, \infty],~ \mu \in \mathbb{R}$, owing to the relation
\begin{equation}\label{Def-weighted Lebesgue norm}
\|f\|_{L^p\left(\langle \cdot\rangle^\mu\right)}:=\left\|\langle \cdot\rangle^\mu f\right\|_{L^p}=\left(\int_{\mathbb{R}^d}\left|f(x)\langle x\rangle^\mu\right|^p \mathrm{~d} x\right)^{1 / p},
\end{equation}
we use the convention $L^p\left(\mathbb{R}^d,\langle x\rangle^{p\mu} dx \right)
=
L^p\left(\langle \cdot\rangle^\mu\right)$.
  If $\mu=0$ we simply write $L^p=L^p\left(\langle \cdot\rangle^0\right)$.
\begin{Def}\label{Def-wt-mdlt-spc}\cite{Guo2019}
	Given \(p,q \in [1,\infty ]\), $\alpha,\beta\in\mathbb{R}$ and \(\{\Box _k\}_{k\in\mathbb {Z}^d}\) a family of frequency-uniform decomposition operators. The (weighted) modulation space norm with respect to \(\{\Box _k\}_k\) is defined as
	\[
	\begin{aligned} \Vert f\Vert _{M^{\alpha,\beta}_{p,q}(\mathbb {R}^d)}
	=
	\big \Vert \langle{k}\rangle^\alpha
	 \Vert \Box _k f\Vert _{L^p({\mathbb {R}^d},\langle\cdot\rangle^\beta)}
	 \big \Vert _{\ell ^q_k({\mathbb {Z}^d})},
	\end{aligned}
	\]
where $\langle{k}\rangle=(1+|k|^2)^{\frac{1}{2}}$ and the weighted $L^p$ norm $L^p({\mathbb {R}^d},\langle\cdot\rangle^\beta)$ is defined   $\|f\|_{L^p\left(\langle \cdot\rangle^\beta\right)}=\left\|\langle \cdot\rangle^\beta f\right\|_{L^p}$.
\end{Def}
Among the modulation spaces, several well-known function spaces arise as special cases studied in  \cite{Cordero2020}\\
Weighted $L^2$-spaces:
\[
M^{0,\beta}_{2,2}
\left(\mathbb{R}^d\right)=L^2\left(\langle \cdot\rangle^\beta, \mathbb{R}^d\right)=\left\{f: f(x)\langle x\rangle^\beta \in L^2\left(\mathbb{R}^d\right)\right\} .
 \]
Fourier-Sobolev spaces:
\[
M^{\alpha,0}_{2,2}
\left(\mathbb{R}^d\right)=H^\alpha\left(\mathbb{R}^d\right)=\left\{f: \hat{f}(\xi)\langle\xi\rangle^\alpha \in L^2\left(\mathbb{R}^d\right)\right\} .
 \]
When $p=q=1$ and $\alpha=\beta=0$, the modulation space $M^{0,0}_{1,1}(\mathbb{R}^d)$  is  the Segal-Feichtinger algebra $\mathcal{S}_0(\mathbb{R}^d)$ \cite{Beyni2020}. Since the  Schwartz space $\mathscr{S}(\mathbb{R}^d)\subset \mathcal{S}_0(\mathbb{R}^d)$, its dual space ${\mathcal{S}_0}^\prime(\mathbb{R}^d)$ is a  Banach space in  $\mathscr{S}^\prime(\mathbb{R}^d)$. For $1\leq p\leq \infty$, there holds
\begin{equation}
\mathscr{S}(\mathbb{R}^d)\subset \mathcal{S}_0(\mathbb{R}^d)\subset
L^p(\mathbb{R}^d)\subset
{\mathcal{S}_0}^\prime(\mathbb{R}^d)
\subset\mathscr{S}^\prime(\mathbb{R}^d)
\end{equation}
and  $\mathcal{S}_0(\mathbb{R}^d)$ is dense in $(L^p(\mathbb{R}^d), \|\cdot\|_{L^p})$ for $1\leq {p}<\infty$ \cite{Feichtg2019}.
The space ${\mathcal{S}_0}(G)$ is both closed under the convolution  with $L^1(G)$ functions and pointwise multiplication by $ A(G) $ functions:
\begin{prop}\cite{Jakobsen2018-FeichtingerAlgebra}
	If $f \in{\mathcal{S}_0}(G)$, then
	\begin{enumerate}
		\item
		for any $h \in L^1(G)$ the function $f \ast h\in\mathcal{S}_0(G)$ and $\|f\ast h\|_{\mathcal{S}_0} \leq$ $\|f\|_{\mathcal{S}_0}\|h\|_{L^1}$,
		\item
		for any $h \in A(G)$ the function $f \cdot h\in\mathcal{S}_0(G)$ and $\|f \cdot h\|_{\mathcal{S}_0} \leq$ $\|f\|_{\mathcal{S}_0}\|h\|_{A}$.
	\end{enumerate}
In particular the inclusion $\mathcal{S}_0(G)\subseteq L^1(G)\cap A(G)$
implies that $\mathcal{S}_0(G)$ is a Banach algebra under pointwise multiplication and convolution.	
\end{prop}
For $\sigma \in {\mathcal{S}_0}^\prime(G)$ (mild   distributions) and $f\in  {\mathcal{S}_0}(G)$, we define frequently used operators through the $({\mathcal{S}_0}^\prime,{\mathcal{S}_0})$ dual pairing:
\begin{equation}\label{def-convolution-S_0}
\begin{aligned}
\left\langle  {\Rsh_x} \sigma, f\right\rangle&:=\left\langle\sigma, {\Rsh_{-x}} f\right\rangle, \quad \forall~ x \in G  \\
\left\langle\sigma^{\checkmark}, f\right\rangle&:=\left\langle\sigma, f^{\checkmark}\right\rangle, \\
\langle g \ast \sigma, f\rangle&:=\left\langle\sigma, g^{\checkmark} \ast f\right\rangle, \quad \forall~  g \in L^1(G)  \\
\langle h \sigma, f\rangle&:=\langle\sigma, h f\rangle, \quad \forall~ h \in A(\hat{G})\\
\langle\hat{\sigma}, f\rangle&:=\langle\sigma, \hat{f}\rangle ,\\
\end{aligned}
\end{equation}
where $ {\Rsh_{x}} f(t):=f(t+x) $ is the  translation, $   f^{\checkmark}(t):=f(-t) $ is the reflection and $\hat{f}$ is the  Fourier transform of $f$. When $ \sigma \in {\mathcal{S}_0}^\prime(\mathbb{R}^d)$ and $f\in  {\mathcal{S}_0}(\mathbb{R}^d) $, the value of their convolution at  $x\in\mathbb{R}^d$ is given by \cite{Feichtg2019}:
\[
(\sigma\ast f)(x)=\langle\sigma,  \Rsh_{x}f^{\checkmark}\rangle.
 \]
 For $\sigma\in{\mathcal{S}_0}^\prime(G)$ and $g\in L^1(G)$, the following  holds for  every test function  $ f \in\mathcal{S}_0(\hat{G}) $:
\begin{equation}\label{Fourier exchange formula}
\langle (g \ast \sigma)^\wedge, f\rangle
=
\left\langle \hat{\sigma}\hat{g},   f
\right\rangle.
\end{equation}
This gives the Fourier exchange formula $(g \ast \sigma)^\wedge=\hat{\sigma}\hat{g}=\hat{g}\hat{\sigma}$ for  $\sigma\in{\mathcal{S}_0}^\prime(G)$ and $g\in L^1(G)$.
\subsubsection{generalized stochastic process and  white noise}
\begin{Def}\cite{Feichtg2014-whiteNoise}
	Let $ H $ be a Hilbert space. A measurable mapping $X: G\rightarrow H$  is called a stochastic process on $\R^d$.
\end{Def}
\begin{Def}\cite{Feichtg2014-whiteNoise}
A bounded linear mapping
\[
\begin{aligned}
X:~ \mathcal{S}_0
(G)&\rightarrow H \\
f&\mapsto \left\langle X, f  \right\rangle
\end{aligned}
 \]
 is called a generalized stochastic process (GSP). The space of all such processes is denoted by $ GSP(G) $.
\end{Def}
We actually can think of $H$ being a class of   $\mathbb{C}$-valued random variables with zero expectation on an arbitrary probability space $(\Omega,\mathcal{F}, \mathbb{P})$, that is, $L^2_0(\Omega,\mathcal{F}, d\mathbb{P})$. For every $\omega\in \Omega$, we define a \textbf{sample path} $X_\omega$ of the GSP $X:~\mathcal{S}_0
(G)\mapsto L^2_0(\Omega, d\mathbb{P})$ via the following linear functional:
\begin{equation}\label{Def-path-GSP}
\left\langle X_\omega, f  \right\rangle
:=
\left\langle X, f  \right\rangle\vert_\omega,~\forall~f\in\mathcal{S}_0
(G).
\end{equation}
In other words, for any $X\in GSP$  as a bounded linear operator from $\mathcal{S}_0
$ to $ L^2_0(\Omega, d\mathbb{P})$, one can define a linear functional $\left\langle X_\omega, \cdot \right\rangle$ in the pathwise sense for almost every sample $\omega$.
\begin{rmk}
 For $X\in GSP(\R^d)$, for almost every $\omega\in \Omega$ the sample path $X_\omega$ induces a linear functional on
 $ \mathcal{S}_0$. This functional, however, is not necessarily continuous.
\end{rmk}
We now give the  definition  of the Fourier transform for GSP.
\begin{Def}\cite{Feichtg2014-whiteNoise}
	Given a GSP $X$ on $G$, we define a GSP $\hat{X}$ on $G$:
\[ 	
	\hat{X}(f):=X(\hat{f}),~ \forall~ f \in \mathcal{S}_0(G).
	 \]
Here $\hat{X}$ is called the Fourier transform (spectral process) to $X$. Similarly, we define a GSP $\check{X}$ on $G$:
\[
\check{X}(f):=X(\check{f}) ,~ \forall~ f \in \mathcal{S}_0(G),
 \]
where $\check{f}$ is the inverse Fourier transform for $f$.
\end{Def}
Since the mappings $f \mapsto \check{f}$ (inverse Fourier transform) and $f \mapsto \hat{f}$ (forward Fourier transform) define isomorphisms on  $\mathcal{S}_0(G)$, respectively, it is clear that $\check{X}$ and $\hat{X}$ are GSPs.
\begin{Def}\cite{Feichtg2014-whiteNoise}
For  $X\in GSP(G)$  and $f \in \mathcal{S}_0(G)$,  the following operations are defined:
\begin{enumerate}
\item
 $\langle{h X, f}\rangle:=\langle X, hf\rangle$, for $h \in A(\hat{G})$;
\item
$\langle{u \ast X, f}\rangle:= \langle{ X, u^{\checkmark}\ast f}\rangle$, for $u\in L^1(G)$.
\end{enumerate}
Since $f \mapsto h f$ and $f \mapsto u\ast f$ define linear and bounded operaters on $\mathcal{S}_0(G)$, it follows that $h X$ and $u\ast X$ define GSPs. By analogy with the computations in \eqref{Fourier exchange formula}, we obtain:
\[
(u \ast X)^{\wedge}=\hat{u} \hat{X},~X\in GSP(G),~ u\in L^1(G).
\]
\end{Def}
Next, we introduce the definition of white noise on $\mathbb{R}^d$.
\begin{Def}\cite{Feichtg2014-whiteNoise}\label{Def-whiteNoise}
A GSP $S:~\mathcal{S}_0(\mathbb{R}^d)\rightarrow H$ is called a centered white noise if
\begin{equation}\label{prop-isometry-whiteNoise}
\exists~c>0, ~\left\langle S(f), S(g) \right\rangle_H =c \int_{\mathbb{R}^d}f\bar{g},~\forall~f, g\in \mathcal{S}_0
(\mathbb{R}^d).
\end{equation}
When $H$ is a closed linear subspace of $L^2_0(\Omega,\mathcal{F}, d\mathbb{P})$ which contains only zero-mean complex Gaussian random variables (centered Gaussian space \cite{Gall2016-Stochastic}), then the white noise $S$ is called the Gaussian white noise. At this point, for any $f\in \mathcal{S}_0$ the random  variable $ S(f) $ obeys the  normal distribution $\text{Norm}\left(0,~ \|f\|^2_{L^2}\right)$ .
\end{Def}
By a density argument, the white noise $S$ can be extended to
 $L^2(\mathbb{R}^d)$ as follows:
\begin{equation}\label{def-pair-noise-L2}
\left\langle S, f \right\rangle
\stackrel{H}{:=}
\lim_{k\rightarrow\infty}
\left\langle S, f_k \right\rangle,
\end{equation}
where $\{f_k\}_k\subset\mathcal{S}_0(\mathbb{R}^d)$ is any sequence converging to $f$ in $L^2(\mathbb{R}^d)$. Thus the white noise $S:~L^2(\mathbb{R}^d)\rightarrow H$ is  an isometry (from deterministic $L^2$ to random $L^2$) up to a scalar multiple \cite{Feichtg2014-whiteNoise}. Specifically, the pairing of the white noise with $L^2$ functions, given in the form of \eqref{def-pair-noise-L2}, is called the \textit{mean-square pairing}.
\begin{Thm}
	The Fourier transform $\hat{S}$ of the white noise $S$ on $\R^d$ also satisfies the isometry property  $\hat{S}: \mathcal{S}_0
	(\mathbb{R}^d)\rightarrow H$, and thus is itself a white noise.
\end{Thm}
\noindent {\bf Proof}: Since the Fourier transform  $\hat{S}$ is a GSP,
then by the isometric property \eqref{prop-isometry-whiteNoise}  we have
\begin{equation}
\left\langle \hat{S}(f), \hat{S}(g) \right\rangle_H
=
\left\langle S(\hat{f}), S(\hat{g}) \right\rangle_H
=c
 \int_{\mathbb{R}^d}\hat{f}\bar{\hat{g}}
= c
 \int_{\mathbb{R}^d}f\bar{g},
\end{equation}
for any $f,~ g\in \mathcal{S}_0(\mathbb{R}^d)$. \hfill\qed \\
Similarly, one can show that the white noise remains the white noise itself under the translations and modulations. The above theorems of GSP  can be generalized to arbitrary LCA groups \cite{Feichtg2014-whiteNoise}.

\subsection{Weighted  Besov and Sobolev spaces}
Since the law of the white noise is invariant by translation (wide sense time stationary), it does not decay at infinity \cite{Feichtg2014-whiteNoise,Chauleur2025}. To deal with this, we introduce weighted Besov spaces and Sobolev spaces.
\begin{Def} \cite{Chauleur2025,Debussche2019}\label{Def-weight-Besov}
Let $\mu, \alpha \in \mathbb{R}$ and $p, q \in[1, \infty)$. The weighted Besov space $B_{p, q, \mu}^\alpha=B_{p, q}^\alpha\left(\langle \cdot\rangle^\mu\right)$ is the set of distribution $u \in\mathscr{S}^{\prime}\left(\mathbb{R}^d\right)$ such that
\[
\|u\|_{B_{p, q, \mu}^\alpha}:=
\left(\sum_{j \geq -1} 2^{j \alpha  q}\left\|\Delta_j u\right\|_{L^p\left( \left\langle \cdot \right\rangle^\mu \right)}^q\right)^{\frac{1}{q}}
<\infty,
 \]
where $\{\Delta_j\}_{j}$ denotes the collection of the  Littlewood-Paley projection operators  \cite{Gavage2007,Chemin2011}.  These are defined by $\Delta_j:=\mathscr{F}^{-1} \eta_j \mathscr{F}$ with $\eta_j \in C_c^\infty(\R^d)$ satisfying the support conditions:
\[\operatorname{Spt} \eta_{-1} \subset\left\{\xi:~ |\xi| \leq 1 \right\},~
\operatorname{Spt} \eta_j \subset\left\{2^{j-1} \leq |\xi| \leq 2^{j+1}\right\},~j\geq 1.
\]
When $p=q=2$, one recovers the  weighted Sobolev spaces $H_\mu^\alpha=B_{2,2, \mu}^\alpha$ with
\begin{equation*}
	\sum_{j \geq -1}
	2^{2j \alpha}\left\|\Delta_j u\right\|_{L^2\left( \left\langle \cdot \right\rangle^\mu \right)}
     ^2
	\simeq
	\| u\|^2_{H^{\alpha}_{\mu}(\mathbb{R}^d)}
	=
	\left\|
	\mathscr{F}^{-1}\langle\cdot\rangle^\alpha \mathscr{F} u
	\right\|^2_{L^2\left( \langle \cdot\rangle^\mu \right)}
	<+\infty.
\end{equation*}
\end{Def}
\begin{rmk}\cite{Tri06}
Set $B_{p, q}^\alpha(\mathbb{R}^d)$ for the unweighted Besov  space $B_{p, q, 0}^\alpha(\mathbb{R}^d)$. A convenient property of the weighted Sobolev space is that the weight can be "pulled in",
\begin{equation}\label{prop-pulled in}
\|u\|_{B_{p, q, \mu}^\alpha(\R^d)}
\simeq
 \left\|\langle\cdot\rangle^\mu u\right\|_{B_{p, q}^\alpha(\R^d)},
\end{equation}
in the sense of equivalent norms, for $\alpha, \mu \in \mathbb{R}$ .
\end{rmk}
Relation \eqref{prop-pulled in} can be used to translate results from the unweighted spaces to their weighted analogues. By the pulled-in property \eqref{prop-pulled in} of the weighted Sobolev norm,  one can establish various properties of the weighted Sobolev spaces under the Fourier transform, translation, modulation, and reflection.
\begin{prop}
The Fourier transform $\mathscr{F}$ interchanges  the time-weight and frequency-weight indices in the weighted Sobolev norm. This duality is concisely expressed as
	\[
	\mathscr{F}H^{\alpha}_{\mu}\simeq  H^{\mu}_{\alpha}.
	\]
In particular, when $\alpha=\mu$, we have $\mathscr{F}H^{\alpha}_{\alpha}\simeq  H^{\alpha}_{\alpha}$, meaning that the Fourier transform $\mathscr{F}$ is an automorphism of  $H^{\alpha}_{\alpha}$.
\end{prop}
\noindent{\textbf{Proof}}: By the pulled-in property \eqref{prop-pulled in}, we have
\begin{equation}
\begin{aligned}
\|f\|_{ H^{\alpha}_{\mu}}=&
\left\|
\mathscr{F}^{-1}\langle\cdot\rangle^{\alpha} \mathscr{F} f
\right\|_{L^2\left( \langle \cdot\rangle^{\mu} \right)}\\
\simeq&
\left\|
\mathscr{F}^{-1}\langle\cdot\rangle^{\alpha} \mathscr{F} \left( \langle \cdot\rangle^{\mu} f\right)
\right\|_{L^2}\\
=&
\left\|
\langle\cdot\rangle^{\alpha} \mathscr{F} \left( \langle \cdot\rangle^{\mu} f \right)
\right\|_{L^2}\\
=&
\left\|
\mathscr{F}  (\langle \cdot\rangle^{\mu} f )
\right\|_{L^2\left( \langle \cdot\rangle^{\alpha} \right)}\\
=&
\left\|
\mathscr{F}(  \langle \cdot\rangle^{\mu} \mathscr{F}^{-1} \left( \mathscr{F} f \right) )
\right\|_{L^2\left( \langle \cdot\rangle^{\alpha} \right)}\\
=&
\| \mathscr{F} f\|_{ H^{\mu}_{\alpha}}, \quad \forall~f\in H^{\alpha}_{\mu}.
\end{aligned}
\end{equation}
 \hfill\qed \\
Moreover, the translation, modulation, and reflection operators are continuous on the weighted Sobolev space ${H}^\alpha_\beta$ for all $\alpha, \beta>0$.
 \begin{prop}
Define the operators
 	\begin{equation*}
 	\begin{aligned}
 	\Rsh_a f(t):=f(t+a),\\
 	\circlearrowleft_a f(t):=e^{i2\pi a\cdot t}f(t),\\
 	f^\checkmark(t):=f(-t),
 	\end{aligned}
 	\end{equation*}	
for  $a\in\R^d$.	
 	Then for any $\alpha,~\beta>0,~a\in\R^d$, and $f\in{H}^\alpha_\beta(\mathbb{R}^d)$, we have:
 	\begin{equation}
 	\begin{aligned}
 	\|\Rsh_a f\|_{H^\alpha_\beta(\mathbb{R}^d)}
 	&\lesssim
 	(1+|a|)^\beta \|f\|_{H^\alpha_\beta(\mathbb{R}^d)},\\
 	\|\circlearrowleft_a f\|_{H^\alpha_\beta(\mathbb{R}^d)}
 	&\lesssim
 	(1+|a|)^\alpha \|f\|_{H^\alpha_\beta(\mathbb{R}^d)},\\
 	\|f^\checkmark\|_{H^\alpha_\beta(\mathbb{R}^d)}
 	&\simeq \|f\|_{H^\alpha_\beta(\mathbb{R}^d)}.\\
 	\end{aligned}
 	\end{equation}
 \end{prop}
\noindent\textbf{Proof}:
We begin with the translation operator. By definition of the weighted Sobolev norm,
 \begin{equation}
 \begin{aligned}
 \|\Rsh_a f\|_{H^\alpha_\beta(\mathbb{R}^d)}
 =&
 \left\|
 \mathscr{F}^{-1}\langle\cdot\rangle^\alpha \mathscr{F} (\Rsh_a f)
 \right\|^2_{L^2\left( \langle \cdot\rangle^\beta \right)}\\
 =&
 \left\|
 \mathscr{F}^{-1}\langle\cdot\rangle^\alpha \circlearrowleft_a (\mathscr{F}  f)
 \right\|^2_{L^2\left( \langle \cdot\rangle^\beta \right)}\\
 =&
 \left\|
 \mathscr{F}^{-1}\langle\cdot\rangle^\alpha \chi_a (\mathscr{F}  f)
 \right\|^2_{L^2\left( \langle \cdot\rangle^\beta \right)},~
 \chi_a(t):=e^{i2\pi a\cdot t}\\	
 =&
 \left\|
 \Rsh_a
 \left(
 \mathscr{F}^{-1}\langle\cdot\rangle^\alpha \mathscr{F}  f
 \right)
 \right\|^2_{L^2\left( \langle \cdot\rangle^\beta \right)}\\
 =&
 \int_{\mathbb{R}^d}
 \left|\left(
 \langle\cdot\rangle^\alpha\hat{ f}
 \right)^{\vee}(t+a) \right|^2 (1+|t|^2)^\beta dt\\
 =&
  \int_{\mathbb{R}^d}
  \left|\left(
  \langle\cdot\rangle^\alpha\hat{ f}
  \right)^{\vee}(t) \right|^2 (1+|t-a|^2)^\beta dt.
  \end{aligned}
  \end{equation}
Using the elementary inequality
 \[
 1 + |t - a|^2 \le
 2(1 + |t|^2)(1 + |a|^2),\quad \forall~t,~a\in\R^d,
 \]
we obtain
 \begin{equation}
 \begin{aligned}
\|\Rsh_a f\|_{H^\alpha_\beta(\mathbb{R}^d)}
 \lesssim&
 (1+|a|)^\beta \|\mathscr{F}^{-1}\langle\cdot\rangle^\alpha \mathscr{F}  f\|_{L^2\left( \langle \cdot\rangle^\beta \right)}\\
 \simeq&
 (1+|a|)^\beta
 \|f\|_{H^\alpha_\beta(\mathbb{R}^d)}.
 \end{aligned}
 \end{equation}
This establishes the continuity of  $\Rsh_a$ on $H^\alpha_\beta$.  Next, for the modulation operator, we use the duality between modulation in the physical space and translation in the frequency space, together with the fact that the Fourier transform is an isomorphism from $H^\alpha_\beta$ to $H^\beta_\alpha$. Then we have
 \begin{equation}
 \begin{aligned}
 \|\circlearrowleft_a f\|_{H^\alpha_\beta(\mathbb{R}^d)}
 \simeq
 \|\mathscr{F}(\circlearrowleft_a f)\|_{H^\beta_\alpha(\mathbb{R}^d)}
 =
 \|\Rsh_{-a} (\mathscr{F} f)\|_{H^\beta_\alpha(\mathbb{R}^d)}.
 \end{aligned}
 \end{equation}
Applying the previously established estimate for translations  gives
 \begin{equation}
 \|\circlearrowleft_a f\|_{H^\alpha_\beta(\mathbb{R}^d)}
 \lesssim
 (1+|a|)^\alpha
 \|\mathscr{F} f\|_{H^\beta_\alpha(\mathbb{R}^d)}
 \simeq
 (1+|a|)^\alpha
 \| f\|_{H^\alpha_\beta(\mathbb{R}^d)},
 \end{equation}
which proves the continuity of the modulation operator.  Finally, the reflection operator $f^\checkmark(t) = f(-t)$ satisfies
\[
\|f^\checkmark\|_{H^\alpha_\beta(\mathbb{R}^d)} = \|f\|_{H^\alpha_\beta(\mathbb{R}^d)},
\]
since both the Lebesgue measure and the weight $\langle t \rangle$ are invariant under  $t \mapsto -t$, and the Fourier transform commutes with reflection up to complex conjugation. Hence it is an isometry on $H^\alpha_\beta(\mathbb{R}^d)$.
 \hfill\qed \\
Thus the weighted Sobolev spaces  ${H}^\alpha_\beta$ with $\alpha,~ \beta>0$ are isometrically  invariant under translations  and modulations-the same property as that of the Segal algebras.
We now recall the duality properties of weighted Besov spaces.
\begin{lem}\label{Lem-weighted Sobolev}\cite{Chauleur2025}
(Duality) For $\alpha,~ \mu \in \mathbb{R}$ and $p, q \in[1, \infty)$, we have the duality
\begin{equation}\label{dual-weighted Sobolev}
\left( B_{p, q, \mu}^\alpha  \right(\R^d))^{\prime}
=
B_{p^{\prime}, q^{\prime}, -\mu}^{-\alpha}(\R^d),
\end{equation}
where $\frac{1}{p}+\frac{1}{p^{\prime}}=1$ and $\frac{1}{q}+\frac{1}{q^{\prime}}=1$.
\end{lem}
Thus, both the Segal algebra $\mathcal{S}_0(\mathbb{R}^d)$ and the weighted Sobolev spaces $H^\alpha_\beta(\mathbb{R}^d)$ ($\alpha,~ \beta > \tfrac{d}{2}$) are invariant under translations, modulations, and the Fourier transform.
\subsection{Relationship between the Segal algebra and the weighted Sobolev spaces}
We now investigate the relationship between the Segal algebra and the weighted Sobolev spaces. It starts with the Fourier-Lebesgue spaces.
\begin{Def}\cite{Beyni2020}
	Let $1 \leq p \leq \infty$. The Fourier-Lebesgue space $\mathscr{F} L^p$ is the Banach space of tempered distributions $f \in \mathscr{S}^{\prime}$ endowed with the norm
	\[ 	\|f\|_{\mathscr{F}L^p}=\|\hat{f}\|_{L^p}.
	\]
	More generally, for $s \geq 0$, we define the weighted Fourier-Lebesgue spaces $\mathscr{F} L^p(\langle\cdot\rangle^s)$ via the norm
	\[ 	\|f\|_{\mathscr{F} L^p(\langle\cdot\rangle^s)}=
	\left\|\langle\cdot\rangle^s \hat{f}\right\|_{L^p}.
	\]
\end{Def}
To explore the relationship between the weighted Sobolev spaces and the Segal algebras, we first recall that the weighted time-frequency  $ L^2 $ spaces $L^2\left(\langle\cdot\rangle^s\right) \cap
\mathscr{F} L^2\left( \langle\cdot\rangle^s\right)$ are embed  into the Segal algebra.
\begin{prop}\cite{Grochenig2001}
	\label{time-frequency L2 vs S_0}
	If $f \in L^2\left(\mathbb{R}^d, \langle\cdot\rangle^s\right)$ and $\hat{f} \in L^2\left(\mathbb{R}^d, \langle\cdot\rangle^s\right)$ for some $s>d$. Then $f \in \mathcal{S}_0(\mathbb{R}^d)$. Equivalently,
	\[
	L^2\left(\mathbb{R}^d, \langle\cdot\rangle^s\right) \cap
	\mathscr{F} L^2\left(\mathbb{R}^d, \langle\cdot\rangle^s\right)
	\subseteq
	\mathcal{S}_0(\mathbb{R}^d).
	\]
	More precisely, for any $ f\in
	L^2\left(\mathbb{R}^d, \langle\cdot\rangle^s\right) \cap
	\mathscr{F} L^2\left(\mathbb{R}^d, \langle\cdot\rangle^s\right) $, we have
	\[
	\|f\|_{\mathcal{S}_0}
	\lesssim
	\|f\|_{L^2\left(\langle\cdot\rangle^s\right)}
	+\|\hat{f}\|_{L^2\left( \langle\cdot\rangle^s\right)}.
	\]
\end{prop}
Consider the Poisson kernel $P_y(x)\simeq \frac{y}{(y^2+|x|^2)^{\frac{1+d}{2}}}$ on  $\mathbb{R}^{1+d}_+$. For all $d<s<\frac{d}{2}+1$, its weighted time-frequency norm satisfies
\begin{equation}
\begin{aligned}
	&\|P_y\|^2_{L^2\left(\langle\cdot\rangle^s\right)}
	+\|\widehat{P_y}\|^2_{L^2\left( \langle\cdot\rangle^s\right)} \\
\simeq&
\int_{\mathbb{R}^d} \frac{y^2\langle x\rangle^{2s}}{(y^2+|x|^2)^{1+d}}	dx
+
\int_{\mathbb{R}^d} e^{4\pi y|\xi|}\langle \xi\rangle^{2s}	d\xi.
\end{aligned}
\end{equation}
Consequently, for $d=1$, the kernel $P_y$ lies in $\mathcal{S}_0(\mathbb{R}^1)$. In fact, both the spaces $L^2\left(\mathbb{R}^d, \langle\cdot\rangle^s\right) $ and $
\mathscr{F} L^2\left(\mathbb{R}^d, \langle\cdot\rangle^s\right)$ can be viewed as special cases of the weighted Besov spaces. Using the Littlewood-Paley decomposition, for all
$f\in H^\alpha_\beta(\mathbb{R}^d)$ with $\alpha,~\beta\geq s>0$, we have
\[
\|f\|^2_{L^2\left(\langle\cdot\rangle^s\right)}
+
\|\hat{f}\|^2_{L^2\left( \langle\cdot\rangle^s\right)}
\lesssim 2
\sum_{j \geq -1} 2^{j 2\alpha  }
\left\|\Delta_j f\right\|_{L^2\left( \left\langle \cdot \right\rangle^\beta \right)}^2
\simeq
\|f\|^2_{H^\alpha_\beta}
<+\infty.
\]
From Proposition \ref{time-frequency L2 vs S_0}, we obtain  the embedding:
\[
H^\alpha_\beta(\mathbb{R}^d)
\subseteq
L^2\left(\mathbb{R}^d, \langle\cdot\rangle^s\right)
\cap
\mathscr{F} L^2\left(\mathbb{R}^d, \langle\cdot\rangle^s\right)
\subseteq
\mathcal{S}_0(\mathbb{R}^d),
\]
for all $\alpha,~\beta\geq s>d$. In this paper we prove a refinement  of the indices of the weighted Sobolev spaces in the above embedding relations. We first recall a pulled-in result concerning the weighted modulation space $M^{\alpha,\beta}_{p,q}(\mathbb {R}^d)$.
\begin{prop}\cite{Cordero2020}\label{prop-pull-in-modulation}
Consider $1 \leq p, q \leq \infty,~ \alpha\geq0$ and $s \in \mathbb{R}$. For $v_s(x)=\langle{x}\rangle^s$, the mapping $f \mapsto v_s f$ is a topological isomorphism from $M^{\alpha,s}_{p,q}(\mathbb {R}^d)$ to $M^{\alpha,0}_{p,q}(\mathbb {R}^d)$, satisfying
\[
 \|f\|_{M^{\alpha,s}_{p,q}}\simeq\|v_s f \|_{M^{\alpha,0}_{p,q}}.
 \]
\end{prop}
This Proposition yields that:
\begin{prop}
	For all  $\alpha,~\beta>\frac{d}{2}$, the embedding
	$
	H^\alpha_\beta(\mathbb{R}^d)
	\subseteq
	\mathcal{S}_0(\mathbb{R}^d)
	$
	holds.
\end{prop}
\noindent{\textbf{Proof}}:
For $f\in M^{\alpha,\beta}_{2,2}(\mathbb{R}^d)$ with $\alpha,~\beta>\frac{d}{2}$, applying the Cauchy-Schwarz inequality twice yields:
\begin{equation}
\begin{aligned}
\|f\|_{\mathcal{S}_0(\mathbb{R}^d)}\simeq&
\sum_{k\in\mathbb {Z}^d}
\Vert\Box _k f\Vert _{L^1(\mathbb {R}^d)} ,\\
\leq&
\sum_{k\in\mathbb {Z}^d}
\Vert\Box _k f\Vert _{L^2(\langle\cdot\rangle^{\beta} )}
\left(\int_{\mathbb{R}^d}\left\langle {x }\right\rangle^{-2\beta} dx \right)^{1/2},\\
\simeq&
\sum_{k\in\mathbb {Z}^d}
\Vert\Box _k f\Vert _{L^2(\langle\cdot\rangle^{\beta} )},\\
\leq&
\left\|
 \langle{k}\rangle^{\alpha}
  \Vert\Box _k f\Vert_{L^2(\langle\cdot\rangle^{\beta} )}
\right\|_{\ell ^2_k({\mathbb {Z}^d})}
  \left(\sum_{k\in\mathbb {Z}^d}
  \langle{k}\rangle^{-2\alpha}
  \right)^{1/2},\\
\lesssim&
\|f\|_{M^{\alpha,\beta}_{2,2}(\mathbb{R}^d)}.
\end{aligned}
\end{equation}
Using the pulled-in property \ref{prop-pull-in-modulation} of weighted modulation spaces and the relation between modulation spaces and Sobolev spaces $M^{\alpha,0}_{2,2}\cong H^\alpha$, we obtain:
\begin{equation}
\|f\|_{\mathcal{S}_0(\mathbb{R}^d)}
\lesssim
\|\langle\cdot\rangle^{\beta}f\|_{M^{\alpha,0}_{2,2}(\mathbb{R}^d)}
\simeq
\|\langle\cdot\rangle^{\beta}f\|_{H^\alpha(\mathbb{R}^d)}
\simeq
\|f\|_{H^\alpha_\beta(\mathbb{R}^d)},
\end{equation}
provided that $\alpha,~\beta>\frac{d}{2}$.
 \hfill\qed \\
\subsection{Regularity of the white noise}
Using the isometry \eqref{prop-isometry-whiteNoise} and the Littlewood-Paley characterization of weighted Besov spaces from \cite{Chauleur2025}, one can obtain weighted H\"older regularity for the white noise.  Here, however, we follow the approach of \cite{Chauleur2025} to establish the weighted Sobolev regularity of the white noise in $H^\alpha_\beta$.  This is natural, since the weighted Sobolev spaces $H^\alpha_\beta(\mathbb{R}^d)$ ($\alpha,\beta > \tfrac{d}{2}$) and the Segal-Feichtinger algebra $\mathcal{S}_0(\mathbb{R}^d)$ are both invariant under translations, modulations, and the Fourier transform.
In what follows, the white noise $S$ is a GSP on $\mathbb{R}^d$ taking values in the  Hilbert space   $H=L^2_0(\Omega,\mathcal{F}, d\mathbb{P})$.
\begin{prop}
	\label{Prop-Sobolev-rglrt-Noise}(weighted Sobolev Regularity)
For $\alpha,~\mu>\frac{d}{2}$, the  white noise $S$  belongs to the weighted Sobolev space $  H^{-\alpha}_{-\mu}(\mathbb{R}^d) $ almost surely.
\end{prop}
\noindent {\bf Proof}:
For almost every $\omega$, the white noise sample path $S_\omega$ defines a linear functional on $ \mathscr{S}(\R^d) $. Thus, for any $x\in\R^d$, the pairings
\[
\left\langle
S_\omega, \eta^{\vee}_j(x-\cdot) \right\rangle,~j\geq-1
 \]
are well-defined.
Since the  white noise $S$ is an isometry from $L^2(\R^d)$ to $L^2(\Omega)$, then for $j\geq-1$ and $\mu>\frac{d}{2}$ we have
\begin{equation*}
	\begin{aligned}
		\mathbb{E}\left\|\Delta_j S\right\|^2_{L^2\left(\langle \cdot\rangle^{-\mu}\right)} &
		=\int_{\R^d} \mathbb{E}_\omega
		\left|\left\langle
		S_\omega, \eta^{\vee}_j(x-\cdot) \right\rangle
		\right|^2
		\langle x\rangle^{-2 \mu} dx,~\eta_j\in {C_c^\infty(\mathbb{R}^d)} \\
		&\simeq
		\left\|\eta^{\vee}_j\right\|_{L^2(\R^d)}^2 \int_{\R^d}\langle x\rangle^{-2 \mu} dx,\\
		&=
		C_\mu 2^{jd} \|\eta^{\vee}\|_{L^2}^2,
	\end{aligned}
\end{equation*}
where  $\Delta_j=\mathscr{F}^{-1} \eta_j \mathscr{F}$ is the frequency  projection as in Definition \ref{Def-weight-Besov}.
Hence for $\alpha,~\mu >\frac{d}{2}$, we get
\begin{equation*}
\mathbb{E}  \|S\|^2_{ H^{-\alpha}_{-\mu}(\R^d)}
\simeq
\mathbb{E}\sum_{j\geq-1}  2^{-2j\alpha}
\left\|\Delta_j S\right\|^2_{L^2\left(\langle \cdot\rangle^{-\mu}\right)}
\simeq
C_\mu
\sum_{j\geq-1}
4^{-j(\alpha-\frac{d}{2})}
<\infty.
\end{equation*}
This implies that for almost every $\omega$,
\begin{equation}
\begin{aligned}
 \|S_\omega\|_{ H^{-\alpha}_{-\mu}}
 =&
 \sup \left\{
 	\left|\left\langle
 	S_\omega, \varphi \right\rangle
 	\right|:~\|\varphi\|_{H^{\alpha}_{\mu}}=1
 	\right\}\\
 	=&
 \sup \left\{
 \left|\left\langle
 S, \varphi \right\rangle\mid_\omega
 \right|:~\|\varphi\|_{H^{\alpha}_{\mu}}=1
 \right\}\\
 <&\infty,\quad \forall~\varphi\in\mathscr{S}(\R^d).
\end{aligned}	
\end{equation}
Thus, owing to the finiteness of its weighted Sobolev norm, the  sample path $S_\omega$ is  a continuous linear functional on $\mathscr{S}(\mathbb{R}^d)$ for almost every $\omega$.  The  mean-square pairing induced by $S$ almost surely takes the form  $\langle\cdot,\cdot\rangle_{H^{-\alpha}_{-\mu}\times H^{\alpha}_{\mu}}$ when the test functions belong to $H^{\alpha}_{\mu}(\mathbb{R}^d)$ with $\alpha,~ \mu>\frac{d}{2}$, and also that of   $\langle\cdot,\cdot\rangle_{\mathscr{S}^{'}\times\mathscr{S}}$.
\hfill\qed \\
\begin{rmk}	
The Gaussian white noise on the non-compact domain $\R^d$ also belongs to the  Sobolev-type space $ H^{-\alpha}_{-\mu}(\R^d)$, which is similar to the unit circle $\partial \mathbb{D}$. However, the former requires an decay weight function $\langle x \rangle^{-\mu}$, which is absent in the compact   setting.
\end{rmk}

\subsubsection{The mean-square convolution of the white noise with an $ L^2$ function}
In this section, we define  convolution of the white noise with $L^2$ functions in the mean-square sense and provide its pointwise functional expression.

Suppose $f\in L^2(\R^d)$. Take a sequence $\{f_k\}_k \subset\mathscr{S}(\R^d)$ such that $f_k{\rightarrow} f$ in  $L^2(\R^d)$. Then for the white noise $S$, the sequence  of random variables
\[
\langle S \ast f_k, \varphi\rangle
=
\left\langle S , {f_k}^{\checkmark}\ast\varphi
\right\rangle,
~ \varphi \in \mathscr{S}(\R^d),
\]
is Cauchy  in $L^2(\Omega, \mathrm{d}\mathbb{P})$ as
\[
\begin{aligned}
\|
\langle S \ast f_k, \varphi\rangle-\langle S \ast f_m, \varphi\rangle
\|_{L^2(\Omega)}
=&
\left\|
\left\langle S , (f_k-f_m)^{\checkmark}\ast\varphi
\right\rangle
\right\|_{L^2(\Omega)},\\
\simeq& \left\|(f_k-f_m)^{\checkmark}\ast\varphi
\right\|_{L^2(\R^d)},\\
\leq&
\|f_k-f_m\|_{L^2(\R^d)}\|\varphi\|_{L^1(\R^d)},
\end{aligned}
\]
where $\varphi \in \mathscr{S}(\R^d)\subseteq L^1(\R^d)$. Thus the convolution of the white noise with an $L^2$ function can be defined as
\begin{equation}
\begin{aligned}
\langle S \ast f, \varphi\rangle
:=
\lim_{k\to\infty}
\langle S \ast f_k, \varphi\rangle
=
\lim_{k\to\infty}
\langle S ,  f_k^\checkmark\ast\varphi\rangle
=:
\langle S , f^\checkmark\ast\varphi\rangle,
\end{aligned}
\end{equation}
for any test function $ \varphi\in\mathscr{S}(\R^d) $.
This limit is independent of the  sequence $\{f_k\}_k\subset\mathscr{S}(\R^d)$ with $f_k \stackrel{L^2}{\rightarrow}  f$. For a.s. $\omega\in\Omega$, we denote
\[
\begin{aligned}
\langle S_\omega \ast f, \varphi\rangle
:=&
\langle S \ast f, \varphi\rangle\arrowvert_\omega\\
=&
(\lim_{k\rightarrow\infty}
\langle S \ast f_k, \varphi\rangle)\arrowvert_\omega\\
=&
(\lim_{k\rightarrow\infty}
\langle S , f_k^\checkmark\ast\varphi\rangle
)\arrowvert_\omega, \quad
\forall~ \varphi \in \mathscr{S}(\R^d).
\end{aligned}
\]
Similarly, the Fourier transform of the convolution between white noise and an $L^2$ function can also be defined via an $L^2(\Omega)$-limit:
\begin{equation*}
\begin{aligned}
\langle (S \ast f)^\wedge, \varphi\rangle
:=
\lim_{k\to\infty}
\langle (S \ast f_k)^\wedge, \varphi\rangle
=
\lim_{k\to\infty}
\langle \hat{S} \hat{ f_k}, \varphi\rangle
=:
\langle \hat{S}\hat{f}  ,\varphi\rangle,
\end{aligned}
\end{equation*}
for any  $ \varphi\in\mathscr{S}(\R^d) $.
As is well known, for $F\in \mathscr{S}^{'}(\R^d)$ and $\varphi\in \mathscr{S}(\R^d)$, the convolution $F\ast \varphi$ belongs to the  slowly growing  smooth function space $\mathscr{O}(\R^d)$ \cite{Neto1973}. Moreover,
\[ \forall~x\in\R^d,~
(F\ast \varphi)(x)=\left\langle F, {\Rsh_x} \varphi^{\checkmark} \right\rangle.
\]
We now show that the convolution of the white noise with an  $L^2$   function exhibits similar properties.
\begin{Thm}\label{Thm-L2-convolution of white noise}
	Let \( S \) be the white noise and \( f \in L^2(\mathbb{R}^d) \). Define
	\[
	(S \star f)(x) := \left\langle S, {\Rsh_x} f^{\checkmark} \right\rangle
	=
	\left\langle S, f(x-\cdot) \right\rangle
	\]	
	for each fixed 	$x \in \mathbb{R}^d$, where the pairing is understood in the sense of taking $L^2(\Omega)$-limits as in \eqref{def-pair-noise-L2}. Then, for almost every \( \omega \in \Omega \), the function $S_\omega \star f$ belongs to $L^2({\mathbb {R}^d},\langle\cdot\rangle^{-d})\subseteq\mathscr{S}^{'}(\R^d)$ and satisfies
	\[
	\langle S_\omega \star f, \varphi\rangle
	=
	\langle S_\omega \ast f, \varphi\rangle,~
	\forall~ \varphi \in \mathscr{S}(\R^d).
	\]
\end{Thm}
See Appendix \ref{prf-L2-convolution} for a proof.\\
The above theorem implies that for every $f\in{L^2}(\R^d)$ and almost all samples $S_\omega$, the pathwise convolution
\[
 S_\omega \ast f \in L^2({\mathbb {R}^d},\langle\cdot\rangle^{-d})\subseteq\mathscr{S}^{'}(\R^d).
 \]
 Moreover, the pathwise convolution satisfies
\[  (S_\omega \ast f)(x)=\left\langle S_\omega, {\Rsh_x} f^{\checkmark} \right\rangle,~\forall~x\in\R^d. \]
\subsubsection{Applications of White Noise Convolution with $L^2$ function}
For any fixed point $(y,x)\in\mathbb{R}^{1+d}_+$, the Poisson kernel $P_{y+\underline{x}}$ is in $L^1\cap L^2(\mathbb{R}^d)$. Let  $H=L^2_0(\Omega,\mathcal{F}, d\mathbb{P})$. Consequently, for any approximating sequence
$\{\varphi_{k}\}\subset \mathcal{S}_0(\mathbb{R}^d)$ with $\varphi_{k}\stackrel{L^2(\mathbb{R}^d)}{\longrightarrow}P_{y+\underline{x}}$, the bilinear pairing
\begin{equation}\label{FormalSolution}
u(y+\underline{x})
:=
\left\langle S, P_{y+\underline{x}}\right\rangle
\stackrel{H}{=}\lim_{k\rightarrow\infty}
\left\langle S, \varphi_{k}\right\rangle
\end{equation}
is  well-defined as an element of $H$. This bilinear pairing \eqref{FormalSolution} will be regarded as the formal solution of Problem \eqref{white-noise-HalfPlane}. The  covariance  of $u$ is given by
\begin{equation}\label{Cov_solution}
\begin{aligned}
\mathcal{C}_u(y+\underline{x}, s+\underline{t})&=
\mathbb{E}
[
\left\langle S, P_{y+\underline{x}}\right\rangle
\overline{
	\left\langle S, P_{ s+\underline{t}}\right\rangle	
}
]
=
\left\langle  P_{y+\underline{x}}, P_{ s+\underline{t}}\right\rangle_{L^2}
=P_{s+y}(\underline{x}-\underline{t}),
\end{aligned}
\end{equation}
where $\mathbb{E}u\equiv0$.
\begin{rmk}
Alternative to the Littlewood-Paley dyadic decomposition, we give a more accurate calculation of the white noise's Bochner-type norm $L^2(\Omega, H^{-\alpha}_{-\mu}(\R))$. For $ \mu>\frac{1}{2} $,  the Bessel potential kernel $G_\mu$ satisfies: 	$\widehat{G_\mu}(\xi)=\left\langle \xi\right\rangle^{-\mu}\in L^2(\mathbb{R})\cap A(\mathbb{R}) $. The $L^2(\Omega, H^{-\alpha}_{-\mu}(\R))$ norm of the white noise is given by
\begin{equation}\label{StatisticalAverageSobolev-EnergyNorm-WhiteNoise}
\begin{aligned}
\mathbb{E}  \|S\|^2_{ H^{-\alpha}_{-\mu}(\R)}
=&
\mathbb{E}
\left\|
\mathscr{F}^{-1}\langle\cdot\rangle^{-\alpha} \mathscr{F} S
\right\|^2_{L^2\left( \langle \cdot\rangle^{-\mu} \right)}\\
\simeq&
\mathbb{E}
\left\|
\langle\cdot\rangle^{-\mu}  S
\right\|^2_{H^{-\alpha}(\R) } \quad
(\text{pulled-in property})\\
=&
\mathbb{E}
\left\|
\left(\langle\cdot\rangle^{-\mu}  S \right)^{\wedge}
\right\|^2_{L^2\left( \langle \cdot\rangle^{-\alpha} \right) }\\
=&
\mathbb{E}
\left\|
G_\mu\ast\hat{S}
\right\|^2_{L^2\left( \langle \cdot\rangle^{-\alpha} \right) } \quad \left(
{G_\mu},~\widehat{G_\mu}\in L^2(\R)
 \right)  \\
=&
\mathbb{E}
\int_{\mathbb{R}}
\left| \left\langle \hat{S}, G_\mu(x-\cdot) \right\rangle  \right|^2 \left\langle {x }\right\rangle^{-2\alpha} dx\\
=&
\int_{\mathbb{R}} \mathbb{E}
\left| \left\langle S, \left(G_\mu(x-\cdot)\right)^{\wedge} \right\rangle  \right|^2 \left\langle {x }\right\rangle^{-2\alpha} dx\\
\simeq&
\|\widehat{G_\mu}\|^2_{L^2(\mathbb{R})} \int_{\mathbb{R}}\left\langle {x }\right\rangle^{-2\alpha} dx\quad
(\text{isometric property})\\
=&
\int_{\mathbb{R}}\left\langle \xi\right\rangle^{-2\mu} d\xi
\int_{\mathbb{R}}\left\langle {x }\right\rangle^{-2\alpha} dx\\
\simeq&
\mathfrak{B}\left(\frac{1}{2}, \alpha-\frac{1}{2}\right)
\mathfrak{B}\left(\frac{1}{2}, \mu-\frac{1}{2}\right).
\end{aligned}
\end{equation}
The norm $L^2(\Omega, H^{-\alpha}_{-\mu}(\R))$ in \eqref{StatisticalAverageSobolev-EnergyNorm-WhiteNoise} can be interpreted as a statistically averaged weighted Sobolev (energy) norm of the white noise.
\end{rmk}
For higher-dimensional cases, analogous to the one-dimensional case, the white noise $S$ satisfies the  norm estimate
\[
\mathbb{E}  \|S\|^2_{ H^{-\alpha}_{-\mu}(\R^d)}
\simeq
\int_{\mathbb{R}^d}\left\langle \xi\right\rangle^{-2\mu}
\int_{\mathbb{R}^d}\left\langle {x }\right\rangle^{-2\alpha} \simeq
\mathfrak{B}\left(\frac{d}{2}, \alpha-\frac{d}{2}\right)
\mathfrak{B}\left(\frac{d}{2}, \mu-\frac{d}{2}\right),
 \]
provided $\alpha,~\mu>\frac{d}{2}$.
\subsection{Foundations of Clifford analysis}
Most of the basic knowledge and notation recalled in this section are referred to \cite{Delanghe1992-Clifford, Qian2019-SingularLipschitz, McIntosh1996}. Let $\mathbf{e}_1, \ldots, \mathbf{e}_n$ be basic elements satisfying the multiplication rule:
\begin{equation}\label{Prop-Clifford-product}
	\mathbf{e}_l \mathbf{e}_j+\mathbf{e}_j \mathbf{e}_l=-2 \delta_{l j},~l, j=1,2, \ldots, n,
\end{equation}
with Kronecker delta function
\[
\delta_{l j}=\left\{
\begin{aligned}
1&,~~ l=j, \\
0&,~~ l \neq  j.
\end{aligned}\right.
\]
In other words,
\begin{equation*}
	\left\{\begin{aligned}
	{\mathbf{e}_j}^2&=-1,~1\leqslant j\leqslant n, \\
		\mathbf{e}_l \mathbf{e}_j&=-\mathbf{e}_j \mathbf{e}_l,~~1\leqslant j< l\leqslant n.
	\end{aligned}\right.
\end{equation*}
Let
\[
\mathbb{R}^n=\left\{\underline{x}=x_1 \mathbf{e}_1+\cdots+x_n \mathbf{e}_n:~x_j \in \mathbb{R},~j=1,2, \ldots, n\right\}
\]
be identical with the usual Euclidean space $\mathbb{R}^n$. Set
\[
\mathbb{R}^{1+n}=\left\{x_0+\underline{x}:~x_0 \in \mathbb{R},~ \underline{x} \in \mathbb{R}^n\right\},
\]
and it is divided into upper- and lower-half spaces
\[
\mathbb{R}_{+}^{1+n}=\left\{x_0+\underline{x}\in\mathbb{R}^{1+n}: x_0>0\right\},~
\mathbb{R}_{-}^{1+n}=\left\{x_0+\underline{x}\in\mathbb{R}^{1+n}: x_0<0\right\}.
\]
An element in $\mathbb{R}_1^n$ is called a vector. The real (or complex) Clifford algebra generated by $\mathbf{e}_1, \mathbf{e}_2, \ldots, \mathbf{e}_n$, denoted by $\mathbb{R}^{n}$ (or $\mathbb{C}^{n}$), is the associative algebra generated by $\mathbf{e}_1, \mathbf{e}_2, \ldots, \mathbf{e}_n$, over the real (or complex) field $\mathbb{R}$ (or $\mathbb{C}$). A general element in $\mathbb{R}^{n}$, therefore, is of the form $x=\sum_\gamma x_\gamma \mathbf{e}_\gamma$, where $\mathbf{e}_\gamma=$ $\mathbf{e}_{i_1} \mathbf{e}_{i_2} \cdots \mathbf{e}_{i_l}$, and $\gamma$ runs over all the ordered subsets of $\{1,2, \ldots, n\}$, namely
\[
\gamma=\left\{1 \leqslant i_1<i_2<\cdots<i_l \leqslant n\right\}, \quad 1 \leqslant l \leqslant n .
\]
Let the conjugate of $x=\sum_\gamma x_\gamma \mathbf{e}_\gamma\in \mathbb{C}^n$ be $\bar{x}=\sum_\gamma \overline{x_\gamma} \bar{\mathbf{e}}_\gamma$, where $\bar{\mathbf{e}}_\gamma=\bar{\mathbf{e}}_{i_1} \bar{\mathbf{e}}_{i_2} \cdots \bar{\mathbf{e}}_{i_l}$ with $\bar{\mathbf{e}}_0=\mathbf{e}_0$ and $\bar{\mathbf{e}}_j=-\mathbf{e}_j$ for $j\neq 0$. Denote the scalar part and the non-scalar part of $ x $ by
\[
\text{Sc}(x):=x_0,\quad \text{NSc}(x):=x-x_0,
\]
respectively. The natural inner product between $x$ and $y$ in $\mathbb{C}^{n}$, denoted by $\langle x, y\rangle$, is the complex number $\sum_\gamma x_\gamma \overline{y_\gamma}$, where $x=\sum_\gamma x_\gamma \mathbf{e}_\gamma$ and $y=\sum_\gamma y_\gamma \mathbf{e}_\gamma$. The norm associated with this inner product is
\[
|x|=\langle x, x\rangle^{\frac{1}{2}}=\left(\sum_\gamma\left|x_\gamma\right|^2\right)^{\frac{1}{2}}= \left(\text{Sc}(\bar{x}x) \right)^{1/2}.
\]
If $x, y$ are vectors, then $  |x y|=|x||y|  $.  The conjugate of a vector $x=x_0+\underline{x}$ is defined as $\bar{x}=x_0-\underline{x}$. An important observation is the inverse formula
\[
x^{-1}=\frac{\bar{x}}{|x|^2}
\]
whenever $x$ is a non-zero vector in $ \mathbb{R}_1^n$. However, for a general Clifford number $x$, the inverse of $x$ may not exist. That is to say, the Clifford algebra is not a division algebra.  In below we will study functions defined in $\mathbb{R}^n$ or $\mathbb{R}^{1+n}$ taking values in $\mathbb{C}^{n}$. So, they are of the form $f(x)=\sum_\gamma f_\gamma(x) \mathbf{e}_\gamma$, where $f_\gamma$ are complex valued functions.

Consider the Dirac operator (generalized Cauchy-Riemann operator)
\begin{equation}\label{def_Dirac-operator}
	{D}=D_0+\underline{D},
\end{equation}
where $D_0=\frac{\partial}{\partial x_0}$ and $\underline{D}=\left(\frac{\partial}{\partial x_1}\right)\mathbf{e}_1+\cdots+\left(\frac{\partial}{\partial x_n}\right) \mathbf{e}_n$ may be considered as the Dirac operator in $\mathbb{R}^n$. To be symmetric, we also write $D_0=\frac{\partial}{\partial x_0}=\left(\frac{\partial}{\partial x_0}\right) \mathbf{e}_0$, with $\mathbf{e}_0=1$. We define the ``left" and``right" roles of the operators $D$ by
\begin{equation}\label{def_left-right-Dirac}
	{D} f=\sum_{j=0}^n \sum_S \frac{\partial f_S}{\partial x_j} \mathbf{e}_j \mathbf{e}_S,\quad
	f {D}=\sum_{j=0}^n \sum_S \frac{\partial f_S}{\partial x_j} \mathbf{e}_S \mathbf{e}_j.
\end{equation}
In particular, if
$ {D} f=0 $
holds  in a domain (open and connected) $\Omega$, then we say that $f$ is left-monogenic in $\Omega$, and if $f {D}=0$ in $\Omega$, we say that $f$ in right-monogenic in $\Omega$. If $f$ is both left-monogenic and right-monogenic, then we say that $\mathrm{f}$ is monogenic. In particular, if a  vector-valued function  $F(x)=f_0(x)+\sum_{j=1}^n f_j(x) \mathbf{e}_j$ is monogenic in $\mathbb{R}^n_1$, then $F$ satisfies the generalized Cauchy-Riemann (C-R) systems \cite{Stein-Weiss1971}. To see this, let
$ 	U\left(x_0+\underline{x}\right)=u_0\left(x_0+\underline{x}\right)-\sum_{j=1}^{n}u_j(x_0+\underline{x}) \mathbf{e}_j$, 
then we have
\begin{equation}\label{Cauchy-Riemann_Equt}
	\begin{aligned}
		{D}U (x_0+\underline{x})
		& =\left(\sum_{k=0}^n \mathbf{e}_k \partial_k \right)\left[u_0\left(x_0+\underline{x}\right) \mathbf{e}_0-
		\sum_{j=1}^{n}u_j(x_0+\underline{x}) \mathbf{e}_j
		\right], \\
		& =
		\sum_{m=0}^{n}
		\frac{\partial u_m}{\partial x_m}+
		\sum_{0 \leq j<k \leq n}\left(\frac{\partial u_j}{\partial x_k}-\frac{\partial u_k}{\partial x_j}\right) \mathbf{e}_j \mathbf{e}_k.
	\end{aligned}
\end{equation}
with $\mathbf{e}_0=1$.  This  also implies that  $U$ is left-monogenic if and only if it is  right-monogenic. Through  the C-R system, the vector-valued monogenic function theory is closely related  to the conjugate harmonic systems of Stein-Weiss \cite{Stein-Weiss1971}.
\section{The Sokhotski-Plemelj formula for the white noise density}\label{Sec-Sokhotski-Plemelj formula}
This section extends the Sokhotski-Plemelj formula to the setting of the white noise density. We proceed in two steps: first, by verifying the monogenicity of the  Cauchy transform of the white noise in the upper half-space; and second, by showing that its jump difference  across the upper and lower half-spaces converges weakly in the mean-square sense to the original white noise as the jump height tends to zero. This ultimately proves that the Cauchy transform serves as an monogenic representation of the white noise.
\subsection{Harmonic and Monogenic Lifting for the White Noise}\label{Sec_Cauchy-trsfm}
We now turn to the monogenicity of the Cauchy transform of the white noise.  Consider the Poisson kernel on $\R^{1+d}_+$,
\begin{equation}\label{Poss-kernel}
P_y(\undl{x})
=
c_d
\frac{y}{ (y^2 + |\undl{x} |^2)^{\frac{d+1}{2}}}
\simeq
\frac{y}{|y+\underline{x}|^{1+d}},
\end{equation}
and its conjugate counterpart
\begin{equation}\label{Conj-Poiss-kernel}
Q_y(\undl{x})
=
c_d
\sum_{j=1}^{d}
\frac{x_j }{(y^2 + |\undl{x}|^2)^{\frac{d+1}{2}}} \textbf{e}_j
=:
\sum_{j=1}^{d}Q^{(j)}_y(\undl{x})\textbf{e}_j
\simeq
\frac{\undl{x}}{|y+\underline{x}|^{1+d}}.
\end{equation}
We introduce the shifted kernels:
\[
\begin{aligned}
P_{y+\undl{x}}(\undl{t})
:=&
P_{y}(\undl{t}-\undl{x})
=P_{y}(\undl{x}-\undl{t})=P_{y+\undl{t}}(\undl{x}),\\
Q^{(j)}_{y+\undl{x}}(\undl{t})
:=&
Q^{(j)}_{y}(\undl{t}-\undl{x})
=-Q^{(j)}_{y}(\undl{x}-\undl{t})
=-Q^{(j)}_{y+\undl{t}}(\undl{x}),~
\forall~1\leq{j}\leq {d} , \\
Q_{y+\undl{x}}(\undl{t}):=&
\sum_{j=1}^{d}Q^{(j)}_{y+\undl{x}}(\undl{t})\textbf{e}_j
=
-\sum_{j=1}^{d}Q^{(j)}_{y+\undl{t}}(\undl{x})\textbf{e}_j
=-Q_{y+\undl{t}}(\undl{x})
=\overline{Q_{y+\undl{t}}(\undl{x})}.
\end{aligned}
\]
In particular, $ P_{y+\undl{0}}(\undl{t})=P_{y}(\undl{t}),~  Q^{(j)}_{y+\undl{0}}(\undl{t})=Q^{(j)}_{y}(\undl{t}) $. Combining the Poisson and conjugate Poisson kernels yields the Cauchy kernel at the point $x:=x_0+\underline{x}\in\R^{1+d}_+$ \cite{Dang2020-Clifford-Hardy,Qian2009-Hilbert-Tsfm,Delanghe1992-Clifford},
\begin{equation}
\begin{aligned}
\frac{1}{{c_d} }
\left(
P_{x_0+\undl{x}}(\undl{t})-Q_{x_0+\undl{x}}(\undl{t})
\right)
=&
\frac{x_0+\underline{x}-\underline{t}}{|x_0+\underline{x}-\underline{t}|^{1+d}}\\
=:&
\mathscr{C}_{x_0}(\underline{x}-\underline{t})
\\
=:&
\mathscr{C}_{x_0+\underline{t}}(\underline{x})\\
=&
\overline{
	\mathscr{C}_{x_0+\underline{x}}(\underline{t})	
}.
\end{aligned}
\end{equation}
For any $f \in L^2(\mathbb{R}^d)$, its Cauchy transform is given by
\begin{equation}
\left\langle f,   \mathscr{C}_{x_0+\underline{x}}\right\rangle_{L^2(\mathbb{R}^d)}
=
\int_{\R^d} f ~\overline{\mathscr{C}_{x_0+\underline{x}}}
=
(f\ast \mathscr{C}_{x_0} )(\underline{x})
=
\int_{\R^d} f(\undl{t}) \mathscr{C}_{x_0}(\underline{x}-\underline{t}) d\undl{t}
.
\end{equation}
Thus, with the above notation, the Cauchy transform can be viewed as a convolution. Moreover, for every $x = x_0 + \underline{x} \in \mathbb{R}^{1+d}_+$, the Cauchy kernel $\mathscr{C}_{x_0 + \underline{x}}$ belongs to $L^2(\mathbb{R}^d)$, and it satisfies
\[
\|\mathscr{C}_{x}\|_{L^2(\R^d)}^2
\simeq
\|P_{x}\|_{L^2}^2+
\sum_{j=1}^{d}
\underbrace{
\|Q^{(j)}_{x}\|^2_{L^2}
}_{
\|R_j (P_{x})\|_{L^2}^2
}	
\simeq 
\|P_{x}\|_{L^2}^2
\simeq\frac{1}{(x_0)^d},
\]
where $R_j$ denotes the Riesz transform associated with the Clifford basis vector $\mathbf{e}_j$.  Consequently, by the isometry \eqref{prop-isometry-whiteNoise}, the pairing
\begin{equation}\label{Clifford-Cauchy-intgr}
\big\langle S_\omega, \, \overline{\mathscr{C}_{x_0 + \underline{x}}} \big\rangle
\end{equation}
is almost surely well-defined.  Since the Cauchy kernel
$\overline{
	\mathscr{C}_{x_0+\underline{x}}(\underline{t})}=\frac{x_0+\underline{x}-\underline{t}}{|x_0+\underline{x}-\underline{t}|^{1+d}}$ is monogenic for $x=x_0+\underline{x} \in \mathbb{R}^{1+d}_+$ , it is natural to ask whether its stochastic pairing  with white noise \eqref{Clifford-Cauchy-intgr} inherits this property. Following Definition \ref{Def-whiteNoise}, let $S$ denote  the white noise on $\R^d$.  We aim to show that, almost surely in $\Omega$, the random function defined by the pairing
\[
x_0+\underline{x}=:x \mapsto \big\langle S_\omega,\, \overline{\mathscr{C}_{x_0 + \underline{x}}} \big\rangle
\]
is monogenic on $\mathbb{R}^{1+d}_+$. By linearity and the structure of the Cauchy kernel, this reduces to verifying that differentiation commutes with the action of $S_\omega$ on the Poisson and conjugate Poisson kernels. Specifically, it suffices to prove that, almost surely in $\Omega$, the following identities hold for all $1 \leq j,k \leq d$:
\begin{equation}\label{diff_scal-Poiss}
\frac{\partial}{\partial x_0} \left\langle S_\omega ,~ P_{x_0+\underline{x}} \right\rangle
=
\left\langle S_\omega ,~ \frac{\partial}{\partial x_0}P_{x_0+\underline{x}} \right\rangle,
\end{equation}
\begin{equation}\label{diff_scal-ConjPoiss}
\frac{\partial}{\partial x_0} \left\langle S_\omega ,~ Q^{(j)}_{x_0+\underline{x}} \right\rangle
=
\left\langle S_\omega ,~ \frac{\partial}{\partial x_0}Q^{(j)}_{x_0+\underline{x}} \right\rangle,
\end{equation}
\begin{equation}\label{diff_vec-Poiss}
\frac{\partial}{\partial x_k} \left\langle S_\omega ,~ P_{x_0+\underline{x}} \right\rangle
=
\left\langle S_\omega ,~ \frac{\partial}{\partial x_k}P_{x_0+\underline{x}} \right\rangle,
\end{equation}
\begin{equation}\label{diff_vec-ConjPoiss}
\frac{\partial}{\partial x_k} \left\langle S_\omega ,~ Q^{(j)}_{x_0+\underline{x}} \right\rangle
=
\left\langle S_\omega ,~ \frac{\partial}{\partial x_k}Q^{(j)}_{x_0+\underline{x}} \right\rangle.
\end{equation}
We unify the treatment of the Poisson and conjugate Poisson kernels by considering a general class of smooth kernels $K_{x_0 + \underline{x}}(\underline{t})$ and establishing conditions under which differentiation with respect to $x$ can be interchanged with the  pairing $\langle S_\omega, K_{x_0 + \underline{x}} \rangle$.
\subsection{Differentiations with respect to the spatial variables}
We now focus on the conditions under which equations \eqref{diff_vec-Poiss} and \eqref{diff_vec-ConjPoiss} hold. We establish the following result, which guarantees that, for almost every sample path of the white noise, spatial differentiation commutes with the pairing between a deterministic test function and the corresponding distribution.
\begin{Thm}\label{Thm-exchange-noise-spatial-variable}
	Let $S$ be the white noise on $\mathbb{R}^d$, and let $f \in C^2(\mathbb{R}^d)$ satisfy $\partial^\nu f \in L^2(\mathbb{R}^d)$ for  $|\nu| \leq 2$.
	Suppose that for  $1 \leq j \leq d$, there exist $\alpha, \beta > \frac{d}{2}$ such that
	\[
	\|\partial_{jj} f\|_{H^\beta_\alpha(\mathbb{R}^d)}<\infty.
	\]
	Then for a.s. $\omega$, the derivative commutes with the pairing:
	\begin{equation}
	\frac{\partial}{\partial x_{j}} \big\langle S_{\omega}, ~f(\undl{x}- \cdot) \big\rangle
	= \big\langle S_{\omega},~ \frac{\partial}{\partial x_{j}} f(\undl{x}  -\cdot) \big\rangle.
	\end{equation}
\end{Thm}
See Appendix \ref{prf-Thm-exchange-noise-spatial-variable} for a proof.\\
This reduces to estimating the second-order derivatives $\partial_{kk} P_{x_0}$ and $\partial_{kk} Q^{(j)}_{x_0}$, where $P_{x_0}$ and $Q^{(j)}_{x_0}$ are given by \eqref{Poss-kernel} and \eqref{Conj-Poiss-kernel}, respectively. Observe that both kernels contain a factor of order $\langle t \rangle^{-(1+d)}$, while the weighted Sobolev norm  involves a weight $\langle t \rangle^{\alpha}$.This motivates studying $\langle t \rangle^{2a}$ and its higher-order derivatives.
\begin{lem}\label{prop-Polynm-Derivative_Structure}
	(Polynomial-Power Derivative Structure)\\
	Let $$f(t)=\langle t \rangle^{2a}= (1 + |t|^2)^a,$$ where $t \in \mathbb{R}^d$ and $a \in \mathbb{R}$. For every multi-index $\gamma \in \mathbb{N}_0^d$, there exists a real polynomial $P_\gamma(t)$ with $\deg( P_\gamma) \leq |\gamma|$ such that:
	\[
	\partial^\gamma f(t) = P_\gamma(t) \cdot (1 + |t|^2)^{a - |\gamma|}
	=P_\gamma(t) \cdot \langle t \rangle^{2a -2 |\gamma|}.	
	\]
	Consequently, there exists a constant $C_\gamma > 0$ such that
	\[
	|\partial^\gamma f(t)| =
	|P_\gamma(t)| \cdot \langle t \rangle^{2a -2 |\gamma|}
	\leq C_\gamma \, \langle t \rangle^{\,2a - |\gamma|}, ~ \forall\, t \in \mathbb{R}^d.
	\]
\end{lem}
See Appendix \ref{prf-prop-Polynm-Derivative_Structure} for a proof.
\begin{rmk}
	Indeed, Lemma \ref{prop-Polynm-Derivative_Structure} admits a natural extension to the function  $ f(t) = \left( t_0^2 + |t|^2 \right)^a $ with $ |t_0| > 0 $, owing to the  estimate
	\[
	\min\{1, t_0^2\} (1+|t|^2)
	\leq
	(t_0^2+|t|^2)
	\leq \max\{1, t_0^2\}  (1+|t|^2), ~ \forall~ t \in \mathbb{R}^d.
	\]
\end{rmk}
We may now treat the Poisson kernel \(P_{x_0}\) and the conjugate Poisson kernel \(Q^{(j)}_{x_0}\) in a unified way. Both take the form
\begin{equation*}
f(\undl{t})=\frac{{P}(\undl{t})}{((x_0)^2 + |\undl{t} |^2)^{\frac{d+1}{2}}},
\end{equation*}
where the polynomial ${P}$ satisfies $\deg({P}) \leq 1$.  By the pulled-in property of  the weighted Sobolev norm $H^{\alpha}_{\beta}$, we have
\[
\left\|
\partial_{jj} f
\right\|_{H^{\alpha}_{\beta}(\R^d)}
\simeq
\left\| \left\langle\cdot \right\rangle^\beta
\partial_{jj} f
\right\|_{H^{\alpha}(\R^d)}
=
\left\| \left(
\left\langle\cdot \right\rangle^\beta
\partial_{jj} f\right)^\wedge
\right\|_{L^2(\left\langle\cdot \right\rangle^\alpha, \R^d)}.
\]
Thus, for $\alpha, \beta > \frac{d}{2}$, this norm is finite provided the Fourier transform $\big( \langle \cdot \rangle^\beta \, \partial_{jj} f \big)^\wedge$ decays sufficiently fast. To this end, we  introduce the constant-coefficient symbol class.
\begin{Def}\label{Def-SymbolClass}\cite{Stein1993}(Symbol class)
	Let $m \in \mathbb{R}$. A function $\sigma \in C^\infty(\mathbb{R}^{d})$ is said to be a {symbol of class $S^m(\mathbb{R}^d)$} if, for every multi-index $\beta$, there exists a positive constant $C_{\beta}$ such that
	\begin{equation}
	| \partial^\beta \sigma( \xi)| \le
	C_{\beta} (1 + |\xi|)^{m - |\beta| }
	\simeq
	C_{\beta} \langle\xi \rangle^{m - |\beta| }
	.
	\end{equation}
\end{Def}
Clearly, the inclusion $S^m(\mathbb{R}^d) \subseteq S^n(\mathbb{R}^d)$ holds whenever $m \le n$.
\begin{rmk}
	Recall from Lemma \ref{prop-Polynm-Derivative_Structure} that the function $\langle \cdot \rangle^{m}$ belongs to the symbol class $S^m(\mathbb{R}^d)$. Consequently, the Poisson kernel  $P_y \in S^{-(1+d)}(\mathbb{R}^d)$.
\end{rmk}
\begin{prop}\label{prop-Symb-product}\cite{Stein1993}
	Regarding symbols, we have the following basic properties:
	\begin{itemize}
		\item {differentiation:} If $f \in S^m$, then $\partial^\beta f \in S^{m-|\beta|}$.
		
		\item {multiplication:} If $f \in S^m$ and $g \in S^{n}$, then their product $fg$ belongs to $S^{m+n}$.
		
		\item {Embedding:} If $f \in S^m$, then $f$ defines a tempered distribution in $\mathscr{S}^{'}(\mathbb{R}^d)$.
	\end{itemize}
\end{prop}
Next, we investigate the decay estimates for the Fourier transform of symbols in $L^2(\mathbb{R}^d)$. For $f\in L^2(\mathbb{R}^d) \cap S^m(\mathbb{R}^d)\subseteq \mathscr{S}^{'}(\mathbb{R}^d)$, its Fourier transform $\hat{f}$ cannot necessarily be represented as an oscillatory integral. We begin by decomposing $f$ using  a smooth partition of unity $\{\eta_j\}_{j\geq-1}$ (as in Definition \ref{Def-weight-Besov}):
\begin{equation*}
f \stackrel{\mathscr{S}^{'}}{=} \sum_{j\geq-1}\eta_j f .
\end{equation*}
Consequently,
\begin{equation*}
\hat{f} = \left( \sum_{j\geq-1}\eta_j f \right)^{\wedge} 
\stackrel{\mathscr{S}^{'}}{=}
\sum_{j\geq-1} \left( \eta_j {f}\right)^{\wedge}.
\end{equation*}
Since each $\eta_j {f}\in C^\infty_c(\mathbb{R}^d)$, its Fourier transform can be expressed as
\begin{equation}
\left( \eta_j {f}\right)^{\wedge}(\xi) = \int_{\mathbb{R}^d}\eta_j(x) {f}(x) e^{i2\pi x\cdot\xi } dx.
\end{equation}
We have the following result regarding the decay of the Fourier transform for smooth cutoffs of the $S^m$ symbol class on dyadic annuli.
\begin{lem}\cite{Stein1993}\label{Stein1993-dyadic-Estmt}
	Suppose $f$ belongs to the symbol class $S^m(\mathbb{R}^d)$. Let $\{\eta_j\}_{j \ge -1}$ be a Littlewood-Paley dyadic cutoff function satisfying Definition \ref{Def-weight-Besov}. Then
	\begin{equation}
	|  \left( \eta_j {f}\right)^{\wedge}(\xi)| \le A_{M} \cdot \frac{1}{|\xi|^{M}} \cdot 2^{j[d + m - M ]},
	\end{equation}
	for all  $M \ge 0$, where the bound $A_{M}$ is independent of $j \ge -1$.
\end{lem}
Based on the Fourier transform decay estimates for the symbol class on the dyadic annuli, we obtain the Sobolev norm estimates for the $L^2$ symbol class.
\begin{Thm}\label{Thm-SymbolFourierEstmt}
	Suppose $f \in L^2(\mathbb{R}^d) \cap S^m(\mathbb{R}^d)$. Then for any $M \ge 0$ and any $\xi \neq 0$, there exists a constant $C_M > 0$ such that
	\[ 
	|\hat{f}(\xi)| \le C_M \frac{1}{|\xi|^M}.
	\]
	Furthermore, for $\alpha > 0$, it holds that $\|f\|_{H^\alpha(\mathbb{R}^d)} < \infty$.
\end{Thm}
\noindent\textbf{Proof}:~
Since $f \in L^2$, its Fourier transform is not necessarily representable as an oscillatory integral. We first apply a decomposition to $f$ using a smooth partition of unity $\{\eta_j\}_{j \ge -1}$ subordinate to the family of cutoff functions (as in Definition \ref{Def-weight-Besov}):
\[ 
f \stackrel{\mathscr{S}'}{=} \sum_{j \ge -1} \eta_j f.
\]
Then for $\xi \neq 0$, by Lemma  \ref{Stein1993-dyadic-Estmt}, we have
\[ 
|  \left( \eta_j {f}\right)^{\wedge}(\xi)| \le A_{M} \cdot \frac{1}{|\xi|^{M}} \cdot 2^{j[d + m - M ]},
\]
for all  $M \ge 0$. Choosing $M > m+d$, we obtain
\[ 
|\hat{f}(\xi)| \leq
\sum_{j\geq-1} | (\eta_j f)^{\wedge}(\xi)| \le A_{M} \frac{1}{|\xi|^{M}} \sum_{j\geq-1} 2^{j(d + m - M)}
\leq C_{M} \frac{1}{|\xi|^{M}}.
\]
Therefore, for $M > \max\{\alpha + \frac{d}{2}, m+d \}$, it follows that
\begin{equation}
\begin{aligned}
\|f\|^2_{H^\alpha(\mathbb{R}^d)}=&
\|\hat{f}\|^2_{L^2(\left\langle \cdot\right\rangle^\alpha, \mathbb{R}^d )}\\
=&
\int_{\mathbb{R}^d}
|\hat{f}(\xi)|^2 (1+|\xi|^2)^\alpha d\xi\\
=&
\int_{|\xi|\leq1} +\int_{|\xi|>1}\\
\lesssim&
\|f\|^2_{L^2}+\int_{|\xi|>1}|\xi|^{2\alpha-2M}d\xi\\
\leq&
\|f\|^2_{L^2}+\frac{C_M}{2M - 2\alpha - d}.
\end{aligned}
\end{equation}
\hfill\qed\\
From Theorem \ref{Thm-SymbolFourierEstmt}, to estimate $\|\partial_{jj} f\|_{H^\alpha_\beta(\mathbb{R}^d)}\simeq\|\langle\cdot \rangle^\beta \partial_{jj} f\|_{H^\alpha(\mathbb{R}^d)}
$, it suffices to verify whether the amplitude function $t \mapsto \langle t \rangle^\beta \partial_{jj} f(t)$ belongs to the $L^2$  symbol class. 
\begin{lem}(Weighted Sobolev Estimate for the  Second-Order Derivatives) \label{lem-SobolevEstimate-second-derivative}
	Let
	\begin{equation}
	f({t}) = \frac{{P}({t})}{(x_0^2 + |{t}|^2)^{\frac{d+1}{2}}},
	\end{equation}
	where $d \geq 1$ and ${P}$ is a polynomial with $\deg({P}) \leq 1$.
	Then, for $\alpha > \frac{d}{2}$ and $\frac{d}{2} < \beta < \frac{d}{2}+2$, we have 
	\[ 	\langle\cdot \rangle^\beta \partial_{jj} f\in{L^2}\cap{S}^{-(d+2)+\beta}(\mathbb{R}^d) \] 
	and 
	\[
	\|\partial_{jj} f\|_{H^\alpha_\beta(\mathbb{R}^d)} 	
	< \infty.
	\]
\end{lem}
\noindent\textbf{Proof}:~
Since the function $t \mapsto \langle t \rangle^\beta$ is smooth, Lemma \ref{prop-Polynm-Derivative_Structure} implies
\[
\bigl| \partial^{\gamma} \langle t \rangle^\beta \bigr| \lesssim \langle t \rangle^{\beta - |\gamma|},
\]
which indicates that $\langle \cdot \rangle^\beta \in S^\beta(\mathbb{R}^d)$. Recall that $f(t) = \frac{P(t)}{(x_0^2 + |t|^2)^{\frac{d+1}{2}}}$, where $P \in S^{\deg(P)}(\mathbb{R}^d)$ with $\deg(P) = m \leq 1$. Applying the product rule for symbol classes (Proposition \ref{prop-Symb-product}), we deduce that $f \in S^{m-(1+d)}(\mathbb{R}^d)$ and $\partial_{jj} f \in S^{m-d-3}(\mathbb{R}^d)$. Consequently,
\[
\langle \cdot \rangle^\beta \partial_{jj} f \in S^{m-d+\beta-3}(\mathbb{R}^d) \subseteq S^{-(d+2)+\beta}(\mathbb{R}^d).
\]
Furthermore, observe that $\langle \cdot \rangle^\beta \partial_{jj} f \in L^2(\mathbb{R}^d)$ provided that $\beta < \frac{d}{2} + 2$. Therefore, by Theorem \ref{Thm-SymbolFourierEstmt}, for $\alpha > \frac{d}{2}$ and $\frac{d}{2} < \beta < \frac{d}{2} + 2$, we have
\[
\| \partial_{jj} f \|_{H^\alpha_\beta(\mathbb{R}^d)} < \infty.
\] \hfill\qed\\
The equations \eqref{diff_vec-Poiss} and \eqref{diff_vec-ConjPoiss}  now follow directly from Theorem \ref{Thm-exchange-noise-spatial-variable} and Lemma \ref{lem-SobolevEstimate-second-derivative}, namely:
\begin{Thm}\label{Thm-diff_vec-Kern}
	Let $S$ denote the white noise on $\mathbb{R}^d$, and let $P_{x_0}$ and $Q^{(j)}_{x_0}$ be the Poisson  and conjugate Poisson kernels, respectively.
	Then  for any $k \in \{1, \dots, d\}$, the following differentiation-pairing
	commutation relations hold almost surely for $\omega$:
	\[
	\frac{\partial}{\partial x_k} \left\langle S_\omega ,~  P_{x_0+\underline{x}} \right\rangle
	=
	\left\langle S_\omega ,~ \frac{\partial}{\partial x_k}P_{x_0+\underline{x}} \right\rangle,
	\]and
	\[ \frac{\partial}{\partial x_k} \left\langle S_\omega ,~ Q^{(j)}_{x_0+\underline{x}} \right\rangle
	=
	\left\langle S_\omega ,~ \frac{\partial}{\partial x_k}Q^{(j)}_{x_0+\underline{x}} \right\rangle.
	\]
\end{Thm}
\subsection{Differentiation with respect to the time variable}
To justify the commutativity between time derivatives with the white noise functional, we need to prove \eqref{diff_scal-Poiss} and \eqref{diff_scal-ConjPoiss}. Analogous to Theorem \ref{Thm-exchange-noise-spatial-variable},  this reduces to establishing the following two estimates.  For $a = a_0 + \underline{a} \in \mathbb{R}^{1+d}_+$ and $\alpha,~ \beta > \frac{d}{2}$, we claim the following  estimates for the Poisson and conjugate Poisson kernels,
\begin{equation*}
\left\| \left( \frac{\partial}{\partial x_0} \right)^2 P_{x_0+\undl{a}}\bigg|_{x_0=a_0+\varepsilon} \right\|_{H^\alpha_\beta(\R^d)}
\lesssim
(1 + |\underline{a}|)^\beta \cdot \left\|\left( \frac{\partial}{\partial x_0} \right)^2
\bigg|_{x_0=a_0+\varepsilon}
P_{x_0}\right\|_{H^\alpha_\beta(\R^d)}.
\end{equation*}
and, for every $1 \leq j \leq d$,
\begin{equation*}
\left\| \left( \frac{\partial}{\partial x_0} \right)^2 Q^{(j)}_{x_0+\underline{a}}\bigg|_{x_0=a_0+\varepsilon} \right\|_{H^\alpha_\beta(\R^d)}
\lesssim
(1 + |\underline{a}|)^\beta \cdot \left\|\left( \frac{\partial}{\partial x_0} \right)^2
\bigg|_{x_0=a_0+\varepsilon}
Q^{(j)}_{x_0}\right\|_{H^\alpha_\beta(\R^d)}.
\end{equation*}
A direct computation shows that both second-order derivatives $\left( \frac{\partial}{\partial x_0} \right)^2
P_{x_0}(\underline{t})$ and $\left( \frac{\partial}{\partial x_0} \right)^2Q^{(j)}_{x_0}(\underline{t})$ can be expressed in the unified form:
\begin{equation}\label{expression-2nd-diff}
g({t})=\frac{{P}({t})}{((x_0)^2 + |{t} |^2)^{\frac{d+5}{2}}},
\end{equation}
where the polynomial ${P}$ satisfies $\deg({P}) \leq 3$. By Lemma \ref{prop-Polynm-Derivative_Structure} and Proposition \ref{prop-Symb-product}, it follows that $g\in {S}^{\deg({P})-({d+5})}(\mathbb{R}^d)\subseteq {S}^{-({d+2})}(\mathbb{R}^d)$ and $\langle \cdot \rangle^\beta g\in {S}^{\beta-({d+2})}(\mathbb{R}^d)$. Moreover, for $\beta<\frac{d}{2}+2$, we have $\langle \cdot \rangle^\beta g\in L^2(\mathbb{R}^d)$. Consequently, the $H^\alpha_\beta(\mathbb{R}^d)$-norm estimate for $g$ follows the same argument as in the case of spatial derivatives.

Therefore, for  $\alpha > \frac{d}{2}$ and $\frac{d}{2} < \beta < \frac{d}{2}+2$, the $H^\alpha_\beta(\mathbb{R}^d)$ norms of the second-order time derivatives of the form  \eqref{expression-2nd-diff} are finite. This establishes \eqref{diff_scal-Poiss} and \eqref{diff_scal-ConjPoiss}, namely:
\begin{Thm}\label{Thm-diff_scal-Kern}
	Let $S$ denote the white noise on $\mathbb{R}^d$, and let $P_{x_0}$ and $Q^{(j)}_{x_0}$ be the Poisson  and conjugate Poisson kernels, respectively.
	Then  the following differentiation-pairing
	commutation relations hold almost surely for $\omega$:
	\[
	\frac{\partial}{\partial x_0} \left\langle S_\omega ,~ P_{x_0+\underline{x}} \right\rangle
	=
	\left\langle S_\omega ,~ \frac{\partial}{\partial x_0}P_{x_0+\underline{x}} \right\rangle,
	\]
	and
	\[
	\frac{\partial}{\partial x_0} \left\langle S_\omega ,~ Q^{(j)}_{x_0+\underline{x}} \right\rangle
	=
	\left\langle S_\omega ,~ \frac{\partial}{\partial x_0}Q^{(j)}_{x_0+\underline{x}} \right\rangle.
	\]
\end{Thm}
Let $S$ be the white noise on $\mathbb{R}^d$. We  establish that  the Dirac operator ${D}$ commutes with the white noise functional almost surely.
\begin{Thm}\label{Thm-analytic-Cauchy}
	The Cauchy transform of the white noise, defined in \eqref{Clifford-Cauchy-intgr} as $\{x\mapsto\left\langle S_\omega, \overline{
		\mathscr{C}_{x}}
	\right\rangle\}$, is almost surely monogenic for $x=x_0+\underline{x} \in \R^{1+d}_+$.
\end{Thm}
\noindent\textbf{Proof}: Let ${D}_x$ be the Dirac operator from \eqref{def_left-right-Dirac}. From Theorems \ref{Thm-diff_vec-Kern} and  \ref{Thm-diff_scal-Kern},  we obtain
\[
\begin{aligned}
{D}_x
\left\langle S_\omega, \overline{
	\mathscr{C}_{x}}
\right\rangle
=&
\left(\sum_{k=0}^d \mathbf{e}_k
\frac{\partial}{\partial {x_k}}
\right)
\left(
\left\langle S_\omega,
P_{x_0+\undl{x}}
\right\rangle
- \sum_{j=1}^{d}
\left\langle S_\omega, Q^{(j)}_{x_0+\undl{x}}
\right\rangle \textbf{e}_j
\right)
\\
=&
\left\langle S_\omega,  {D}_x P_{x_0+\undl{x}}
\right\rangle
-
\left\langle S_\omega,  {D}_x Q_{x_0+\undl{x}}
\right\rangle
\\
=&
\left\langle S_\omega,  {D}_x \overline{
	\mathscr{C}_{x}}
\right\rangle  \\
=&0,\quad 
\forall~x=x_0+\underline{x}\in \R^{1+d}_+,
\end{aligned}
\]
almost surely for $\omega$.
\hfill\qed\\
It follows from Theorem \ref{Thm-analytic-Cauchy} that the Poisson transform 
\[
u_\omega(y + \underline{x})
:= \left\langle S_\omega, P_{y + \underline{x}} \right\rangle
= \frac{{c_d}}{2} \left(
\left\langle S_\omega, \overline{\mathscr{C}_{y + \underline{x}}} \right\rangle
+
\left\langle S_\omega, \mathscr{C}_{y + \underline{x}} \right\rangle
\right)
\]
is almost surely harmonic in $\mathbb{R}^{1+d}_+$ owing to the decomposition of the Laplacian
$  \Delta=D\overline{D}=
(D_0+\underline{D})(D_0-\underline{D}) $. We may regard
$\left\langle S_\omega, \overline{\mathscr{C}_{y + \underline{x}}} \right\rangle$
as the \textit{canonical monogenic completion} of the  sample path $S_\omega$.

\subsection{Vertical weak convergence of the Poisson transform}
\label{Sec_Radial weak convergence}
This section proves that the white noise convolution with the Poisson kernel converges vertically to the white noise itself in the sense of $L^2(\Omega,d\mathbb{P})$. Note that
\[
\begin{aligned}
\|S\ast P_{y}\|_{L^2(\Omega,L^2(\mathbb{R}^d, \left\langle\cdot \right\rangle^\alpha ) )}^2
=&
\int_{\R^d}\mathbb{E}_\omega|(S_\omega\ast P_{y})(x)|^2
\left\langle{x} \right\rangle^{2\alpha} dx\\
\simeq&
\frac{1}{y^d}
\int_{\R^d}\left\langle{x} \right\rangle^{2\alpha} dx\\
\lesssim&\frac{1}{y^d},  \quad \alpha<-\frac{d}{2}
\end{aligned}
\]
has already established that the function
\[
\left\lbrace
(x,\omega)\mapsto
(S_\omega\ast P_{y})(x)
\right\rbrace
\in L^2(\Omega,L^2(\mathbb{R}^d, \left\langle\cdot \right\rangle^\alpha ) ),~\alpha<-\frac{d}{2}.
\]
Hence, for a.s. $\omega\in\Omega$ and $y>0$, the function
\[ \left\{x\mapsto
(S_\omega\ast P_{y})(x)
\right\} {\in} L^2(\mathbb{R}^d, \left\langle\cdot \right\rangle^\alpha )\subseteq
\mathscr{S}^\prime(\mathbb{R}^d),~\alpha<-\frac{d}{2}. \]
In this case, the action of $S_\omega\ast P_{y}$ on $ \mathscr{S}(\mathbb{R}^d)$ can be expressed by an integral
\[
\left\langle S_\omega\ast P_y, \varphi \right\rangle_{\mathscr{S}^\prime\times\mathscr{S}}
=
\int ( S_\omega\ast P_y) \varphi,\quad
\forall~\varphi\in \mathscr{S}(\mathbb{R}^d).
\]
And for $\alpha>\frac{d}{2}$, the following relations hold:
\begin{equation}
\begin{aligned}
\mathbb{E}
\left|\left\langle S_\omega\ast P_{y}, \varphi \right\rangle_{\mathscr{S}^\prime\times\mathscr{S}}
\right|^2
\leq&
\mathbb{E}_\omega
\| S_\omega\ast P_{y} \|^2_{L^2(\left\langle{\cdot}\right\rangle^{-\alpha} )}
\|\varphi\|^2_{L^2(\left\langle\cdot\right\rangle^{\alpha} )}\\
=&
\int \mathbb{E}_\omega
|\left\langle S_\omega,P_{y+\underline{x}} \right\rangle |^2
\left\langle{x}\right\rangle^{-2\alpha}dx
\|\varphi\|^2_{L^2(\left\langle\cdot\right\rangle^{\alpha} )}\\
\simeq&
\frac{1}{y}\int\left\langle{x}\right\rangle^{-2\alpha}dx
\|\varphi\|^2_{L^2(\left\langle\cdot\right\rangle^{\alpha} )}~(\text{isometric property})\\
\simeq&
\frac{C_\alpha }{y}
\|\varphi\|^2_{L^2(\left\langle\cdot\right\rangle^{\alpha} )}
\end{aligned}
\end{equation}
Before proving the weak convergence, we deal with the change of order of integration.
\begin{lem}\label{Pure-Quadratic-Term}
	Let $S:~L^2(\R^d)\rightarrow L^2(\Omega,d\mathbb{P})$ be the  white noise, and let $ \varphi\in \mathscr{S}(\mathbb{R}^d) $. Then for $y>0$, the function
	\[
	(x,t,\omega)\mapsto
	(S_\omega\ast P_{y})(x)
	\overline{(S_\omega\ast P_{y})(t)}
	\varphi(x) \overline{\varphi(t)}
	\]
	belongs to $L^1(\Omega\times\R^d_x\times\R^d_t)$.
\end{lem}
See Appendix \ref{prf-Pure-Quadratic-Term} for a proof.\\
\begin{lem}\label{Cross Term}
	Let $S:~L^2(\R^d)\rightarrow L^2(\Omega,d\mathbb{P})$ be the  white noise, and let $ \varphi\in \mathscr{S}(\mathbb{R}^d) $. Then for $y>0$, the function
	\[
	(x, \omega)\mapsto
	(S_\omega\ast P_{y})(x)
	\overline{\left\langle S_\omega,  \varphi\right\rangle}
	\varphi(x)
	\]
	belongs to $L^1(\Omega\times\R^d_x)$.
\end{lem}
See Appendix \ref{prf-Cross Term} for a proof.\\
\begin{Thm}\label{Thm-weakly converge}
	The convolution $S_\omega\ast P_y$ weakly converges to the white noise $S$ in $L^2(\Omega,d\mathbb{P})$ as $y\rightarrow0_+$. Specifically, for any $\varphi\in \mathscr{S}(\mathbb{R}^d)$,  there holds
	\begin{equation}
	\lim_{y\rightarrow0_+}
	\mathbb{E}\left|
	\left\langle {S\ast P_y} ,  \varphi\right\rangle -
	\langle {S, \varphi}\rangle
	\right|^2
	=0.
	\end{equation}
\end{Thm}
See Appendix \ref{prf-Thm-weakly converge} for a proof.
\begin{rmk}
	From the monogenic representation perspective, Theorem \ref{Thm-weakly converge} implies that for any $d\geq1 $,
	\[
\frac{{c_d}}{2}
 \left\langle
	\left\langle S, \overline{\mathscr{C}_{y + \underline{x}}} \right\rangle
	+
	\left\langle S, \mathscr{C}_{y + \underline{x}} \right\rangle
	, 
	\varphi \right\rangle
	\xrightarrow[y \to 0^+]{L^2(\Omega, d\mathbb{P})}
	\langle S, \varphi \rangle,
	\quad \forall\, \varphi \in \mathscr{S}(\mathbb{R}^d).
	\]
	Here, $\left\langle S, \overline{\mathscr{C}_{y + \underline{x}}} \right\rangle$ and $\left\langle S, \mathscr{C}_{y + \underline{x}} \right\rangle=-
	\left\langle S,  \overline{ \mathscr{C}_{-y + \underline{x}} }\right\rangle$ are almost surely monogenic for $y + \underline{x}\in\mathbb{R}^{1+d}_+$ and $-y+ \underline{x}\in\mathbb{R}^{1+d}_-$, respectively.
	This can be viewed as a generalization of the \textit{Sokhotski-Plemelj} formula to the case where the density is the white noise.
\end{rmk}
Now we call that $\left\langle S, \overline{\mathscr{C}_{x}} \right\rangle $ is a \textit{ monogenic representation} of the white noise $S_\omega$,  while the  {monogenic jump} 
\[ 
	\left\langle S, \overline{\mathscr{C}_{x}} \right\rangle
	+
	\left\langle S, \mathscr{C}_{x} \right\rangle
	=
	\left\langle S, \overline{\mathscr{C}_{x}} \right\rangle
	-
	\left\langle S,  \overline{ \mathscr{C}_{-\overline{x}} }\right\rangle
 \]
 is called a  \textit{harmonic representation} of it. Since the difference between any two analytic representations of a tempered distribution must be an entire function \cite{Bremermann1965},  the analytic representation of the white noise is likewise not unique. This leads to the following corollary:
\begin{cor}
	For any  $G\in \text{Mon}(\mathbb{R}^{1+d})$, the sum of the Cauchy transform of the white noise and the function $G$, denoted by $F_\omega(x)=\frac{{c_d}}{2}\left\langle S_\omega, \overline{\mathscr{C}_{x}} \right\rangle+G(x)$, remains an monogenic  representation of white noise. Moreover, it satisfies the
	\[
	\left\langle 
	F(x )
	- F(-\bar{x} )
	, ~
	\varphi(x) 
	\right\rangle
	\xrightarrow[x_0 \to 0^+]{L^2(\Omega, d\mathbb{P})}
	\langle S,~
	 \varphi \rangle,
	\quad \forall\, \varphi \in \mathscr{S}(\mathbb{R}^d),
	\]
where $x=x_0+\underline{x}\in\mathbb{R}^{1+d}_+.$ 
Below, we will refer to the monogenic and harmonic representations of white noise i.e.  $\left\langle S_\omega, \overline{\mathscr{C}_{x}} \right\rangle$ and $\left\langle S_\omega, P_x \right\rangle$ as analytic and harmonic white noise, respectively.
\end{cor}
\section{A new weighted Bergman space for the  monogenic representation of the white noise}\label{Sec_A new Bergman spac}
With the monogenic representation of the white noise  established, we now seek its corresponding monogenic function space. We begin by considering the case of one complex variable.  Motivated by the fact that its counterpart in the unit disc belongs to the weighted Bergman space \cite{Yang2024-whiteNoise}, we naturally turn to the  Bergman space on the upper half-plane $\mathbb{C}_+$.  We first recall the definitions of the weighted holomorphic and the harmonic Bergman spaces on  $\mathbb{C}_+$. We will use the following normalized invariant measure $d\nu$ on  $\mathbb{C}_+$.
\begin{equation}\label{measure-UpperPlane}
d\nu(w)=\frac{1}{\pi}\frac{1}{(2y)^2}dxdy=
\frac{1}{i2\pi}\frac{1}{(2\Im(w))^2}dw d\bar{w}.
\end{equation}
For $\alpha>-1$, the weighted Bergman space $B^2_\alpha(\mathbb{C}_+)$ on the upper half-plane is defined as the following.
\begin{Def}\cite{Vasilevski2008-BergmanSpace}
	\label{Def-holomorphic-Bergman}
	A function $u$ belongs to the weighted holomorphic Bergman space $B^2_\alpha(\mathbb{C}_+)$ for $\alpha>-1$ on the upper half-plane if $u$ is holomorphic on $\mathbb{C}_+$ and belongs to $L^2(\mathbb{C}_+,d\nu_\alpha)$
	where
	\begin{equation*}
	d\nu_\alpha(w)=(1+\alpha) (2\Im(w))^{2+\alpha}  d\nu(w)
	= \frac{1+\alpha}{i2\pi}
	(2\Im(w))^{\alpha}
	dwd\bar{w},
	\end{equation*}
	and
	\begin{equation*}
	\|u\|_{B^2_\alpha(\mathbb{C}_+)}=
	\left( \int_{\mathbb{C}_+} |u(w)| d\nu_\alpha(w) \right) ^{1/2}.
	\end{equation*}	
\end{Def}
\begin{Def}\cite{Axler_Harmonic-funct-thrm,Loaiza-2015}
	\label{Def-harm-Bergman}For  $\alpha>-1$, the weighted harmonic Bergman space ${b}^2_\alpha(\mathbb{R}^2_+)$  is the closed subspace of $L^2(\mathbb{R}^2_+, d\nu_\alpha)$ consisting of all complex valued harmonic functions and
	\begin{equation*}
	\|u\|_{b^2_\alpha(\mathbb{R}^2_+)}\simeq
	\left( \int_{\mathbb{R}^2_+} |u(x,y)| y^\alpha dxdy \right) ^{1/2}.
	\end{equation*}
\end{Def}
\subsection{Weighted monogenic/harmonic Bergman space}
The classical weighted Bergman space $B^2_\alpha(\mathbb{R}^2_+)$ can not be adopted to the Poisson transform of the white noise  \eqref{Poisson-WhiteNoise} since its Bochner type norm
\[
\begin{aligned}
\mathbb{E}_\omega
\int_{\mathbb{R}^2_+}
\left|
\left({S}_\omega\ast P_{y} \right) (x)
\right|^2
y^\alpha dxdy
\simeq&
\int_{\mathbb{R}^2_+}
\|P_y\|_{L^2(\mathbb{R})}^2
y^\alpha dxdy\\
\simeq&
\int_0^\infty y^{\alpha-1}dy\int_{\mathbb{R}}dx
\end{aligned}
\]
is infinity. The $B^2_\alpha(\mathbb{R}^2_+)$ weight function $\varrho_\alpha(y)=y^\alpha$ measures the distance of
the point $w$ to the boundary of the upper half-plane, which  can not cancel the boundary singularity arising in the point at infinity $\{y+\underline{x}\in\mathbb{R}^2_+:~|y+\underline{x}|\rightarrow \infty \}$. To deal with this, we need to  construct a new weight function to eliminate  the  singularity. Now we consider the following fractional linear (Mobius) transformation
\begin{equation}\label{MobiusTsfm-HalfPlane}
\begin{aligned}
\Phi:~\mathbb{C}_+ &\longrightarrow\mathbb{D}\\
w &\longmapsto z=\frac{w-i}{1-i w}
\end{aligned}
\end{equation}
maps the upper half-plane onto the unit disk, and its inverse transformation is given by
\begin{equation}\label{inver-MobiusTsfm}
\begin{aligned}
\Phi^{-1}:~  \mathbb{D} &\longrightarrow \mathbb{C}_+\\
z &\longmapsto w=\frac{z+i}{1+i z}.
\end{aligned}
\end{equation}
This gives the relationship between the area measure of  $\mathbb{C}_+ $ and $\mathbb{D} $,
\begin{equation}
d z d \bar{z}
=\left(\frac{\sqrt{2}}{1-i w}\right)^2\left(\frac{\sqrt{2}}{1+i \bar{w}}\right)^2  d w d \bar{w}.
\end{equation}
Then the inner product for $f,~g\in B^2_{\alpha}(\mathbb{D})$ becomes
\begin{equation}\label{Keep-innerProduct}
\left\langle f,  g \right\rangle_{B^2_{\alpha}(\mathbb{D})}
=
\int_{\mathbb{D}} f(z) \overline{g(z)} (1-|z|^2)^\alpha   d z d \bar{z}
=
\int_{\mathbb{C}_+} \mathcal{U}f(w) \overline{\mathcal{U}g(w)} \rho_\alpha(w)   d w d \bar{w},
\end{equation}
under the Mobius transformation $\Phi$, where the  weight function is given by
\begin{equation}\label{calculate-WeightFunct}
\begin{aligned}
\rho_\alpha(w)=& \left( 1- |\Phi(w)|^2 \right)^\alpha \\
=&
\left( \frac{4 \Im(w)}{(1-iw)(1+i\bar{w})} \right)^\alpha\\
=&
\left(
\frac{4y}{(1+y)^2+x^2}
\right)^\alpha
\end{aligned}
\end{equation}
with $w=x+iy$. In fact, substituting the standard Cayley transform $\frac{i-w}{i+w}=-i\Phi(w)$ into \eqref{calculate-WeightFunct} yields the same weight function. Here, the operator $\mathcal{U}$ associated with $\Phi$  is given by
\begin{equation}\label{Op-Unitary}
\mathcal{U}f(w)
=\frac{2}{(1-i w)^2} f [\Phi(w)]
= \frac{2}{(1-i w)^2} f\left(\frac{w-i}{1-i w}\right)
,~
\forall~f \in B^2_{\alpha}(\mathbb{D}).
\end{equation}
Geometrically, this weight function still measures the distance from a point  $ w=x+iy $ in the upper half-plane to its boundary. Here, the upper half-plane is considered with its point at infinity compactified.
\begin{itemize}
	\item
	For any fixed $x\in \mathbb{R}$, as $y\rightarrow 0_+$, we have
	\[
	\rho_\alpha(x+iy)\leq (4y)^\alpha
	\rightarrow 0 .
	\]
	\item
	As $R:=\sqrt{|x|^2+y^2}\rightarrow\infty$, we have
	\[
	\rho_\alpha(x+iy)\leq
	\left(
	\frac{4R}{1+R^2}
	\right)^\alpha\rightarrow0.
	\]	
\end{itemize}
Therefore, a weighted analytic Bergman space $\mathcal{B}^2_\alpha(\mathbb{C}_+)$ is introduced by the  weight function $\rho_\alpha$,
\begin{equation}
\mathcal{B}^2_\alpha(\mathbb{C}_+)
=
\{u\in\text{Hol}(\mathbb{C}_+):~\int_{\mathbb{C}_+}|u(w)|^2\rho_\alpha(w)dwd\bar{w}<\infty \},
\end{equation}
with the inner product
\begin{equation}
\left\langle u,  v \right\rangle_{\mathcal{B}^2_\alpha(\mathbb{C}_+)}=
\left( \int_{\mathbb{C}_+}
u(w) \overline{v(w)}
\rho_\alpha(w)dwd\bar{w} \right) ^{1/2} .
\end{equation}
Since the transformation $\Phi$ is bi-injection, the operator
$ \mathcal{U}: B^2_\alpha(\mathbb{D})\rightarrow
\mathcal{B}^2_\alpha(\mathbb{C}_+) $
defined by \eqref{Op-Unitary}  is unitary by the relation \eqref{Keep-innerProduct} and its inverse
\begin{equation}\label{inver-Unitary}
\begin{aligned}
\mathcal{U}^{-1}:~ \mathcal{B}^2_\alpha(\mathbb{C}_+) &\rightarrow B^2_\alpha(\mathbb{D})\\
F
&\mapsto
\frac{2}{(1+i z)^2} F [\Phi^{-1} (z)  ].
\end{aligned}
\end{equation}
We also can define the harmonic analog of the space $\mathcal{B}^2_\alpha(\mathbb{C}_+)$, which is defined by
\begin{equation}
\mathfrak{b}^2_\alpha(\mathbb{R}^2_+)
=
\{u\in\text{Harm}(\mathbb{R}^2_+):~\int_{\mathbb{R}^2_+}|u(x,y)|^2\rho_\alpha(x+iy)dxdy<\infty \}.
\end{equation}
The weight function  $\rho_\alpha$, induced by the fractional linear transformation $\Phi$, effectively suppresses the singularity arising from the space-time integral of the Cauchy transform of the white noise.
\begin{Thm}
	\label{Thm-Stch-Hol-Bergman}
The Cauchy transform of the  white noise, defined by
\[ 
F_\omega(z)
=
\left\langle S_\omega, \mathscr{C}_z \right\rangle
= \left\langle S_\omega(t),~ \frac{1}{t-z} \right\rangle, \quad z=x+iy \in \mathbb{C}_+,
 \]
belongs to the stochastic weighted analytic Bergman space $ L^2(\Omega, \mathcal{B}^2_\alpha(\mathbb{C}_+)) $ if and only if $\alpha > 1 $.
\end{Thm}
\noindent {\bf Proof}: For a.s. $\omega$, the analyticity of $F_\omega(z) $ follows from the  Theorem \ref{Thm-analytic-Cauchy}.  Then we need compute  the $ L^2(\Omega,\mathcal{B}^2_\alpha(\mathbb{C}_+) ) $ norm of $F_\omega(z)$:
\begin{equation}\label{Stochastic-AnltcBergman-estmt}
\begin{aligned}
&\|F\|_{L^2(\Omega,\mathcal{B}^2_\alpha(\mathbb{C}_+) )}^2\\
\simeq&
\int_{\mathbb{R}^2_+}
\mathbb{E}_\omega
\left(
|\left({S}_\omega\ast P_{y} \right) (x)|^2
+
|\left( {S}_\omega\ast Q_{y}\right) (x)|^2
\right)
\rho_\alpha(x+iy)dxdy,\\
\simeq&
\int_{\mathbb{R}^2_+}
\left(
\left\langle P_{y+\underline{x}}, P_{y+\underline{x}} \right\rangle
+
\left\langle Q_{y+\underline{x}}, Q_{y+\underline{x}} \right\rangle
\right)
\rho_\alpha(x+iy)dxdy,\\
=&
2\int_{\mathbb{R}^2_+}
\underbrace{
	\left\|  P_{y+\underline{x}} \right\|^2_{L^2(\mathbb{R})}
		 }_{
	= P_{2y}(\underline{0})~\simeq~ \frac{1}{y}	
}
\left(
\frac{4y}{(1+y)^2+x^2}
\right)^{\alpha}  dxdy,\\
\simeq&
\int_{0}^{\infty}
\frac{y ^{\alpha-1}}{(1+y)^{2\alpha}}
\int_{\R}
\frac{1}{\left( 1+ \left(  \frac{|x|}{1+y}\right)^2\right)^\alpha}
 dx
dy,\\
\simeq&
\int_{0}^{\infty}
\frac{y ^{\alpha-1}}{(1+y)^{2\alpha-1}}
dy, \\
\leq&\int_{0}^{\infty}
\frac{1}{(1+y)^{\alpha}}
dy,\\
=&\frac{1}{\alpha-1}, \quad \alpha>1,
\end{aligned}
\end{equation}
where $\left\langle Q_{y+\underline{x}}, Q_{y+\underline{x}} \right\rangle=\left\langle P_{y+\underline{x}}, P_{y+\underline{x}} \right\rangle$ follows from the anti-symmetry of the Hilbert transform. Hence  the $ L^2(\Omega,\mathcal{B}^2_\alpha(\mathbb{C}_+) ) $ norm of $F$ is equivalent to  $ \|u \|^2_{L^2(\Omega,  \mathcal{B}^2_\alpha(\mathbb{R}^2_+))} $. \hfill\qed\\
As seen from the calculation in \eqref{Stochastic-AnltcBergman-estmt}, however, the conventional weighted Bergman space defined in Definition \ref{Def-holomorphic-Bergman}, with the  weight function $\rho_\alpha(x+iy)=y^\alpha$, is not suitable for the analytic and harmonic representations of the white noise. This is due to the divergence of the norm:
\[ 
\|F\|_{L^2(\Omega,{B}^2_\alpha(\mathbb{R}^2_+)) )}^2\simeq
\int_{\mathbb{R}^2_+}
\left\|  P_{y+\underline{x}}\right\|^2_{L^2(\mathbb{R})} 
 y^\alpha d\underline{x}dy
\simeq
\int_0^\infty  y^{\alpha-1}
 dy \int_{\mathbb{R}}
 d\underline{x}
 =\infty.
 \] 
A right weight function that does the job is, as
  found in the computation \eqref{Stochastic-AnltcBergman-estmt}, $\rho_\alpha(x+iy)=\left(
  \frac{4y}{(1+y)^2+x^2}
  \right)^\alpha$.
\subsection{Weighted Monogenic Bergman Space}
{Inspired by the weight function derived from conformal mappings in complex analysis of  one  variable, we generalize them to higher dimensions via:}
\[ \rho_\alpha(y+\underline{x}) = \left( \frac{4y}{(1+y)^2 + |\underline{x}|^2} \right)^\alpha,\quad |\underline{x}|^2 = \sum_{j=1}^{d} |x_j|^2. \]
Consequently, we introduce the corresponding weighted Bergman spaces for monogenic and harmonic white noise as follows:
\[ 
\begin{aligned}
\mathcal{B}^2_\alpha(\mathbb{R}^{1+d}_+)
&= \{ u \in \text{Mon}(\mathbb{R}^{1+d}_+) : \int_{\mathbb{R}^{1+d}_+} |u(y+\underline{x})|^2 \rho_\alpha(y+\underline{x}) \, d\underline{x} \, dy < \infty \}, \\
\mathfrak{b}^2_\alpha(\mathbb{R}^{1+d}_+)
&= \{ u \in \text{Harm}(\mathbb{R}^{1+d}_+) : \int_{\mathbb{R}^{1+d}_+} |u(y+\underline{x})|^2 \rho_\alpha(y+\underline{x}) \, d\underline{x} \, dy < \infty \},
\end{aligned}
\]
We obtain the following hypercomplex generalization of Theorem \ref{Thm-Stch-Hol-Bergman}:
 \begin{Thm}
 	\label{Thm-Stch-Mono-Bergman}
The monogenic white noise, defined by
\[ 
F_\omega(x)
=
\left\langle S_\omega, \overline{\mathscr{C}_{x}} \right\rangle
= \left\langle S_\omega(t),~ \frac{x-\underline{t}}{|x-\underline{t}|^{1+d}} \right\rangle, \quad x=x_0+\underline{x} \in \mathbb{R}^{1+d}_+,
\]
belongs to the stochastic weighted monogenic  Bergman space $ L^2(\Omega, \mathcal{B}^2_\alpha(\mathbb{R}^{1+d}_+)) $ if and only if $\alpha > d $.
 \end{Thm}
\noindent {\bf Proof}: This follows similarly to the calculations in \eqref{Stochastic-AnltcBergman-estmt}. The modified range for the index $\alpha$ arises because the integral in \eqref{Stochastic-AnltcBergman-estmt} takes the following higher-dimensional form:
\[ 
\begin{aligned}
&\int_{\mathbb{R}^2_+}
	\left\|  P_{y+\underline{x}} \right\|^2_{L^2(\mathbb{R}^d)}
\left(
\frac{4y}{(1+y)^2+x^2}
\right)^{\alpha}  dxdy \\
\simeq &
\int_{0}^{\infty}
\frac{y ^{\alpha-d}}{(1+y)^{2\alpha}}
\underbrace{\int_{\R^d}
	\left( 1+ \left(  \frac{|x|}{1+y}\right)^2\right)^{-\alpha}
	dx}_{\simeq (1+y)^d,~\alpha>\frac{d}{2}}
dy\\
\simeq&
\int_{0}^{\infty}
\frac{y ^{\alpha-d}}{(1+y)^{2\alpha-d}}dy\\
\leq&\frac{1}{\alpha-1},\quad \alpha>d\geq1. 
\end{aligned}
 \]
\hfill \qed\\ 
\begin{cor}
	Let $\rho$ be a non-negative weight function satisfying an upper half-space integrability condition
	\begin{equation}\label{condition-weight-random}
	\int_{\mathbb{R}^{1+d}_+}
	\frac{\rho(y+\underline{x})}{y^d} \, d\underline{x} \, dy < \infty.
	\end{equation}
	Then, the monogenic white noise $F_\omega(y+\underline{x}) = \left\langle S_\omega, \overline{\mathscr{C}_{y+\underline{x}}} \right\rangle$  belongs to $L^2(\Omega, L^2(\mathbb{R}^{1+d}_+, \rho(y+\underline{x}) \, d\underline{x} \, dy))$.
\end{cor}
\noindent {\bf Proof}: By direct calculation, we have
\[ 
\begin{aligned}
&\|F\|_{L^2(\Omega,  L^2(\mathbb{R}^{1+d}_+,~ \rho(y+\underline{x})d\underline{x}dy))}^2\\
\simeq&
\int_{\mathbb{R}^{1+d}_+}
\mathbb{E}_\omega
|\left\langle S_\omega, \overline{\mathscr{C}_{y+\underline{x}}} \right\rangle|^2
\rho(y+\underline{x})d\underline{x}dy,\\
\simeq&
\int_{\mathbb{R}^{1+d}_+}
\left\|  P_{y+\underline{x}}\right\|^2_{L^2(\mathbb{R}^d)}
\rho(y+\underline{x}) 
d\underline{x} dy\\
\simeq&
\int_{\mathbb{R}^{1+d}_+}
\frac{\rho(y+\underline{x}) }{y^d}
d\underline{x} dy\\
<&\infty.
\end{aligned}
 \]
\hfill \qed\\
Thus, the condition \eqref{condition-weight-random} is sufficient to guarantee that the monogenic white noise belongs to the associated stochastic weighted space $L^2(\Omega, L^2(\mathbb{R}^{1+d}_+, \rho(y+\underline{x}) \, d\underline{x} \, dy))$.
\section{Priori estimates for the harmonic white noise}
\label{Sec_Priori estimates}
In this section, we  obtain new norm estimates comparing the weighted Bergman space norms of the harmonic white noise in the upper half-space with the weighted Sobolev norms of their boundary restrictions.
\subsection{the Case for $\mathbb{D}$: Equivalent  characterization of the weighted Bergman norm }
Let $\{a_k\}_{k\in\mathbb{Z}}$ be a sequence of independent standard complex normal  random variables, i.e. $\{a_k\}_{k} \stackrel{\text{i.i.d.}}{\sim} \text{Norm}_\mathbb{C}(0, 1)$. Let $\mathscr{D}(\partial\mathbb{D})$ denote the locally convex linear topological  space of complex-valued smooth functions on the unit circle, equipped with the seminorms
\[ 
[f]_\alpha=\sup_{x\in\partial \mathbb{D}} \left|\left( \frac{d}{dx}\right) ^\alpha f(x)\right|,\quad \alpha\in \mathbb{N}_0.
\]
Let $\mathscr{D}^{\prime}(\partial\mathbb{D})$ be its dual space. For details, we refer to \cite{Ambrosio2020}.
Using an equivalent definition of the Gaussian white noise on the unit circle, it can be expressed as a Fourier series \cite{Veraar2010-Regularity-GaussianNoise}.
\begin{equation}\label{ortho-expan-WhiteNoise}
S_\omega=\sum_{k\in\mathbb{Z}} a_k(\omega
) e_k
\end{equation}
in 	$\mathscr{D}^{\prime}\left(\partial\mathbb{D}\right)$ with $e_k(\theta)=e^{i k \theta}$, i.e. $ \left\langle  S_\omega, g \right\rangle =\sum_{k\in\mathbb{Z}}a_k \overline{\hat{g}_{k}} $ for $g \in  \mathscr{D}\left(\partial\mathbb{D}\right)$. Recall the related Sobolev spaces.
\begin{Def}\label{def-negative-index-Sobolev}\cite{Veraar2010-Regularity-GaussianNoise}
	For $s \in \mathbb{R}$, the  periodic
	Sobolev space on $\partial\mathbb{D}$ is defined as
	\begin{equation}\label{def-negative-Sobolev-norm-torus}
	H^s(\partial\mathbb{D})=
	\{
	f\in \mathscr{D}^{\prime}(\partial\mathbb{D}):~\|f\|^2_{H^s(\partial\mathbb{D})}
	=
	\sum_{k\in\mathbb{Z}}\left(1+|k|^2\right)^{s}|\hat{f}_k|^2
	<\infty
	\},
	\end{equation}	
	where $\{\hat{f}_k\}_k$ are the Fourier coefficients of $f$.
\end{Def}		
It is well-known that the sample paths of the Gaussian white noise on the unit circle almost surely belong to Sobolev space $H^{-s}(\partial\mathbb{D})$ for any $s>\frac{1}{2}$ \cite{Veraar2010-Regularity-GaussianNoise}. Thus, the Fourier coefficients of the Gaussian white noise almost surely satisfy
\begin{equation}\label{SobolevNorm-whiteNoise-circle}
\sum_{k\in\mathbb{Z}}
\left(1+|k|^2\right)^{-s}|a_k(\omega)|^2<\infty,{\quad} s>\frac{1}{2}.
\end{equation}		
The Poisson transform of the  Gaussian white noise is given by the Abel summation of its Fourier series:
\begin{equation}\label{solution-unitDisk-AbelSum}
u_\omega(z)=\sum_{k\in\mathbb{Z}}
a_k(\omega) r^{|k|} e^{ik\theta},{\quad}z=re^{i\theta} \in \mathbb{D}.
\end{equation}
By \cite{Yang2024-whiteNoise}, for every $\alpha>0$, the Poisson transform \eqref{solution-unitDisk-AbelSum} almost surely belongs to the weighted Bergman space $B^2_\alpha(\mathbb{D})$.
We now estimate the weighted Bergman norm for individual sample path of  \eqref{solution-unitDisk-AbelSum}.  For almost every $\omega$, by the norm equivalence established in \cite{Qu2019-WeightedBergman}, we have
\begin{equation}
\|u_\omega\|^2_{B^2_\alpha(\mathbb{D})}
\simeq
\sum_{k\in\mathbb{Z}} |a_k(\omega)|^2
\left(1+|k|\right)^{-(1+\alpha)}
\simeq
\|S_\omega\|^2_{H^{-\frac{1+\alpha}{2}}(\mathbb{T})}
<\infty,
\end{equation}
for $\alpha>0$. Taking the expectation of above    norms  yields
\begin{equation}
\mathbb{E}
\|u\|^2_{B^2_\alpha(\mathbb{D})}
\simeq
\mathbb{E}
\|S\|^2_{H^{-s}(\mathbb{T})}
\simeq
\sum_{k\in\mathbb{Z}}
\frac{\mathbb{E}_\omega|a_k(\omega)|^2}{(1+|k|)^{(1+\alpha)}}
\simeq
\sum_{k=0}^{\infty}\frac{1}{(1+k)^{2s}}
<\infty,
\end{equation}
with $s=\frac{1+\alpha}{2}>\frac{1}{2}$. 
\subsection{the Case for $\mathbb{R}^{1+d}_+$: the weighted Bergman   estimate}
While $L^2$ (energy)-estimates generally follow from the Plancherel Theorem, the non-separable coupling of space-time variables in our weighted Bergman norm's weight function obstructs such direct approaches.  The main conclusions of this section are as follows:
\begin{Thm}\label{Thm_estmt-Sobolev-Bergman}
	Let $ S $ be the white noise on $\R^d$ and   $P_y$ the Poisson kernel.  Denote the Poisson transform of $S$ by \[ \mathscr{P}(S)(y+\underline{x}):=(S\ast P_y)({x}).   \]
	Suppose that the parameters satisfy $\beta > \frac{d}{2}$ and $\tau>2\beta+d+1$. Then, there exists a constant $C$, depending only on $\alpha$ and $\tau$, such that for almost every  $\omega$, we have
	\[
	\left\| \mathscr{P}(S_\omega) \right\|_{\mathcal{B}^2_{\tau}(\mathbb{R}^{1+d}_+)}
	\leq
	C
	\left\| S_\omega \right\|_{H^{-\beta}_{-\beta}(\mathbb{R}^d)}.
	\]
	Here, the norm of the weighted Bergman space $\mathcal{B}^2_{\tau}$ is induced by the weight function $\rho_{\tau}(y+\underline{x}) =\left(
	\frac{4y}{(1+y)^2+|\underline{x}|^2}
	\right)^\tau$. 
\end{Thm}
Before proving Theorem \ref{Thm_estmt-Sobolev-Bergman} we need some preparations.  We first observe a scaling property for the  Sobolev norms.
\begin{lem}
	For any $\lambda>0$ and given $ g\in
	L^2\left(\left\langle \frac{\cdot}{\lambda}
	\right\rangle^{-\alpha} ,~\R^d\right)$, we have
	\begin{equation}\label{scale-equ-Sobolev}
	\|g\|^2_{L^2\left(\langle \frac{\cdot}{\lambda}\rangle^{-\alpha},~\R^d\right)}\simeq
	\frac{1}{\lambda^d}
	\left\|
	\hat{g}\left( \frac{\cdot}{\lambda}\right) \right\|^2_{H^{-\alpha}(\R^d)}.
	\end{equation}
\end{lem}
\noindent {\bf Proof}: For any $\lambda>0$, it follows that
\begin{equation}
\begin{aligned}
\|g\|^2_{L^2\left(\langle \frac{\cdot}{\lambda}\rangle^{-\alpha} \right)}&
=
\int_{\R^d}
|g(x)|^2 \left(1+\left|\frac{x}{\lambda}\right|^2\right)^{-\alpha}dx,
\\
=&
\lambda^d
\int_{\R^d}
|g(\lambda x)|^2 \left(1+|x|^2\right)^{-\alpha}dx,
\\=&
\lambda^d
\int_{\R^d}
\left|
\frac{1}{\lambda^d}
\left[
\hat{g}\left(\frac{\cdot}{\lambda}\right)
\right]^{\vee} (x)
\right|^2 \left\langle {x}\right\rangle ^{-2\alpha}dx,
\\\simeq&
\frac{1}{\lambda^d}
\left\|
\hat{g}\left( \frac{\cdot}{\lambda}\right) \right\|^2_{H^{-\alpha}(\R^d)}. 
\end{aligned}
\end{equation}
\quad 
\hfill\qed
\subsubsection{Dyadic frequency decomposition of the weighted  $L^2$-norm}
For almost every white noise path $S_\omega \in H^{-\gamma}_{-\mu}(\mathbb{R}^d)$, with $\gamma,~\mu>\frac{d}{2}$,  we have 
\begin{equation}
\label{Norm-Bergm-weight-L2}
\begin{aligned}
\|
\mathscr{P}(S)
\|^2_{\mathcal{B}^2_{\tau}(\mathbb{R}^{1+d}_+)}
\simeq&
\int_{\R^{1+d}_+}
\left|(S_\omega \ast P_y)(x) \right|^2
\frac{y^\tau}
{\left((1+y)^2+|x|^2\right)^\tau}
dx dy,\\
=&
\int_{0}^{\infty} y^\tau
\int_{\R^d}
\left|(S_\omega \ast P_y)(x)\right|^2  \left(
\frac{1}{(1+y)^2+|x|^2}
\right)^\tau dx dy,\\
=&
\int_{0}^{\infty} \frac{y^\tau}{(1+y)^{2\beta}}
\int_{\R^d}
\left| (S_\omega \ast P_y)(x) \right|^2  \left(
\frac{1}{1+ \left( \frac{|x|}{1+y}\right)^2}
\right)^\tau dx dy,\\
=&
\int_{0}^{\infty} \frac{y^\tau}{(1+y)^{2\tau}}
\left\|
S_\omega \ast P_y
\right\|^2_{L^2\left(\langle \frac{\cdot}{1+y}\rangle^{-\tau}\right)}  dy,\quad 
\end{aligned}
\end{equation}
with $ \tau>d $.  The expression in  \eqref{Norm-Bergm-weight-L2} represents  the \textit{weighted time integral of the spatial weighted 
$L^2$(energy)-norm for the harmonic white noise}. 
The non-separable coupling of space-time variables in the weight function prevents  the direct application of the  Plancherel Theorem for estimating the weighted $L^2$ norm. We therefore adopt the Littlewood-Paley decomposition, expressing the norm as the weighted Besov norm to estimate contributions from each dyadic frequency block. By the property of weighted Besov space \cite{Tri06}, we can express the $L^2\left(\langle \frac{\cdot}{1+y}\rangle^{-\tau}\right)$ norm into the  Besov norm as
\begin{equation}
\label{BergmanNorm-High-Low-freqc}
\left\|
S_\omega \ast P_y
\right\|^2_{L^2\left(\langle \frac{\cdot}{1+y}\rangle^{-\tau}\right)}
=
\sum_{j= -1}^{\infty}
\|
\Delta_j
(S_\omega \ast P_y)
\|^2_{L^2\left(\langle \frac{\cdot}{1+y}\rangle^{-\tau}\right)}.
\end{equation}
Inspired by the works  \cite{Muckenhoupt1983,Strichartz1967}, we convert the weighted Lebesgue norm of a convolution into the Fourier-Sobolev norm of a product.  Note that by Theorem \ref{Thm-L2-convolution of white noise}, $S_\omega\ast{P_y}\in\mathscr{S}^{'}(\mathbb{R})$ holds almost surely. Therefore, we have
\[
\left( \Delta_{j}
(S_\omega\ast{P_y})\right)^\wedge=
\left( {\eta^{\vee}_j}\ast(S_\omega\ast{P_y})\right)^\wedge
=\eta_j (S_\omega\ast{P_y})^{\wedge}=\eta_j \widehat{P_y} \widehat{S_\omega}.
\]
Observe that $\widehat{P_y}(\xi)=\exp(-2\pi y|\xi|)$ is not differentiable at $\xi=0$,  whereas it is smooth on $\mathbb{R}\setminus\{0\}$. Therefore, we decompose \eqref{BergmanNorm-High-Low-freqc} into low and high-frequency parts as follows:
\begin{equation}
\begin{aligned}
&\|
\Delta_{-1}
(S_\omega \ast P_y)
\|^2_{L^2\left(\langle \frac{\cdot}{1+y}\rangle^{-\tau}\right)}
+
\sum_{j= -1}^{\infty}
\|
\Delta_j
(S_\omega \ast P_y)
\|^2_{L^2\left(\langle \frac{\cdot}{1+y}\rangle^{-\tau}\right)},
\\
=&
\sum_{|{j^{'}}+1|\leq 1}
\|
\Delta_{-1}
\Delta_{j^{'}}
(S_\omega \ast P_y)
\|^2_{L^2\left(\langle \frac{\cdot}{1+y}\rangle^{-\tau}\right)}\\
+&
\sum_{j=0}^{\infty}
\sum_{|j^{'}-j|\leq 1}
\|
\Delta_j
\Delta_{j^{'}}
(S_\omega \ast P_y)
\|^2_{L^2\left(\langle \frac{\cdot}{1+y}\rangle^{-\tau}\right)},
\end{aligned}
\end{equation}
where $\Delta_j=\Delta_j (\Delta_{j-1}+\Delta_j+\Delta_{j+1})$.
\subsubsection{The low-frequency estimate $ \Delta_{-1} $ }
To estimate the
\begin{equation}
\sum_{{j^{'}}:~|{j^{'}}+1|\leq 1}
\|
\Delta_{-1}
\Delta_{j^{'}}
(S_\omega \ast P_y)
\|^2_{L^2\left(\langle \frac{\cdot}{1+y}\rangle^{-\tau}\right)}
\end{equation}
with $ \tau>\frac{d}{2}$, we rewrite the above expression in terms of the Fourier-Sobolev norm
\begin{equation}
\begin{aligned}
\left\|
\Delta_{-1}
\Delta_{j^{'}}
(S_\omega \ast P_y)
\right\|
^2_{L^2\left(\langle \frac{\cdot}{\lambda}\rangle^{-\tau}\right)}
=
\frac{1}{\lambda^d}
\left\|
\left(\eta_{-1}  \widehat{P_y}\eta_{j^{'}} \widehat{S_\omega} \right)
\left(\frac{\cdot}{\lambda}\right)
\right\|
^2_{H^{-\tau}(\R^d)},
\end{aligned}
\end{equation}
where the scaling factor $\lambda:=1+y>1$ and $\widehat{P_y}(\xi)=\exp(-2\pi y|\xi|)$. The  scaling factor $\frac{1}{(1+y)^d}$ will also provide a  decay term for the subsequent time integration. Before proceeding, we recall the following crucial estimate for Sobolev norms of products as the following.
\begin{prop}\label{SobolevEstimate-product}
	\cite{Gavage2007}
	For all $s$ and $t$ with $s+t>0$, if $u \in H^s(\R^d)$ and $v \in H^t(\R^d)$ then the product $uv$ belongs to $H^r(\R^d)$ for all $r \leq \min (s, t)$ and $r<s+t-d / 2$. Furthermore, there exists a constant $C$ (depending only on $r, s, t$ and $d$ ) such that
	\begin{equation}
	\|u v\|_{H^r(\R^d)} \leq C\|u\|_{H^s(\R^d)}\|v\|_{H^t(\R^d)}.
	\end{equation}
\end{prop}
Since $\eta_{j}$ belong to $C_c^\infty(\R^d)$, they belong to the Sobolev space $H^s$ for  arbitrary $s>0$. Then we  need  a  {more precise} Sobolev norm estimate of  $\widehat{P_y}$.
\begin{lem}\label{lem-Poiss-weight-L2}
	For any $y>0$ and $r<\frac{d+2}{2}$, there holds
	\begin{equation}
	\|\widehat{P_y}\|^2_{H^{r}(\R^d)}
	\lesssim C_r
	\begin{cases}
	\frac{1}{y^{d-2r}},~&y>1,\\
	\frac{1}{y^{d}},~&0<y<1,\\
	\end{cases}
	\end{equation}
\end{lem}
See Appendix \ref{prf-lem-Poiss-weight-L2} for a proof.\\
For the rescaled Fourier transform of the Poisson kernel
\[
\widehat{P_y}
\left(\frac{\xi}{\lambda}\right)
=\exp\left(-2\pi \frac{y}{\lambda}|\xi|\right)
=\widehat{P_{\frac{y}{\lambda} }}(\xi)
\]
with $ \lambda=1+y $, the following estimate holds when  $r<\frac{d+2}{2}$,
\begin{equation}
\begin{aligned}
\left\|
\widehat{P_y} \left(\frac{\cdot}{\lambda}\right)
\right\|^2_{H^{r}(\R^d)}
=&
\| \widehat{P_{\frac{y}{\lambda} }}
\|^2_{H^{r}(\R^d)} \\
\leq & C_r
\begin{cases}
\left(\frac{y}{1+y}\right)^{-d},~&0<y<1,\\
1,~&y>1,
\end{cases}\\
\leq& C_r
\left( \frac{1}{y^d}+1\right) .
\end{aligned}
\end{equation}
\begin{cor}\label{lem-cutoff-sobolev}
	Let $m_{-1} := \eta_{-1} \widehat{P_y}$, where $\eta_{-1} \in C^\infty_c(\mathbb{R}^d)$ is the  cut-off function and $P_y$ denotes the Poisson kernel. We set $\lambda = 1+y$ with $y > 0$. For  parameters  $\tau >\frac{d+2}{2} > \beta > 0$, the multiplication operator $\widetilde{m_{-1}}\left( \frac{\cdot}{\lambda} \right)$ induced by  $m_{-1}\left( \frac{\cdot}{\lambda} \right)$ is bounded from $H^\tau(\R^d)$ to $H^\beta(\R^d)$, satisfying
	\[  
	\left\|
	\widetilde{m_{-1}}\left( \frac{\cdot}{\lambda} \right)
	\right\|^2_{H^\tau(\R^d){\rightarrow} H^\beta(\R^d)}
	\lesssim
	\left\|  
	m_{-1}\left( \frac{\cdot}{\lambda} \right)
	\right\|^2_{H^r(\R^d)} 
	\lesssim 
	C_r
	\left(y^{d}+y^{-d} \right) .
	\]
	with $r\in(0,\frac{d+2}{2} )$. 
\end{cor}
See Appendix \ref{prf-lem-cutoff-sobolev} for a proof.\\
Using the duality principle, then we have
\begin{equation}\label{}
\begin{aligned}
\left\|
m_{-1}
\left(\frac{\cdot}{\lambda}\right)  ~g
\right\|
^2_{H^{-\tau}}
\lesssim&
\left(y^{d}+y^{-d} \right)
\left\| g
\right\|^2_{H^{-\beta}},~
\forall~ g\in H^{-\beta}(\R^d)
\end{aligned}
\end{equation}
with $\tau>\frac{d+2}{2}>\beta>0$. 
Now we give  the low-frequency estimate at any fixed time.
\begin{prop}
{For $\tau >\frac{d+2}{2}> \beta > \frac{d}{2}$}, and the white noise sample path $S_\omega \in H^{-\alpha^{'}}_{-\beta^{'}}(\mathbb{R}^d)$ with $\beta^{'} > \frac{d}{2}$ and set $\lambda=1+y$, we have the low-frequency estimate at the fixed time $y>0$:
	\begin{equation}\label{Low-frequ-Estmt-fixed-time}
	\left\|
	\Delta_{-1}
	\Delta_{j^{'}}
	(S_\omega \ast P_y)
	\right\|
	^2_{L^2\left(\langle \frac{\cdot}{\lambda}\rangle^{-\tau}\right)}
	\lesssim
	\frac{(1+y)^{2(\beta+d)}}{y^d}
	\|\Delta_{j^{'}} S_\omega\|^2_{L^2
		\left(
		\langle
		\cdot
		\rangle
		^{-\beta}
		\right)},
	\end{equation}
	where $j^{'}\geq-1$ and satisfies $|j^{'}+1|\leq1$.
\end{prop}
\noindent {\bf Proof}: For $\tau>\frac{d+2}{2}>\beta>\frac{d}{2}$, we have
\begin{equation*}
\begin{aligned}
\left\|
\Delta_{-1}
\Delta_{j^{'}}
(S_\omega \ast P_y)
\right\|
^2_{L^2\left(\langle \frac{\cdot}{\lambda}\rangle^{-\tau}\right)}
=&
\frac{1}{\lambda^d}
\left\|
\left(m_{-1}  \eta_{j^{'}} \widehat{S_\omega} \right)
\left(\frac{\cdot}{\lambda}\right)
\right\|
^2_{H^{-\tau}(\R^d)},\\
\lesssim&
\frac{1}{\lambda^d}
\left(y^d+\frac{1}{y^d} \right) \lambda^d
\|\Delta_{j^{'}} S_\omega\|^2_{L^2
	\left(\langle \frac{\cdot}{\lambda}\rangle
	^{-\beta}\right)},\\
=&
\frac{1+y^{2d}}{y^d}
\|\Delta_{j^{'}} S_\omega\|^2_{L^2
	\left(\langle \frac{\cdot}{\lambda}\rangle
	^{-\beta}\right)},\\
\simeq&
\frac{(1+y)^{2d}}{y^d}
\|\Delta_{j^{'}} S_\omega\|^2_{L^2
	\left(\langle \frac{\cdot}{\lambda}\rangle
	^{-\beta}\right)},\\
\lesssim&
\frac{(1+y)^{2(\beta+d)}}{y^d}
\|\Delta_{j^{'}} S_\omega\|^2_{L^2
	\left(\langle
	\cdot
	\rangle
	^{-\beta}\right)}, 
\end{aligned}
\end{equation*}
where $\|\Delta_{j^{'}} S_\omega\|^2_{L^2
	\left(\langle \frac{\cdot}{\lambda}\rangle
	^{-\beta}\right)}
\leq \lambda^{2\beta}
\|\Delta_{j^{'}} S_\omega\|^2_{L^2
	\left(\langle
	\cdot
	\rangle
	^{-\beta}\right)}
$ for $\lambda>1$.
\hfill\qed \\
\subsubsection{An estimate for the high-frequency  part  }
To estimate the
\begin{equation}
\sum_{j^{'}:~|j^{'}-j|\leq 1}
\|
\Delta_{j}
\Delta_{j^{' }}
(S_\omega \ast P_y)
\|^2_{L^2\left(\langle \frac{\cdot}{1+y}\rangle^{-\beta}\right)}
\end{equation}
for $ \beta>\frac{d}{2} $, we only need to study
\begin{equation}
\begin{aligned}
\left\|
\Delta_{j}
\Delta_{j^{'}}
(S_\omega \ast P_y)
\right\|
^2_{L^2\left(\langle \frac{\cdot}{\lambda}\rangle^{-\beta}\right)}
\leq&
\lambda^{2\beta}
\left\|
\Delta_{j}
\Delta_{j^{' }}
(S_\omega \ast P_y)
\right\|
^2_{L^2\left(\langle \cdot\rangle^{-\beta}\right)},\\
=&
\lambda^{2\beta}
\left\|
\eta_{j}  \widehat{P_y}\eta_{j^{' }} \widehat{S_\omega}
\right\|
^2_{H^{-\beta}(\R^d)},
\end{aligned}
\end{equation}
with $\lambda=1+y$. Denote by $m_j=\eta_{j}  \widehat{P_y}$ for $j\geq0$. In contrast to low-frequency estimate, the multiplier $m_j $ exhibits no singularity in high-frequency dyadic blocks away from the origin and possesses arbitrary order differentiability.

We begin by investigating the estimates for higher-order partial derivatives of the Fourier transform $\widehat{P_y}$ in the high-frequency part.
\begin{prop}\label{prop-high-dertv-exp}
For any $M > 0$ and $\alpha \in \mathbb{N}^d$, there is a constant $C(\alpha, d)>0$ independent of $\xi$ such that for all $|\xi| \geq M$,
	\begin{equation} 
	\left|\partial^\alpha e^{-|\xi|}\right|
	\leq
	C(\alpha, d)
	e^{-|\xi|}.
	\end{equation}	
\end{prop}
See Appendix \ref{prf-prop-high-dertv-exp} for a proof.
\begin{cor}\label{Cor-Estmt-exp-funct}
For any $M > 0$ and $|\xi| \geq M$, the higher-order partial derivatives of the scaled exponential function $e^{-y|\xi|}$ satisfy the following estimate for all $y > 0$:
	\begin{equation}
	\forall~\alpha \in \mathbb{N}^d, \quad 
	|(\partial_\xi)^\alpha e^{-y|\xi|}|
	\lesssim
	y^{|\alpha|}e^{-y|\xi|}.
	\end{equation}	
\end{cor}
Given the higher-order partial derivatives of $\widehat{P_y}$ at high frequencies, we can obtain the corresponding Fourier-Sobolev norm estimates.
\begin{lem}\label{lem-high-freq-multpl}
	For $\beta>\frac{d}{2}$ and $j\geq0$,  the multiplication operator $\widetilde{m_j}$ induced by 
	$m_j =\eta_{j}  \widehat{P_y} $
	is bounded from $H^\beta(\R^d)$ to $H^\beta(\R^d)$ with
	\[ 
	\|\widetilde{m_j}
	\|^2_{H^{K}\rightarrow H^{K}}
	\lesssim
	\left\|
	m_j
	\right\|^2_{H^{K}}
	\lesssim
	\frac{2^{j(d-1)}}{y}
	\sum_{|k|\leq{K}}
	\frac{1}{2^{2j|k|}}
	\sum_{l\leq{k}}
	(2^j y)^{2|l|}
	e^{-2\pi y 2^j}
	,
	\]
	for $K\in \mathbb{N}_+$ with $K>[\beta]+1 $.
\end{lem}
See Appendix \ref{prf-lem-high-freq-multpl} for a proof.\\
By the duality, for $\beta>\frac{d}{2}$ we have
\begin{equation}\label{High-frequ-Estmt-fixed-time}
\begin{aligned}
\left\|
\Delta_{j}
\Delta_{j^{'}}
(S_\omega \ast P_y)
\right\|
^2_{L^2\left(\langle \frac{\cdot}{\lambda}\rangle^{-\beta}\right)}
\leq&
\lambda^{2\beta}
\left\|
m_j \eta_{j^{'}} \widehat{S_\omega}
\right\|
^2_{H^{-\beta}(\R^d)}\\
\lesssim&
\lambda^{2\beta}
\left\|
m_j
\right\|
^2_{H^{K}(\R^d)}
\| \Delta_{j^{'}} S_\omega \|^2_{L^2\left(\langle \cdot\rangle^{-\beta}\right)}
\end{aligned}
\end{equation}
with $\lambda=1+y>1$.
We call \eqref{High-frequ-Estmt-fixed-time} the high-frequency estimate at any fixed time $y>0$.
\subsubsection{Combining Estimates for High and Low Frequency Parts}
The proof of Theorem \ref{Thm_estmt-Sobolev-Bergman} now follows by combining the high- and low-frequency estimates.\\
 \noindent {\bf Proof  of  Theorem} \ref{Thm_estmt-Sobolev-Bergman}:~
For $\tau>\frac{d+2}{2}>\beta>\frac{d}{2}$, we integrate  the low-frequency estimate \eqref{Low-frequ-Estmt-fixed-time}  as
\begin{equation}\label{Low-frequ-Est-with-time-integral}
\begin{aligned}
&\int_{0}^{\infty}
\|
\Delta_{-1}
\Delta_{j^{'}}
(S_\omega \ast P_y)
\|^2_{L^2\left(\langle \frac{\cdot}{1+y}\rangle^{-\tau}\right)}
\frac{y^\tau}{(1+y)^{2\tau}}
dy\\
\lesssim&
\int_{0}^{\infty} \frac{y^{\tau-d}}{(1+y)^{2(\tau-\beta-d)}}
dy
\|\Delta_{j^{'}} S_\omega\|^2_{L^2
	\left(\langle
	\cdot
	\rangle
	^{-\beta}\right)},\\
\lesssim&
\|\Delta_{j^{'}}S_\omega\|^2_{L^2
	\left(\langle
	\cdot
	\rangle
	^{-\beta}\right)},\\
\lesssim&
2^{-2j^{'}\beta}
\|\Delta_{j^{'}} S_\omega\|^2_{L^2
	\left(\langle
	\cdot
	\rangle
	^{-\beta}\right)},
\quad |j^{'}+1|\leq 1,
\end{aligned}
\end{equation}
if $\tau>2\beta+d+1$. 

For  $\tau>\beta>\frac{d}{2}$, integrating the high-frequency estimate \eqref{High-frequ-Estmt-fixed-time} with respect to $  \frac{y^\tau}{(1+y)^{2\tau}}dy$ yields
\begin{equation}\label{High-frequ-Estmt}
\begin{aligned}
&\int_{0}^{\infty}
\left\|
\Delta_{j}
\Delta_{j^{'}}
(S_\omega \ast P_y)
\right\|
^2_{L^2\left(\langle \frac{\cdot}{\lambda}\rangle^{-\tau}\right)}
\frac{y^\tau}{(1+y)^{2\tau}}dy\\
\lesssim&
\int_{0}^{\infty}
\frac{y^\tau \lambda^{2\tau}}{(1+y)^{2\tau}}
\left\|
m_j
\right\|
^2_{H^{K}(\R^d)}
\underbrace{
	\| \Delta_{j^{'}} S_\omega \|^2_{L^2\left(\langle \cdot\rangle^{-\tau}\right)}
}_{
\leq~\| \Delta_{j^{'}} S_\omega \|^2_{L^2\left(\langle \cdot\rangle^{-\beta}\right)},~\tau>\beta
}
dy,
\\
\lesssim&
\int_{0}^{\infty}
y^\tau
\sum_{|k|\leq{K}}
\frac{2^{j(d-1)}}{2^{2j|k|}}
\sum_{l\leq{k}}
(2^j y)^{2|l|}
e^{-2\pi y 2^j}
\frac{dy}{y}~
\| \Delta_{j^{'}} S_\omega \|^2_{L^2\left(\langle \cdot\rangle^{-\beta}\right)},\\
=&
\sum_{|k|\leq{K}}
\frac{2^{j(d-1)}}{2^{2j|k|}}
\sum_{l\leq{k}}
\int_{0}^{\infty}
(2^j y)^{2|l|}
e^{-2\pi y 2^j} y^\tau
\frac{dy}{y}
~
\| \Delta_{j^{'}} S_\omega\|^2_{L^2\left(\langle \cdot\rangle^{-\beta}\right)},\\
\lesssim&
\sum_{|k|\leq{K}}
\frac{2^{j(d-1)}}{2^{2j|k|}}
\sum_{l\leq{k}}
\frac{1}{2^{j\tau}}
\Gamma(2|l|+\tau-1)
\| \Delta_{j^{'}} S_\omega\|^2_{L^2\left(\langle \cdot\rangle^{-\beta}\right)},\\
\lesssim&
2^{-j(\tau-d+1)} \| \Delta_{j^{'}} S_\omega\|^2_{L^2\left(\langle \cdot\rangle^{-\beta}\right)},\\
\lesssim&
2^{-2j^{'}\beta} \| \Delta_{j^{'}} S_\omega\|^2_{L^2\left(\langle \cdot\rangle^{-\beta}\right)},
\end{aligned}
\end{equation}
if we choose $\tau>2\beta+d-1$.

To adapt to the  index $\tau$ such that  $\tau>2\beta+d+1$ in the low-frequency estimate  \eqref{Low-frequ-Est-with-time-integral}, we derive from the \eqref{High-frequ-Estmt}  that for any $j\geq-1$: 
\begin{equation}\label{high-frequ-estmt-correction}
\int_{0}^{\infty}
\left\|
\Delta_{j}
\Delta_{j^{'}}
(S_\omega \ast P_y)
\right\|
^2_{L^2\left(\langle \frac{\cdot}{\lambda}\rangle^{-\tau}\right)}
\frac{y^\tau}{(1+y)^{2\tau}}dy
\lesssim
2^{-2j^{'}\beta} \| \Delta_{j^{'}} S_\omega\|^2_{L^2\left(\langle \cdot\rangle^{-\beta}\right)}.
\end{equation}
Combining the low-frequency estimate \eqref{Low-frequ-Est-with-time-integral} with the high-frequency estimate \eqref{high-frequ-estmt-correction}, the weighted Bergman norm $\|\mathscr{P}(S_\omega) \|^2_{\mathcal{B}^2_{\tau}(\R^{1+d}_+)}$ in \eqref{Norm-Bergm-weight-L2} is bounded by
\begin{equation}
\begin{aligned}
\sum_{j= -1}^{\infty}
\sum_{|{j^{'}}-j|\leq 1}
2^{-2j^{'}\beta} 
\| \Delta_{j^{'}} S_\omega\|^2_{L^2\left(\langle \cdot\rangle^{-\beta}\right)}
\simeq&
\sum_{j= -1}^{\infty}
2^{-2j\beta} \| \Delta_{j} S_\omega\|^2_{L^2\left(\langle \cdot\rangle^{-\beta}\right)}\\
\simeq&
\|S_\omega\|^2_{H^{-\beta}_{-\beta}(\R^d)}
\end{aligned}
\end{equation}
with $2^{-j\beta}\simeq2^{-j^{'}\beta}$ when   $|j^{'}-j|\leq1$. \hfill\qed

In fact, it follows from \eqref{Norm-Bergm-weight-L2} that the conclusion of Theorem \ref{Thm_estmt-Sobolev-Bergman} remains valid for any weight function $\rho_\tau$ satisfying the decay condition:
\begin{equation}\label{eq:decay-condition}
0<\rho_\tau(y+\underline{x})\lesssim
\left(
\frac{y}{(1+y)^2+|\underline{x}|^2}
\right)^\tau,
\end{equation}
for  $\tau>2\beta+d+1$.  The upper bound in \eqref{eq:decay-condition} can be replaced by $\left( \frac{y}{c+y^2+|\underline{x}|^2} \right)^\tau$ for any $c>0$, since $c+y^2\simeq(\sqrt{c}+y)^2\simeq(1+y)^2$. A typical example is $\rho_\tau(y+\underline{x}):=y^\tau e^{-(y^2+|\underline{x}|^2)}$.
Furthermore, the decay condition \eqref{eq:decay-condition} yields the upper half-space integrability  condition \eqref{condition-weight-random}. 

\section{Conclusions}
We present harmonic and monogenic representations for  white noise in the Euclidean  spaces through establishing the high-dimensional Sokhotski-Plemelj formula. The monogenicity of the associated random field is  established by  using harmonic analysis techniques and the properties of weighted Sobolev spaces. As an extension of H. G. Feichtinger's white noise theory, this approach bridges monogenic functions and the white noise,  thereby links hypercomplex and stochastic analysis. It also offers guidance for constructing monogenic representations of general stationary stochastic processes. Furthermore, we construct a new class of Bergman spaces and derive the relationship between the weighted Bergman norms in the upper half-spaces and the weighted Sobolev norms of the  white noise.
\appendix
\section{Proof of theorem \ref{Thm-L2-convolution of white noise}}\label{prf-L2-convolution}
By the Definition, we have
\[
\begin{aligned}
\mathbb{E} \int_{\mathbb{R}^d} \frac{\left| (S \star f)(x) \right|^2}{(1+|x|)^{2d}} dx
=
\int_{\mathbb{R}^d}
\frac{\mathbb{E}_\omega\left| (S_\omega \star f)(x) \right|^2}{(1+|x|)^{2d}}
dx
\lesssim\| f \|^2_{L^2(\R^d)}.
\end{aligned}
\]
Hence, for almost every \( \omega \in \Omega \), the function \( S_\omega \star f \)
belongs to \( L^2({\mathbb {R}^d},\langle\cdot\rangle^{-d}) \). We now estimate the $L^2(\Omega)$ norm of
\[ \langle S_\omega \star f, \varphi\rangle
-
\langle S_\omega \ast f, \varphi\rangle. \] Take a sequence $\{f_k\}_k\subset\mathscr{S}(\R^d)$ such that to $f_k\to f$ in $L^2(\mathbb{R}^d)$. Then by the Minkowski inequality, we have
\begin{equation}\label{Random mean distance}
\begin{aligned}
&\|
\langle S \star f, \varphi\rangle-\langle S \ast f_k, \varphi\rangle
\|_{L^2(\Omega)}\\
=&
\left(
\mathbb{E}_\omega
\left|
\int_{\R^d}
\left( S_\omega \star f-S_\omega \ast f_k\right) \varphi
\right|^2
\right)^{\frac{1}{2}},\\
\leq&
\int_{\R^d}
\| (S \star f)(x)- (S \ast f_k)(x)\|_{L^2(\Omega)} |\varphi(x)| dx,\\
=&
\int_{\R^d}
\| \langle S , {\Rsh_x} f^{\checkmark}\rangle
-
\langle S , {\Rsh_x} f_k^{\checkmark}\rangle
\|_{L^2(\Omega)} |\varphi(x)| dx,\\
\simeq&
\int_{\R^d}
\|  {\Rsh_x} f^{\checkmark}- {\Rsh_x} f_k^{\checkmark}
\|_{L^2(\R^d)} |\varphi(x)| dx,\\
=&
\|  f-f_k \|_{L^2(\R^d)}
\int_{\R^d} |\varphi(x)| dx
\xrightarrow[k\to\infty]{}0.
\end{aligned}
\end{equation}
On the other hand, $
\langle S \ast f_k, \varphi\rangle
\stackrel{L^2(\Omega)}{\longrightarrow}
\langle S \ast f, \varphi\rangle$ as $k\to\infty$. Therefore, for almost every $\omega\in\Omega$,
\[
\langle S_\omega \ast f, \varphi\rangle
=
\langle S_\omega \star f, \varphi\rangle
=
\int
\left\langle S_\omega, {\Rsh_x} f^{\checkmark} \right\rangle \varphi(x)dx
. 
\]
 \hfill\qed
\section{Proof of theorem \ref{Thm-exchange-noise-spatial-variable}}\label{prf-Thm-exchange-noise-spatial-variable}
For  a.s. $\omega$ and fixed $\undl{a}\in\R^d$, for $0<\varepsilon<1$ there exist $\theta\in(0,1)$ such that
\begin{equation}\label{DiffQuotient-derivative}
\begin{aligned}
&\Big\langle S_{\omega},~
\frac{f(\undl{a}+\varepsilon \vec{e_j} - \cdot) - f(\undl{a} - \cdot)}{\varepsilon}
- \frac{\partial}{\partial x_j}\bigg|_{\undl{x}=\undl{a}} f(\undl{x} - \cdot)
\Big\rangle\\
=&
\frac{\varepsilon}{2}
\Big\langle S_{\omega},~
\left(\frac{\partial}{\partial x_j}\right)^2\bigg|_{\undl{x}=\undl{a}+\theta\varepsilon\vec{e_j}}
f(\undl{x} - \cdot)
\Big\rangle
\end{aligned}
\end{equation}
by the second-order Taylor expansion. The sample paths of the white noise $ S_\omega $ almost surely belong to the weighted Sobolev space $  H^{-\alpha}_{-\beta}(\mathbb{R}^d) $
for any $\alpha,~\beta>\frac{d}{2}$. We need to estimate the  $  H^{\alpha}_{\beta}(\mathbb{R}^d) $ norm of the second-order partial derivative as
\[
\left(\frac{\partial}{\partial x_j}\right)^2\bigg|_{\undl{x}=\undl{a}+\theta\varepsilon\vec{e_j}}
f(\undl{x} - \undl{t})
=
\partial_{jj} f(\undl{a}+\theta\varepsilon\vec{e_j}-\undl{t}).
\]
Thus we have
\[
\begin{aligned}
\left\|
\left(\frac{\partial}{\partial x_j}\right)^2\bigg|_{\undl{x}=\undl{a}+\theta\varepsilon\vec{e_j}}
f(\undl{x} - \cdot)
\right\|_{H^{\alpha}_{\beta}}
=&
\left\|
\partial_{jj} f(\undl{a}+\theta\varepsilon\vec{e_j}-\cdot)
\right\|_{H^{\alpha}_{\beta}}\\
\lesssim&
(1+|\undl{a}+\theta\varepsilon\vec{e_j}|)^\beta
\left\|
\partial_{jj} f
\right\|_{H^{\alpha}_{\beta}}\\
\leq&
(2+|\undl{a}|)^\beta
\left\|
\partial_{jj} f
\right\|_{H^{\alpha}_{\beta}}.
\end{aligned}
\]
Therefore, when
$ \left\|
\partial_{jj} f
\right\|_{H^{\alpha}_{\beta}(\R^d)}
<\infty $
holds for some $\alpha,~\beta>\frac{d}{2}$,   the expression
\eqref{DiffQuotient-derivative} is an  $H^{\alpha}_{\beta}{\times}H^{-\alpha}_{-\beta}$-type pairing. Hence, it is bounded  by
\[
(2+|\undl{a}|)^\beta \varepsilon
\left\|
\partial_{jj} f
\right\|_{H^{\alpha}_{\beta}}
\]
from the duality inequality of $H^{\alpha}_{\beta}{\times}H^{-\alpha}_{-\beta}$, which converges to $0$
as $\varepsilon\rightarrow0$.  \hfill\qed
\section{Proof of Lemma \ref{prop-Polynm-Derivative_Structure}}\label{prf-prop-Polynm-Derivative_Structure}
We proceed by induction on $k = |\gamma|$.

\textbf{Base case ($k = 0$).}
When $\gamma = 0$, we have $\partial^0 f = f = (1 + |t|^2)^a$.
Set $P_0(t) \equiv 1$, which is a polynomial of degree $0$. Then
\[
\partial^0 f(t) = P_0(t) \cdot (1 + |t|^2)^{a - 0},
\]
so the claim holds for $k = 0$.

\textbf{Inductive hypothesis.}
Assume that for some $k \geq 0$, the statement holds for all multi-indices $\iota$ with $|\iota| = k$; that is, there exists a polynomial $P_\iota(t)$ with $\deg (P_\iota) \leq k$ such that
\[
\partial^\iota f(t) = P_\iota(t) \cdot (1 + |t|^2)^{a - k}.
\]

\textbf{Inductive step ($|\gamma| = k+1$).}
Let $\gamma$ be a multi-index with $|\gamma| = k+1$. Then there exists an index $j \in \{1, \dots, d\}$ and a multi-index $\beta$ with $|\iota| = k$ such that $\gamma = \iota + \vec{e_j}$, where $\vec{e_j}$ is the $j$-th standard basis vector in $\mathbb{N}_0^d$. Hence,
\[
\partial^\gamma f = \partial_{t_j} \big( \partial^\iota f \big).
\]
By the inductive hypothesis,
\[
\partial^\iota f(t) = P_\iota(t) \cdot (1 + |t|^2)^{a - k}.
\]
Differentiating with respect to $t_j$ and applying the product rule gives
\[
\partial^\gamma f(t) = 
(\partial_{t_j} P_\iota(t)) \cdot (1 + |t|^2)^{a - k}
+ P_\iota(t) \cdot \partial_{t_j} \big[ (1 + |t|^2)^{a - k} \big].
\]
Since $\partial_{t_j}(1 + |t|^2) = 2 t_j$, we compute
\[
\partial_{t_j} \big[ (1 + |t|^2)^{a - k} \big]
= (a - k) (1 + |t|^2)^{a - k - 1} \cdot 2 t_j
= 2(a - k) t_j (1 + |t|^2)^{a - k - 1}.
\]
Substituting back,
\[
\begin{aligned}
\partial^\gamma f(t) =& \partial_{t_j} P_\iota(t) (1 + |t|^2)^{a - k}
+ 2(a - k) t_j P_\iota(t) (1 + |t|^2)^{a - k - 1}\\
=& 
\Big[ \partial_{t_j} P_\iota(t)(1 + |t|^2) + 2(a - k) t_j P_\iota(t) \Big] \cdot (1 + |t|^2)^{a - (k+1)}.
\end{aligned}
\]
Define
$ 
P_\gamma(t) := \partial_{t_j} P_\iota(t)(1 + |t|^2) + 2(a - k) t_j P_\iota(t).
 $
We now verify the degree bound. Since $\deg (P_\iota) \leq k$, we have:
\[
\begin{cases}
\deg(\partial_{t_j} P_\iota(t) ) \leq k - 1,\\
\deg\big((\partial_{t_j} P_\iota(t) )(1 + |t|^2)\big) \leq (k - 1) + 2 = k + 1,\\
\deg(t_j P_\iota (t)) \leq 1 + k = k + 1.\\
\end{cases}
\]
Thus, $P_\gamma(t)$ is a polynomial with $\deg (P_\gamma) \leq k + 1$, and
\[
\partial^\gamma f(t) = P_\gamma(t) \cdot \langle t \rangle^{2a -2 |\gamma|}.
\]
This completes the inductive step.  \hfill\qed 
\section{Proof of lemma \ref{Pure-Quadratic-Term}}\label{prf-Pure-Quadratic-Term}
Write $ S_\omega[f]= \left\langle S_\omega, f\right\rangle $ for short. Taking absolute values of the integrand and applying the Fubini-Tonelli theorem, we have
\begin{equation*}
\begin{aligned}
&\mathbb{E}_\omega \int \int
\left|  (S_\omega\ast P_{y})(x)
\overline{(S_\omega\ast P_{y})(t)}
\varphi(x) \overline{\varphi(t)}
\right| dx dt\\
=&
\int \int
\mathbb{E}_\omega \left| S_\omega[P_{y+\underline{x}}] S_\omega[P_{y+\underline{t}}]\right|
|\varphi(x)| |\varphi(t)| dx dt,\\
\leq&
\int \int
\left(
\mathbb{E} \left| S_\omega[P_{y+\underline{x}}]\right|^2
\mathbb{E}
\left|
S_\omega[P_{y+\underline{x}}]
\right|^2
\right)^{\frac{1}{2}}
|\varphi(x)| |\varphi(t)| dx dt,\\
\simeq&
\int \int
\|P_{y+\underline{x}}\|_{L^2}
\|P_{y+\underline{t}}\|_{L^2}
|\varphi(x)| |\varphi(t)| dx dt,~
(\text{isometric property})\\
\simeq&
\frac{1}{y}\|\varphi\|_{L^1}^2.
\end{aligned}
\end{equation*}
 \hfill\qed
\section{Proof of lemma \ref{Cross Term}}\label{prf-Cross Term}
Taking absolute values of the integrand and applying the Fubini-Tonelli theorem, we obtain
\begin{equation*}
\begin{aligned}
&\mathbb{E}_\omega \int
\left|
(S_\omega\ast P_{y})(x)
\overline{\left\langle S_\omega,  \varphi\right\rangle}
\varphi(x)
\right| dx \\
=&
\int
\left(
\mathbb{E}_\omega
\left|S_\omega[P_{y+\underline{x}}]
S_\omega[\varphi]\right|
\right)
|\varphi(x)| dx ,\\
\leq&
\int
\left(
\mathbb{E} \left| S[P_{y+\underline{x}}]\right|^2
\mathbb{E}
\left|
S[\varphi]\right|^2
\right)^{\frac{1}{2}}
|\varphi(x)|  dx ,\\
\simeq&
\|\varphi\|_{L^2}
\int
\|P_{y+\underline{x}}\|_{L^2}
|\varphi(x)| dx ,~
(\text{isometric property})\\
\simeq&
\frac{1}{\sqrt{y}}\|\varphi\|_{L^1}\|\varphi\|_{L^2}.
\end{aligned}
\end{equation*}
 \hfill\qed
\section{Proof of theorem \ref{Thm-weakly converge}}\label{prf-Thm-weakly converge}
Note that
\[
\begin{aligned}
&\mathbb{E}\left|
\left\langle S\ast P_y ,  \varphi\right\rangle -
\langle {S, \varphi}\rangle
\right|^2\\
=&
\mathbb{E}\left|
\left\langle S\ast P_y,  \varphi\right\rangle
\right|^2
+
\mathbb{E}| \langle {S, \varphi}\rangle|^2
\\
&-\mathbb{E}
\left(
\left\langle S\ast P_y ,  \varphi\right\rangle \overline{ \langle {S, \varphi}\rangle}
\right)
-
\mathbb{E}
\left(
 \overline{\left\langle S\ast P_y ,  \varphi\right\rangle }
 \langle {S, \varphi}\rangle
\right)
.
\end{aligned}
\]
We estimate the pure quadratic terms and the cross term separately. Estimate of the pure quadratic terms:
For $0<y<1$, we have
\begin{equation}
\begin{aligned}
\mathbb{E}\left|
\left\langle S\ast P_y ,  \varphi\right\rangle
\right|^2
=&
\mathbb{E}_\omega
\int S_\omega[P_{y+\underline{x}}] \varphi(x) dx
\int \overline{S_\omega[P_{y+\underline{t}}] \varphi(t)} dt
,\\
=&
\int  \int \mathbb{E}
\left[
S_\omega[P_{y+\underline{x}}]
\overline{S_\omega[P_{y+\underline{t}}]}
\right]  \varphi(x) dx
\overline{\varphi(t)}  dt,\\
=&
\int \int P_{2y}(t-x)  \varphi(x) dx
\overline{\varphi(t)} dt,\\
=&  \left\langle P_{2y}\ast \varphi,  \varphi\right\rangle_{L^2},
\end{aligned}
\end{equation}
where the interchange of integrals is justified by Lemma \ref{Pure-Quadratic-Term}. Estimate of the cross term:
\begin{equation}
\begin{aligned}
\mathbb{E}_\omega
\left(
\left\langle S_\omega\ast P_y ,  \varphi\right\rangle \overline{S_\omega[\varphi]}
\right)
=&
\mathbb{E}_\omega
\left(
\int
S_\omega[P_{y+\underline{x}}]\varphi(x) \overline{S_\omega[\varphi]} dx
\right)
\\
=&
\int \mathbb{E}
\left(
S[P_{y+\underline{x}}] \overline{S[\varphi]}
\right)
\varphi(x) dx,\\
=&
\int
(P_{y}\ast \varphi) (x)
\varphi(x) dx,\\
=&
\left\langle P_{y}\ast \varphi,  \varphi\right\rangle_{L^2},
\end{aligned}
\end{equation}
where the interchange of integrals is justified by Lemma  \ref{Cross Term}.  Combining the above and using the identity approximation property of the Poisson kernel, we conclude
\begin{equation}
\begin{aligned}
&\mathbb{E}\left|
\left\langle {S\ast P_y} ,  \varphi\right\rangle -
\langle {S, \varphi}\rangle
\right|^2\\
=&
\left\langle P_{2y}\ast \varphi,  \varphi\right\rangle -
 \left\langle P_{y}\ast \varphi,  \varphi\right\rangle
 -\overline{\left\langle P_{y}\ast \varphi,  \varphi\right\rangle}
+
\|\varphi\|^2_{L^2} ,\\
 \xrightarrow[y\rightarrow 0_+]{~}
 &
\|\varphi\|^2_{L^2}-2\|\varphi\|^2_{L^2}+\|\varphi\|^2_{L^2}=0.
\end{aligned}
\end{equation}
 \hfill\qed
\section{Proof of lemma  \ref{lem-Poiss-weight-L2}}\label{prf-lem-Poiss-weight-L2}
For $r<\frac{d+2}{2}$, we have 
\begin{equation}
\begin{aligned}
\|\widehat{P_y}\|^2_{H^{r}(\R^d)}\simeq&
y^2
\int_{\R^d}
\frac{(1+|x|^2)^r}{(y^2+|x|^2)^{1+d}} dx,\\
=&
y^2\left(
\int_{|x |\leq y}
+
\int_{|x |>y}
\right) .\\
\end{aligned}
\end{equation}
For the estimate of the integral $\int_{|x |\leq y} $, we have
\begin{equation}
\begin{aligned}
\int_{\{x\in\R^d:~|x |\leq y\}}
\frac{(1+|x|^2)^r}{(y^2+|x|^2)^{1+d}} dx
\leq&
\frac{(1+y^2)^r }{y^{2(1+d)}}\int_{|x |\leq y}dx
\lesssim&
\frac{(1+y^2)^r }{y^{2+d}}.
\end{aligned}
\end{equation}
For the estimate of the integral $\int_{|x |> y} $, we consider two cases. For $y>1$, we have
\begin{equation}
\begin{aligned}
\int_{\{x\in\R^d:~|x |> y\}}
\frac{(1+|x|^2)^r}{(y^2+|x|^2)^{1+d}} dx
\leq&
\int_{|x |> y}
\frac{(1+|x|^2)^r}{(1+|x|^2)^{1+d}} dx,\\
=&
\int_{|x |> y}
\frac{1}{(1+|x|^2)^{(1+d)-r}} dx,\\
\lesssim&
\int_{ y}^\infty
\rho^{2r-2(1+d)} \rho^{d-1} d\rho,\\
\simeq&
\frac{1}{y^{d+2-2r}},
\end{aligned}
\end{equation}
with $r<\frac{d+2}{2}$. For $0<y<1$, it follows that
\begin{equation}
\begin{aligned}
\int_{\{x\in\R^d:  |x |> y\}}
\frac{(1+|x|^2)^r}{(y^2+|x|^2)^{1+d}} dx
=&
\int_{y<|x |\leq 1}
+
\int_{|x |> 1},\\
\leq&
2^r
\int_{ y<|x |\leq 1}
\frac{1}{|x |^{2(1+d)}}
dx
+
2^r
\int_{ |x |> 1}
\frac{|x |^{2r}}{|x |^{2(1+d)}}
dx, \\
\leq& C_r
\left({\frac{1}{y^{2+d}}}+1 \right) ,
\end{aligned}
\end{equation}
with $r<\frac{d+2}{2}$. In summary,
\begin{equation}
\|\widehat{P_y}\|^2_{H^{r}(\R^1)}
\leq C_r 
\begin{cases} y^2
\frac{1}{y^{d+2-2r}}
=\mathcal{O}\left(\frac{1}{y^{d-2r}}\right),~&y>1,\\
y^2
\left( 
{\frac{1}{y^{2+d}}}+1
\right) 
=\mathcal{O}\left(\frac{1}{y^{d}}\right),~&0<y<1,\\
\end{cases}
\end{equation}
with $r<\frac{d+2}{2}$. \hfill\qed\\
\section{Proof of corollary  \ref{lem-cutoff-sobolev}}\label{prf-lem-cutoff-sobolev}
For the Sobolev $H^s$ norm of cut-off function  $\eta_{-1}(\frac{\cdot}{\lambda})$ with any $s>0$, we have
\begin{equation}
\begin{aligned}
\left\|
\eta_{-1}
\left(\frac{\cdot}{\lambda} \right)
\right\|^2_{H^{s}(\R^d)}
=&
\lambda^d
\left\|
\eta_{-1}^{\vee}
\right\|^2_{L^2\left(\R^d,~ \langle \frac{\cdot}{\lambda}\rangle^s\right)},\\
=&
\lambda^d
\int_{\R^d}
|\eta_{-1}^{\vee}(x )|^2
\left(1+\left|\frac{x }{\lambda} \right|^2 \right)^s
dx,\\
\leq&
\lambda^d
\int_{\R^d}
|\eta_{-1}^{\vee}( x )|^2
\left(1+ |x |^2 \right)^s
dx,  ~\lambda=1+y>1 \\
\leq&
C_s \lambda^d
\int_{\R^d}
\left |
\frac{1}{(1+|x|)^{s+d}}
\right|^2
\left(1+ |x |^2 \right)^s
dx,~\eta_{-1}\in \mathscr{S}(\R^d)\\
\lesssim&
C_s \lambda^d.
\end{aligned}
\end{equation}
By the Sobolev estimate of the product of functions in property \ref{SobolevEstimate-product},  there holds
\begin{equation}
\begin{aligned}
\left\| m_{-1}  \left(\frac{\cdot}{\lambda} \right)
\right\|^2_{H^{r}}
\lesssim&
\left\|
\eta_{-1}
\left(\frac{\cdot}{\lambda} \right)
\right\|^2_{H^{s}(\R^1)}
\left\|
\widehat{P_y} \left(\frac{\cdot}{\lambda} \right)
\right\|^2_{H^{r}(\R^1)}
\\
\lesssim&
C_s C_r  \lambda^d  \left(1+y^{-d} \right),
\\
\lesssim&
C_{s,r}  \left(1+y^{d}+y^{-d}\right),\\
\simeq&
C_{s,r}  \left(y^{d}+y^{-d}\right)
,
\end{aligned}
\end{equation}
for $r\in(0,\frac{d+2}{2} )$ if we choose  $s>\frac{d+2}{2} >r$. Then for $\tau>\frac{d+2}{2} >\alpha>0$ , we have
\begin{equation}\label{SblvEstmt-product}
\begin{aligned}
\left\|
m_{-1}
\left(\frac{\cdot}{\lambda}\right)  ~g
\right\|
^2_{H^{\alpha}}
\lesssim&
\left\|  m_{-1}  \left(\frac{\cdot}{\lambda}\right) \right\|^2_{H^{r}}
\left\| g
\right\|^2_{H^{\tau}},~
\forall~g\in H^{\tau}(\R^d)\\
\lesssim&
C_{s,r}
\left(y^{d}+y^{-d} \right)
\left\| 
g
\right\|^2_{H^{\tau}}.  
\end{aligned}
\end{equation}
 \hfill\qed
\section{Proof of proposition \ref{prop-high-dertv-exp}}\label{prf-prop-high-dertv-exp}
We proceed by induction on the order of the multi-index $|\alpha|$.

\textbf{Base case ($|\alpha| = 1$).} In this case,    $\partial^\alpha = \frac{\partial}{\partial \xi_j}$ for some $j \in \{1, \dots, d\}$. Thus,
\[ 
\left| 
\frac{\partial}{\partial \xi_j} e^{-|\xi|}
\right| 
=
\left| \frac{\xi_j} {|\xi|}\right|
e^{-|\xi|}
\leq e^{-|\xi|}.
\]

\textbf{Inductive hypothesis ($|\alpha| \leq k$).}
Assume that the estimate holds for all multi-indices $\alpha$ with $|\alpha|\leq {k}$, i.e.,
\[ 
|\partial^\alpha e^{-|\xi|}|
\leq
C(\alpha, d)
e^{-|\xi|}. 
\]

\textbf{Inductive step ($|\varsigma| = k+1$).} 
Consider a  multi-index $\varsigma$ with $|\varsigma| = k + 1$. We can write  $\partial^\varsigma = \frac{\partial}{\partial \xi_j} \circ \partial^\alpha$ for some index $j$ and a multi-index $\alpha$ with $|\alpha| = k$.  Applying the Leibniz rule \eqref{LeibnizRule}, we obtain:
\[ 
\begin{aligned}
\partial^{\varsigma} e^{-|\xi|} 
&= \partial^\alpha 
\left( \frac{\partial}{\partial \xi_j}  e^{-|\xi|} \right) \\
&= \partial^\alpha 
\left( 
\frac{\xi_j} {|\xi|}
e^{-|\xi|}
\right) \\
&= 
\sum_{\beta:~|\beta|\leq|\alpha|} \binom{\alpha}{\beta}
\partial^{\alpha-\beta} \left(  \frac{\xi_j} {|\xi|}\right) 
\partial^\beta ( e^{-|\xi|}). 
\end{aligned}
\]
It then follows that:
\[ 
\begin{aligned}
\left| 
\partial^\varsigma e^{-|\xi|}
\right| 
&\lesssim  
\sum_{|\beta|\leq|\alpha|}
\left| 
\partial^{\alpha-\beta} \left(  \frac{\xi_j} {|\xi|}\right) 
\right| 
\underbrace{\left| 
	\partial^\beta  e^{-|\xi|}
	\right|}_{ 
\leq
C(\beta, d)
e^{-|\xi|},~  |\beta| \leq  k
	} \\
&\leq 
C(\alpha, d)
e^{-|\xi|}
\sum_{|\beta|\leq|\alpha|}
\left| 
\partial^{\alpha-\beta} \left(  \frac{\xi_j} {|\xi|}\right) 
\right| \\
&=
C(\alpha, d)
e^{-|\xi|}
\sum_{|\beta|\leq|\alpha|}
\left| 
\partial^{\alpha-\beta}
\frac{\partial}{\partial \xi_j}  {|\xi|}
\right| \\
&\lesssim
C(\varsigma, d)
e^{-|\xi|}.
\end{aligned}
\] 
This estimate relies on the fact that  $\{\xi \mapsto |\xi|\}$ belongs to the symbol class $S^1(\mathbb{R}^d)$ for $\xi \neq 0$ \cite{Alinhac2007_Pseudo-differential,Stein1993}. Consequently, it satisfies 
\[
\forall \beta \in \mathbb{N}^d, \quad \left| \partial_\xi^\beta |\xi| \right|
\leq C_{\beta} (1 + |\xi|)^{1 - |\beta|}
\leq C_{\beta} (1 + M)^{1 - |\beta|}
\lesssim C_{\beta},
\]
whenever $|\xi| \geq M > 0$. This concludes the proof. \hfill $\square$
\section{Proof of lemma   \ref{lem-high-freq-multpl}}\label{prf-lem-high-freq-multpl}
Set $K \in \mathbb{N}$ with $K>[\beta]+1$, then it follows that
\begin{equation}
\begin{aligned}
\left\|
m_j
\right\|
^2_{H^{\beta}(\R^d)}
\leq&
\left\|
m_j
\right\|
^2_{H^K(\R^d)},\\
\simeq&
\sum_{k:~|k|\leq{K}}
\left\|
\partial^k(\eta_{j}  \widehat{P_y} )
\right\|
^2_{L^2(\R^d)},\\
=&
\sum_{k:~|k|\leq{K}}
\left\|
\sum_{l:~l\leq{k}}
\binom{k}{l}
\partial^l \widehat{P_y} \partial^{k-l} \eta_j
\right\|^2_{L^2}, \\
\lesssim&
\sum_{|k|\leq{K}} \sum_{l\leq{k}}
\left\|
\partial^l \widehat{P_y} \partial^{k-l} \eta_j
\right\|^2_{L^2}.
\end{aligned}
\end{equation}
By Corollary \ref{Cor-Estmt-exp-funct}, the estimates for the higher-order  derivatives of  $\widehat{P_y}$ in the high-frequency region are given by
\begin{equation}
\begin{aligned}
\left|
\partial^l \widehat{P_y}(\xi)
\right|
\lesssim&
 y^{|l|} e^{-2\pi y|\xi|},\quad \forall~\xi \in \operatorname{Spt} \eta_j=\{\xi:~|\xi|\simeq 2^j\},
\end{aligned}
\end{equation}
and the higher-order derivative estimates for the cutoff function $\eta_j$ satisfy
\begin{equation}
\left| \partial^{l} \eta_j(\xi) \right|
\leq  \frac{C_{l}}{2^{j|l|}},\quad \operatorname{Spt} \eta_j=\{\xi:~|\xi|\simeq 2^j\},
\end{equation}
for any $l \in \mathbb{N}$ \cite{Chemin2011}.
Thus we have
\begin{equation}
\begin{aligned}
\left\|
m_j
\right\|
^2_{H^{K}(\R^d)}
\lesssim&
\sum_{|k|\leq{K}} \sum_{l\leq{k}}
\left\|\partial^{k-l} \eta_j\right\|_{L^\infty(\R^d)}^2
\left\|
\partial^l \widehat{P_y} 
\right\|^2_{L^2(\operatorname{Spt} \{\eta_j\}   )}
\\
\lesssim&
\sum_{|k|\leq{K}}
\frac{1}{2^{2j|k|}}
\sum_{l\leq{k}}
(2^j  y)^{2|l|}
\underbrace{\int_{2^{j-1}\leq|\xi|\leq2^{j+1}}
	e^{-4\pi y|\xi|}d\xi}_
{\simeq\int_{2^{j-1}}^{2^{j+1}}
	e^{-4\pi y\rho} \rho^{d-1}d\rho
	},\\
\lesssim&
\sum_{|k|\leq{K}}
\frac{1}{2^{2j|k|}}
\sum_{l\leq{k}}
(2^j  y)^{2|l|}
\left( 
2^{(j+1)(d-1)}\int_{2^{j-1}}^{2^{j+1}}
e^{-4\pi y\rho}d\rho
\right) 
,\\
\lesssim&
\frac{1}{y}
\sum_{|k|\leq{K}}
\frac{2^{j(d-1)}}{2^{2j|k|}}
\sum_{l\leq{k}}
(2^j y)^{2|l|}
e^{-2\pi y 2^j}
,
\end{aligned}
\end{equation}
with  $K>[\beta]+1$. Note that  $H^\beta(\R^d)$ is a Banach algebra for any $\beta>\frac{d}{2}$ \cite{Chemin2011}, then
\begin{equation}
\left\|
m_j g
\right\|
^2_{H^{\beta}}
\lesssim
\left\|
m_j
\right\|^2_{H^{K}}
\left\| g
\right\|^2_{H^{\beta}},~
\forall~g\in H^\beta(\R^d).
\end{equation}
 \hfill\qed
\section*{Acknowledgements}
The authors would like to thank H. G.  Feichtinger for his careful reading of the manuscript and constructive comments.  
The research is supported by National Science and Technology Major Project (2027ZD01999504), the Major Project of Guangzhou National Laboratory (Grant No. GZNL2024A01004), the Science and Technology Development Fund of Macau SAR (Grant No. 0015/2025/AFJ), and the Science and Technology Development Fund of Macau SAR (Grant No. 0002/2024/RDP). 

\newpage

\end{CJK*}
\end{sloppypar}

\end{document}